\documentclass[11pt]{article}
\usepackage{epsfig,makeidx,latexsym,amsmath,amssymb,amsfonts,color,srcltx,multirow,url}
\usepackage{tikz,graphics,wrapfig,times}
\usepackage[colorlinks=true,breaklinks=true,bookmarks=true,urlcolor=blue,
     citecolor=blue,linkcolor=blue,bookmarksopen=false,draft=false]{hyperref}
\usepackage{fullpage}
\makeatletter
\def\imod#1{\allowbreak\mkern10mu({\operator@font mod}\,\,#1)}
\makeatother
\newtheorem{proposition}{Proposition}[section]
\newtheorem{corollary}{Corollary}[section]
\newtheorem{theorem}{Theorem}[section]
\newtheorem{example}{Example}[section]
\newtheorem{remark}{Remark}[section]

\newtheorem{definition}{Definition}[section]

\newtheorem{lemma}{Lemma}[section]

\def\three?{7}
\def\four?{4}
\def\ten?{10}
\def\myp{$\!\!\!~^)$}

\def\norm2to2{{\|\cdot\|_{2,2}}}

\def\Opt{\hbox{\rm Opt}}

\def\Conv{\hbox{\rm  conv}}
\def\conv{\hbox{\rm  conv}}
\def\clconv{\overline{\conv}}

\def\cone{\hbox{\rm  cone}}

\def\Id{{\hbox{\rm Id}}}
\def\Tr{{\mathop{\hbox{\rm  Tr}}}}

\def\cB{{\cal B}}
\def\cC{{\cal C}}

\def\cF{{\cal F}}

\def\cK{{\cal K}}
\def\cL{{\cal L}}

\def\cS{{\cal S}}

\def\cV{{\cal V}}

\def\cZ{{\cal Z}}

\def\Ext{{\mathop{\hbox{\rm  Ext}}}}

\def\tilde{\widetilde}

\def\Ker{{\hbox{\rm  Ker}}}
\def\ker{{\hbox{\rm  Ker}}}
\def\Im{{\hbox{\rm  Im}}}

\def\bR{{\mathbb R}}

\def\Z{{\mathbb Z}}

\def\qed{\ \hfill$\square$\par\smallskip}

\def\intt{\hbox{\rm int}}
\def\intt{\hbox{\rm  int}}

\def\solS{\cS(A,\cK,\cB)}
\def\viPi{\Pi(A,\cK,\cB)}
\def\coneC{C(A,\cK,\cB)}
\def\coneCm{C_m(A,\cK,\cB)}

\def\coneCs{C_s(A,\cK,\cB)}

\def\mybeq{\begin{equation}\begin{array}{c}}
\def\myeeq{\end{array}\end{equation}}

\def\mypict3{\epsfxsize=220pt\epsfysize=80pt\epsfbox}

\def\bR{{\mathbb R}}

\newcommand{\ese}{\end{eqnarray*}}
\newcommand{\bse}{\begin{eqnarray*}}

\newcommand{\epr}{\hfill$\diamondsuit$\vspace{0.25cm}\par}

\newcommand{\myfootnote}[1]{\footnote{\myp#1}$~^)$}
\newcommand{\myproof}[1]{\textbf{Proof.}\quad #1 \qed}

\definecolor{DarkBlue}{rgb}{0.1,0.1,0.6}
\definecolor{Red}{rgb}{0.9,0.0,0.1}

\begin{document}
\title{On Minimal Valid Inequalities for Mixed Integer Conic Programs}

\author{Fatma K{\i}l{\i}n\c{c}-Karzan\thanks{Carnegie Mellon University, Pittsburgh, Pennsylvania 15213, USA, {\tt fkilinc@andrew.cmu.edu}}}
 
\date{March 31, 2015}

\maketitle
\begin{abstract}
We study disjunctive conic sets involving a general regular (closed, convex, full dimensional, and pointed) cone $\cK$ such as the nonnegative orthant, the Lorentz cone or the positive semidefinite cone. In a unified framework, we introduce $\cK$-minimal inequalities and show that under mild assumptions, these inequalities together with the trivial cone-implied inequalities are sufficient to describe the convex hull.
We study the properties of $\cK$-minimal inequalities by establishing algebraic necessary conditions for an inequality to be $\cK$-minimal. This characterization leads to a broader algebraically defined class of $\cK$-sublinear inequalities. We establish a close connection between $\cK$-sublinear inequalities and the support functions of sets with a particular  structure. This connection results in practical ways of showing that a given inequality is $\cK$-sublinear and $\cK$-minimal. 

Our framework generalizes some of the results from the mixed integer linear case. It is well known that the minimal inequalities for mixed integer linear programs are generated by sublinear (positively homogeneous, subadditive and convex) functions that are also piecewise linear. This result is easily recovered by our analysis. Whenever possible we highlight the connections to the existing literature. However, our study unveils that such a cut generating function view treating the data associated with each individual variable independently is not possible in the case of general cones other than nonnegative orthant, even when the cone involved is the Lorentz cone.
\end{abstract}


\section{Introduction}

A \emph{Mixed Integer Conic Program} (MICP) is an optimization program of the form
$$
\Opt=\inf_{x\in E}\left\{\langle c, x\rangle:
Ax=b,~ x\in\cK,~ x\in\cZ \right\} \eqno{(MICP)}
$$
where $\cK$ is a \emph{regular} (full-dimensional, closed, convex and pointed) cone in a finite dimensional Euclidean space $E$ with an inner product $\langle\cdot,\cdot\rangle$, $c\in E$ is the objective vector, $b\in \bR^m$ is the  right hand side vector, $A:E\rightarrow \bR^m$ is a  linear map, and $\cZ$ is a set imposing certain structural restrictions on the variables $x$. Examples of regular cones include the nonnegative orthant $\bR^n_+:=\{x\in\bR^n:~x_i\geq 0 ~\forall i=1,\ldots,n\}$, the Lorentz cone $\mathcal{L}^n:=\{x\in\bR^n:~x_n\geq\sqrt{x_1^2+\ldots+x_{n-1}^2}~\}$, and the positive semidefinite cone $\mathcal{S}^n_+:=\{x\in\bR^{n\times n}:~ a^Tx\ \! a\geq 0 ~\forall a\in\bR^n, ~x=x^T \}$ and their direct products. When $E=\bR^n$, the most common form of structural restrictions is integrality $x_i\in\Z \mbox{ for all } i\in I $ where $I \subset \{1,\ldots, n\}$ is the index set of integer variables. We assume that all of the data involved with MICP, i.e., $c,b,A$ is rational. 

Mixed Integer Linear Programs (MILPs) arise as a special case of MICP where $\cK$ is the nonnegative orthant.  
Conic constraints include various specific convex constraints such as linear, convex quadratic, eigenvalue, etc., and hence, offer significant representation power over linear constraints (see \cite{BenTal_Nemirovski_01} for a detailed introduction to conic programming and its applications in various domains). 
Allowing discrete decisions in addition to the conic constraints further enhances the representation power of MICPs. 
While MILPs offer an incredible representation power, various optimization problems involving risk constraints and discrete decisions give rise to MICPs. Robust optimization and stochastic programming paradigms, or more broadly decision making under uncertainty domain,  encompasses many examples of MICPs such as portfolio optimization with fixed transaction costs  in finance \cite{Drewes_Pokutta_10,Lobo_Fazel_Boyd_07}, and stochastic joint location-inventory models \cite{Atamturk_Berenguer_Shen_12}. 
Moreover, the most powerful relaxations to many combinatorial optimization problems are based on conic (in particular semidefinite) relaxations (see 
\cite{Goemans_97} for a survey on this topic).  Reintroducing the integer variables back into these relaxations yields exact mixed integer conic programming formulations of these problems with tighter continuous relaxations. 
Besides, MILPs have been heavily exploited for approximating non-convex nonlinear optimization problems arising in several important applications in many diverse fields. For a wide range of these problems, MICPs offer tighter relaxations and thus potentially a better overall algorithmic performance. 
Therefore, MICPs have gained considerable interest.

{
The literature on solving MICPs is growing rapidly. On one hand, clearly, any method for general nonlinear integer programming applies to MICPs as well. 
A significant body of work has extended known techniques from MILPs to nonlinear integer programs. These include the Reformulation Linearization Technique (see \cite{Sherali_Tuncbilek_95,Sherali_Adams_98} and references therein), 
Lift-and-Project  and Disjunctive Programming methods \cite{Balas_79,Balas_Ceria_Cornuejols_93,Burer_Saxena_12,Modaresi_Kilinc_Vielma_13,Saxena_Bonami_Lee_08,Sherali_Shetti_80,Stubbs_Mehrotra_99,Stubbs_Mehrotra_02}, and the lattice-free set paradigm \cite{Bienstock_Michalka_13}. In addition to these, several papers  \cite{Kojima_Tuncel_00,Lasserre_01,Saxena_Bonami_Lee_11,Saxena_Bonami_Lee_10} introduce hierarchies of convex (semidefinite programming) relaxations in higher dimensional spaces. These relaxations quickly become impractical due their exponentially growing sizes and the difficulty of projecting them onto the original space of variables. 
Another stream of research \cite{Abhishek_etal_06,Bonami_etal_08,Drewes_09,Qualizza_Belotti_Margot_12,Tawarmalani_Sahinidis_02,Tawarmalani_Sahinidis_04,Tawarmalani_Sahinidis_05,Vielma_Ahmed_Nemhauser_08} is on the development of linear outer approximation based branch-and-bound algorithms for nonlinear integer programming. 
While they have  
the advantage of fast and easy to solve relaxations, 
the bounds from these approximations may not be as strong as desired. Moreover, adding too many inequalities that are similar to each other  may lead to numerical instability.

Exploiting the conic structure when present, as opposed to general convexity,  paves the way for developing algorithms with much better performance. Particularly in the case of MILPs, this has led to very successful results.  Despite the lack of effective warm-start techniques, efficient interior point methods exist for $\cK=\mathcal{L}^n$ or $\cK=\mathcal{S}^n_+$ \cite{BenTal_Nemirovski_01}. Therefore, supplying the branch-and-bound tree with the corresponding continuous conic relaxation at the nodes and deriving cutting planes to strengthen these relaxations have gained considerable interest recently. 
In this vein, \c{C}ezik and Iyengar \cite{Cezik_Iyengar_05} developed valid inequalities for MICPs with general regular cones by extending Chvatal-Gomory (C-G) integer rounding cuts \cite{Nemhauser_Wolsey_88}. 
In a recent and fast growing literature, several authors \cite{Andersen_Jensen_13,Atamturk_Narayanan_10,Atamturk_Narayanan_11,Belotti_Goez_Polik_Terlaky_12,Bienstock_Michalka_13,BKK14,Dadush_Dey_Vielma_11_2,Drewes_09,Drewes_Pokutta_10,KKY14_ipco,KKY14,Modaresi_Kilinc_Vielma_13,Vielma_Ahmed_Nemhauser_08,YC14} study MICPs involving Lorentz cones, $\cK=\mathcal{L}^n$, and suggest valid inequalities.  
}

This growing demand for solving MICPs has led many commercial software packages such as CPLEX \cite{cplex}, Gurobi \cite{Gurobi}, and MOSEK \cite{Mosek} to recently expand their features and include technology to solve MICPs. Nevertheless, the theory and algorithms for solving MICPs are still in their infancy \cite{Atamturk_Narayanan_10}.
Currently, the most promising approaches to solve MICPs are based on the extension of cutting plane techniques \cite{Atamturk_Narayanan_10,Atamturk_Narayanan_11,Billionnet_Sourour_Marie-Christine_09,Bonami_11,Buchheim_Caprara_Lodi_12,Cezik_Iyengar_05,Drewes_09,Kilinc_Linderoth_Luedtke_10,Stubbs_Mehrotra_99} in combination with conic continuous relaxations and branch-and-bound algorithms. While numerical performance of these techniques is still under investigation, evidence from MILPs indicates that adding a small yet essential set of strong cutting planes  is key to the success of such a procedure. Yet, except very specific and simple cases, the strength (redundancy, domination, etc.) of the corresponding valid inequalities has not been evaluated in the case of MICPs. 
This is in sharp contrast to the MILP case, where the related questions have been studied extensively. In particular, the feasible region of an MILP with rational data is a polyhedron and the facial structure of a polyhedron (its faces, and facets) is very well understood.  Various ways of proving whether or not a given linear inequality is necessary in the description of the convex hull of the feasible set of MILP, e.g., a facet, are well established \cite{Nemhauser_Wolsey_88}. In addition to this, a new framework to establish minimality and extremality of valid inequalities for certain generic infinite relaxations of MILPs (see \cite{Conforti_Cornuejols_Zambelli_11} and references therein) as well as their relations to facets in certain simplified settings \cite{Cornuejols_Margot_09} are developing rapidly. 
Thus far, results in this vein are lacking in the MICP context. Consequently,   establishing a theoretical framework to measure the necessity and strength of cutting planes in the MICP context remains a natural and important question. Our goal in this paper is to address this question.

In this paper we study the closed convex hull of a \emph{disjunctive conic set}, that is, the union of finitely or infinitely many conic sets in the original space of variables. We are mainly motivated by the fact that most cutting planes used in MILP can be viewed in the context of  disjunctive programming 
and such a general \emph{disjunctive conic programming} framework  encompasses MICPs. 
Our approach is based on identifying an appropriate dominance concept among valid linear inequalities  and then extending the minimality definition from the MILP context to the disjunctive conic framework. When the underlying cone $\cK$ is taken as the nonnegative orthant, the MILP counterparts of our results and further developments for MILPs were studied extensively in the literature. 
Despite this extensive literature for $\cK=\bR^n_+$, to the best of our knowledge there is no literature on this topic in the general conic case with an arbitrary regular cone $\cK$. 
We contribute to the literature by introducing minimal inequalities for disjunctive conic sets and performing a systemic study of their properties in a unified manner for all regular cones $\cK$. We establish the  sufficiency of minimal inequalities along with necessary conditions, sufficient conditions as well as practical tools for testing whether or not a given inequality is minimal. 

Our derivations are based on a finite dimensional problem instance. This is in contrast to much of the literature on minimal inequalities for MILPs initiated by \cite{Gomory_Johnson_72_1, Gomory_Johnson_72_2, Johnson_74}. 
In a practical cutting plane procedure for solving MILPs and/or MICPs, one is indeed faced with a  problem in a finite dimensional space. Thus, we believe that this is not a limitation but rather a contribution to the corresponding MILP literature. Besides, to the best of our knowledge, the extensions of other well-known regular cones such as $\cL^n$ and $\cS^n_+$ to the infinite dimensional spaces are not well defined. Hence, an infinite relaxation seems to be more meaningful when the associated cone is the nonnegative orthant.
Therefore, our study does not rely on and differ substantially from the majority of previous literature in the MILP context that relies on infinite relaxations. 
{ 
Furthermore, we note that a conic view with a regular polyhedral cone can be valuable in the MILP context as well.

We demonstrate that some of the results from MILP setup naturally extend to MICPs.  
In this regard, our approach ties back to  the cornerstone paper of Johnson \cite{Johnson_81} as well as the recent work of Conforti et al.\@ \cite{Conforti_Cornuejols_etal_13}. 
In particular,  when $\cK=\bR^n_+$, our results show that minimal inequalities can be properly related to support functions that generate cut coefficients\myfootnote{Informally, these are referred as  \emph{cut generating functions}. A cut generating function generates the coefficient of a variable in a cut using only information of the instance pertaining to this variable. See \cite{Conforti_Cornuejols_etal_13} and section \ref{sec:ConnectionToLatticeFreeSetsAndCGFs} for an extended discussion.}, and these functions are sublinear (subadditive, positively homogeneous, and convex) and piecewise linear. This connection in the case of $\cK=\bR^n_+$ together with the sufficiency of minimal inequalities for describing the closed convex hull of disjunctive conic sets highlights the roots of  functional strong duality results for MILPs.
For other regular cones, we show that there exist extreme inequalities, which cannot be generated from any cut generating function when we straightforwardly extend the definition of cut generating functions to MICPs. 
Whenever possible, we highlight these connections to the existing literature.
}

\subsection{Preliminaries and Notation}\label{problem}
 
Let $(E,\langle \cdot,\cdot\rangle)$ be a finite dimensional Euclidean space with inner product $\langle \cdot,\cdot\rangle$. 
Let $\cK\subset E$ be a \emph{regular} (full-dimensional, closed, convex and pointed) cone. Note that when every $\cK_i\subset E_i$ for $i=1,\ldots,k$ is a regular cone, then their direct product  $\tilde{\cK}=\cK_1\times...\times\cK_k$ is also a regular cone in the Euclidean space $\tilde{E}=E_1\times...\times E_k$ with inner product $\langle \cdot,\cdot\rangle_{\tilde{E}}$, which is the sum of the inner products $\langle \cdot,\cdot\rangle_{E_i}$. Therefore, without loss of generality we only focus on the case with a single regular cone $\cK$.

In this paper, given a linear map $A:E\rightarrow \bR^m$, a regular cone $\cK\subset E$, and a \emph{nonempty} set of  right hand side vectors $\cB\subseteq \bR^m$,  we study the following \emph{disjunctive conic set} defined by $A$, $\cK$, and $\cB$: 
\[
\solS:=\{x\in \cK: ~A x \in\cB \}.
\]
We are mainly interested in determining the properties of strong valid linear inequalities describing the closed convex hull of $\solS$. We would like to emphasize that we do not impose any structural assumptions on $ A$ and $\cB$. In particular, $A$ is an arbitrary linear map from $E$ to $\bR^m$ and $\cB$ is an arbitrary set of vectors in $\bR^m$. Note that the set $\cB$ can be finite or infinite, structured such as lattice points or completely unstructured.  
In order to avoid trial cases, we assume that $\solS\neq\cK$, in particular $\cK\not\subseteq\{x\in E:\, Ax\in\cB\}$, and $\solS\neq\emptyset$, i.e., there exists $b\in\cB$ and $x_b\in\cK$ satisfying $Ax_b=b$.

For a given set $S$, we denote its topological interior with $\intt(S)$, its closure with $\overline{S}$, and  its boundary with $\partial S=\overline{S}\setminus\intt(S)$. 
We use $\conv(S)$ to denote the convex hull of $S$, $\clconv(S)$ for its closed convex hull, and $\cone(S)$ to denote the cone generated by the set $S$. We denote the kernel of a linear map $A:E\rightarrow \bR^m$ by $\Ker(A)=\{u\in E:~Au=0\}$, and its image by $\Im(A)=\{Au:~ u\in E\}$. 
We use $A^*$ to denote the conjugate linear map \myfootnote{
When we consider the standard Euclidean space $E=\bR^n$, a linear map $A:\bR^n\rightarrow \bR^m$ is just an $m\times n$ real-valued matrix, and its conjugate is given by its transpose, $A^*=A^T$. 

Also, let us consider the space of symmetric $n\times n$ matrices $E=\mathcal{S}^n$. We use $\Tr(\cdot)$ to denote the trace of a matrix, i.e., the sum of its diagonal entries. When $E=\mathcal{S}^n$, it is natural to specify a linear map $A:\mathcal{S}^n\rightarrow \bR^m$ as  a collection $\{A^1,\ldots,A^m\}$ of $m$ matrices from $\mathcal{S}^n$ such that
\[
AZ = (\Tr(ZA^1);\ldots;\Tr(ZA^m)):\mathcal{S}^n\rightarrow \bR^m.
\]
In this case, the conjugate linear map $A^*:\bR^m\rightarrow \mathcal{S}^n$ is given by
\[
A^*y=\sum_{j=1}^m y_jA^j, ~~ y=(y_1;\ldots;y_m)\in\bR^m. 
\]
}   
given by the identity 
\[
y^TAx = \langle A^*y, x\rangle ~~\forall (x\in E, y\in \bR^m). 
\]

We use $\langle \cdot,\cdot\rangle$ notation for inner product in Euclidean space $E$, and proceed with usual dot product notation with transpose for the inner product in $\bR^m$. We assume all vectors in $\bR^m$ are given in column form.  

For a given cone $\cK\subset E$, we let $\Ext(\cK)$ denote the set of its extreme rays, and use $\cK^*$ to denote its dual cone given by 
\[
\cK^* :=\left\{y\in E: ~\langle x, y\rangle \geq 0 ~~\forall x \in \cK \right\}.
\]
Whenever the cone $\cK$ is regular, so is $\cK^*$. 
\par
Given a regular cone $\cK$, a relation $a-b\in\cK$ (also denoted by $a\succeq_{\cK}b$) is called {\sl conic inequality} between $a$ and $b$. Such a relation indeed preserves the major properties of the usual coordinate-wise vector inequality $\geq$. We denote the strict conic inequality by $a\succ_{\cK}b$ to indicate that $a-b \in \intt(\cK)$. In the sequel, we refer to a constraint of the form $Ax-b\in\cK$ as a {\sl conic inequality constraint} or simply {\sl conic constraint} and also use $Ax\succeq_{\cK} b$ interchangeably in the same sense.
\par 
There are three important regular cones common to most MICPs, namely the nonnegative orthant $\bR^n_+$, the Lorentz cone $\mathcal{L}^n$, and the positive semidefinite cone $\mathcal{S}^n_+$. In the first two cases, the corresponding Euclidean space $E$ is just $\bR^n$ with dot product as the corresponding  inner product. In the last case, $E$ becomes the space of symmetric $n\times n$ matrices with Frobenius inner product $\langle x, y \rangle = \Tr(xy^T)$. These three regular cones are also self-dual, that is, $\cK^*=\cK$.
\par
Notation $e^i$ is used for the $i^{th}$ unit vector of $\bR^n$, and $\Id$ for the  identity map in $E$. When $E=\bR^n$, $\Id$ is just the $n\times n$ identity matrix $I_n$.
 
\subsection{Motivation and Connections to MICPs}\label{ConnectionsToMICP}

While the \emph{disjunctive conic set} $\solS$ can be of interest by itself, here we provide a few examples to highlight our naming choice and the significance of this framework. In particular, we show that these sets  naturally represent the feasible regions of MICPs as well as some natural  relaxations for them. 

We start with the following example transformation that generalizes the usual \emph{disjunctive programming} from the polyhedral (linear) case  \cite{Balas_71,Balas_72,Balas_79,Balas_10}  to the one with conic constraints.
\begin{example}\label{ex:disjunctive}
Suppose that we are given a finite collection of convex sets of the form $C_i=\{x\in\cK:\, A^i x\succeq_{\cK_i} b^i \}$ for $i\in\{1,\ldots,\ell\}$, where $\cK\subset \bR^n$ and $\cK_i\subset\bR^{m_i}$ are regular cones, $A^i$ are $m_i\times n$ matrices, and $b^i\in\bR^{m_i}$. Then  $\bigcup_{i\in \{1,\ldots,\ell\}} C_i$  can be represented in the form of $\solS$ as follows:
\begin{equation*}
\left\{x\in \bR^n:\, \underbrace{\begin{pmatrix}(A^1)^T \\ (A^2)^T \\ \vdots \\ (A^\ell)^T\end{pmatrix} x }_{:= Ax} \in \underbrace{\left\{\begin{pmatrix}\{b^1\}+\cK_1 \\  \bR^{m_2}  \\  \vdots \\ \bR^{m_\ell} \end{pmatrix} \bigcup \begin{pmatrix} \bR^{m_1}  \\    \{b^2\}+\cK_2  \\ \vdots \\ \bR^{m_\ell} \end{pmatrix} \bigcup \begin{pmatrix} \bR^{m_1}  \\ \bR^{m_2}   \\ \vdots  \\ \{b^m\}+\cK_m \end{pmatrix}\right\} }_{:=\cB},
~x\in\cK \right\}.
\end{equation*}

When $\cK=\bR^n_+$ and $\cK_i=\bR^{m_i}_+$ for all $i=1,\ldots,\ell$, then $\bigcup_{i\in \{1,\ldots,\ell\}} C_i$ is the well-known disjunctive set representing the union of polyhedra \cite{Balas_71,Balas_72,Balas_79,Balas_10}. 

Moreover, when $\cK$ is a general regular cone but $\cK_i=\bR_+$ for all $i=1,\ldots,\ell$, then the set $\solS$ models multi-term disjunctions on the cone $\cK$. 
\epr
\end{example}
In fact, the multi-term disjunction structure of Example \ref{ex:disjunctive} allows us to model removal of  any polyhedral lattice-free set such as triangle, quadrilateral or cross disjunction from a regular cone (or its cross-section) by appropriately selecting the cones $\cK_i$, the matrices $A^i$, and the vectors $b^i$. Besides, every convex set $Q\in E$ can be regarded as the cross-section of a convex cone in $E\times \bR$ given by $\cK_Q:=\cone(\{(x,1)\in E\times \bR:~ x\in Q \})$ and the hyperplane $H=\{(x,\lambda)\in E\times \bR:~\lambda=1\}$. Yet, the resulting cone $\cK$ may not be regular in general.

Our next set of examples highlight the connection of $\solS$ with the feasible sets of MICPs and their relaxations.
\begin{example}\label{ex:MICPform1}
Suppose that we are given the following conic optimization problem with integer variables  
\begin{equation}\label{eq:MICP1}
\Opt=\inf_{x\in\bR^n}\left\{c^T x:~
\tilde{A}x=b,~ x\in\cK,~ x_i\in\Z \mbox{ for all } i=1,\ldots,\ell  \right\}. 
\end{equation}
By defining   
\[
A=\left[\begin{array}{c}\tilde{A} \\ I_n  \end{array}\right],
\quad \mbox{and}\quad 
\cB=\left\{ \left( \begin{array}{c} b \\ \Z^{\ell} \\ \bR^{n-\ell} \end{array} \right) \right\},  
\] where $I_n$ is the $n\times n$ identity matrix, we can convert this problem into optimizing the same linear function over $\solS$, i.e., $\Opt=\inf_{x\in\bR^{n}}\left\{ c^T x:~ Ax\in\cB,~ x\in\cK  \right\}$. \epr
\end{example}

\begin{example}\label{ex:MICPform2}
Let us also consider another MICP of form 
\begin{equation}\label{eq:MICP2}
\Opt:=\inf_{y\in\bR^n}\left\{\tilde{c}^Ty:~
\tilde{A}y-b\in\tilde{\cK},~ y_i\in\Z \mbox{ for all } i=1,\ldots,\ell  \right\}, 
\end{equation}
where $\tilde{\cK}$ is a regular cone in the Euclidean space $E$. 
Then, by introducing new variables $y^+, y^-$, and setting 
\[
x=\left(\begin{array}{l} y^+\\  y^- \end{array}\right),
\quad \cK=\bR^n_+\times\bR^n_+, 
\quad c=\left(\begin{array}{c} \,\tilde{c}\\  -\tilde{c} \end{array}\right),
\quad 
A=\left[\begin{array}{cc}\tilde{A}\; &~ -\tilde{A} \\ I_n & ~ -I_n \end{array}\right],
\quad \mbox{and}\quad 
\cB=\left\{ \left( \begin{array}{c} b+\tilde{\cK} \\ \Z^{\ell} \\ \bR^{n-\ell} \end{array} \right) \right\},  
\] where $I_n$ is the $n\times n$ identity matrix, once again we can  precisely represent  this problem into disjunctive conic form. \epr
\end{example}

There is an important structural difference between the disjunctive conic sets  arising in Examples \ref{ex:MICPform1} and \ref{ex:MICPform2}:  The cone $\cK$ of $\solS$ in Example  \ref{ex:MICPform1} is rather general, in particular it can be any regular cone. On the other hand,  the resulting cone used in Example  \ref{ex:MICPform2} after the transformation is a very specific one, it is the nonnegative orthant. There are two important distinctions between a general regular cone and the specific case of nonnegative orthant that will appear in our discussions later on in section \ref{sec:subadditive}. These are, first, the nonnegative orthant is decomposable, i.e., it does not introduce correlations among variables,  and second, all of its extreme rays are orthogonal to each other.

\begin{example}\label{ex:MICPform3}
Let us revisit Example \ref{ex:MICPform2} and investigate the following alternative  disjunctive conic form given in a lifted space by a single additional variable, $t\in\bR$, as follows
\[
x=\left(\begin{array}{l} y\\  t \end{array}\right),
\quad \cK=\left\{(y;t)\in\bR^n\times\bR:~ \tilde{A} y - b t \in\tilde{\cK} \right\},
\]
together with 
\[ 
\quad c=\left(\begin{array}{c} \,\tilde{c}\\  0 \end{array}\right),
\quad 
A=\left[\begin{array}{cc} I_\ell & ~ 0\\ 0 & 1  \end{array}\right],
\quad \mbox{and}\quad 
\cB=\left\{ \left( \begin{array}{c} \Z^{\ell} \\ 1 \end{array} \right) \right\},  
\] where $I_\ell$ is the $\ell\times \ell$ identity matrix. The resulting optimization problem over this disjunctive conic set  is also exactly equivalent to \eqref{eq:MICP2}.  

Analogous transformations are  possible for Examples \ref{ex:disjunctive} and \ref{ex:MICPform1} as well. 
\epr
\end{example}

\begin{remark}
The transformation given in Example \ref{ex:MICPform2} may seem more attractive in comparison to that of Example \ref{ex:MICPform3} because the final disjunctive conic form $\solS$ in Example \ref{ex:MICPform2} possesses very simple conic structure $\cK=\bR^{2n}_+$. On the other hand, the  transformation used in Example \ref{ex:MICPform3} not only gets us to a disjunctive conic form with fewer additional variables but also the new cone $\cK$ encodes important structural information about the problem such as the linear map $\tilde{A}$ and the vector $b$. 

As we detail in section \ref{sec:minimal}, the cone $\cK$ plays a critical role in identifying dominance relations among valid inequalities for $\solS$. In particular, our minimality notion is explicitly based on the ordering defined by the dual cone $\cK^*$. As a result of this, any structural information encoded in $\cK$ is quite useful in identifying the properties of extremal inequalities. 
In fact, this opens up new possibilities even for the well-studied case of MILPs, which we discuss in  Remark \ref{rem:HowToDefineMinimality}.
\epr  
\end{remark}

In Examples \ref{ex:MICPform1}-\ref{ex:MICPform3}, we provide disjunctive conic sets $\solS$ to exactly represent the corresponding feasible sets of MICPs. This indicates that the explicit description of the resulting $\clconv(\solS)$ is often not easy to characterize. An alternative use of our disjunctive conic framework in the context of MICPs is to obtain and study disjunctive conic form relaxations that are practical, yet still nontrivial and useful. One possibility for obtaining such relaxations in the form of $\solS$ is to iteratively add the integrality requirements by changing $\bR$ to $\Z$ in the description of the set $\cB$ corresponding to a variable $x_i$.

Another option for developing relaxations in disjunctive conic form is based on a more practical separation problem. Suppose that in Example \ref{ex:MICPform1} we have obtained a feasible solution $\hat x$ to the continuous relaxation of MICP, yet $\hat{x}\notin\clconv(\solS)$. For this example, the following disjunctive conic set $\solS$ can be exploited to identify valid inequalities that cut off $\hat x$. Consider $d\in \Z^n$ and $r_0\in\Z$ such that $d_i=0$ for all $i=\ell+1,\ldots,n$ and $r_0<\sum_{i=1}^\ell d_i \hat{x}_i < r_0+1$. Then the split disjunction induced by $\sum_{i=1}^\ell d_i x_i \leq r_0~\vee~\sum_{i=1}^\ell d_i x_i \geq r_0+1$ is valid for the feasible set of the optimization problem \eqref{eq:MICP1}, whereas the current solution $\hat x$ violates it.  
Given such a split disjunction, the question of obtaining cuts separating  $\hat x$ is equivalent to studying $\clconv(\solS)$ where
\[
A=\left[\begin{array}{c}\tilde{A} \\ d^T   \end{array}\right],
\quad \mbox{and}\quad 
\cB=\left\{ \left( \begin{array}{c} b \\ r_0-\bR_+ \end{array} \right) \bigcup \left( \begin{array}{c} b \\ r_0+1+\bR_+ \end{array} \right) \right\}.  
\] 
In particular, the inequality description of this $\clconv(\solS)$ will contain cuts for the original MICP separating $\hat{x}$. 
The same reasoning also applies in the case of Example \ref{ex:MICPform2}, e.g., such a split disjunction in this case can be represented by defining $\solS$ with
\[
x=\left(\begin{array}{l} y^+\\  y^- \end{array}\right),
~~ \cK=\bR^n_+\times\bR^n_+, 
~~
A=\left[\begin{array}{cc}\tilde{A}\; &~ -\tilde{A} \\ d^T & ~ -d^T \end{array}\right],
~~ \mbox{and}~~
\cB=\left\{ \left( \begin{array}{c} b+\tilde{\cK} \\ r_0-\bR_+ \end{array} \right) \bigcup \left( \begin{array}{c} b+\tilde{\cK} \\ r_0+1+\bR_+ \end{array} \right) \right\}.
\]
 
We stress that  in our discussion above $\hat x$ is not restricted to be an extreme point solution. In many cases in MICPs, $\hat x$ will be obtained by solving a continuous relaxation of MICP via interior point methods. Therefore, it will not necessarily be an extreme point solution.  
Nevertheless, our framework is flexible enough as it allows us to study the separation of an arbitrary point $\hat{x}\notin\clconv(\solS)$.  In contrast,  most of the MILP literature, and almost all of the so-called cut-generating function literature, focuses on separating extreme point solutions. 
This main focus on the separation of extreme point solutions in theory and practice of MILPs is because the overwhelming choice for solving the linear programming relaxations is the simplex algorithm and it leads to extreme point solutions $\hat x$. 
In the MILP literature, by translation of the associated point $\hat x$ and the feasible set, this separation problem is often cast as separating the origin from the convex hull of a set of points.

Nonetheless, the theoretical framework of disjunctive programming in MILP does provide general techniques to separate non-extreme-point solutions in the same manner as discussed above. 
Thus, exact representations and relaxations of the above forms have been studied in a number of other contexts in the specific case of $\cK=\bR^n_+$. In particular, when we additionally assume that $\cB$ is finite, we immediately arrive at the \emph{disjunctive programming} framework of Balas \cite{Balas_79}.   
Furthermore, Johnson  \cite{Johnson_81} has studied the set  $\cS(A,\bR^n_+,\cB)$ when $\cB$ is a finite list under the name of \emph{linear programs with multiple right hand side choice}.  
In another closely related recent work, Conforti et al.\@ \cite{Conforti_Cornuejols_etal_13} study $\solS$ with $\cK=\bR^n_+$ and possibly an infinite set $\cB$ such that $\cB\neq\emptyset$, is closed and $0\notin\cB$, and demonstrate that Gomory's \emph{corner polyhedron}  \cite{Gomory_69} as well as some other problems such as \emph{linear programs with complementarity restrictions} \cite{Judice_Sherali_06} can be viewed in this framework. In contrast to \cite{Balas_79},  Johnson \cite{Johnson_81} 
studies the characterizations of minimal inequalities, and Conforti et al.\@ \cite{
Conforti_Cornuejols_etal_13} study minimal cut generating functions. 
We discuss connections of these and our study in section \ref{sec:ConnectionToLatticeFreeSetsAndCGFs}.

Finally, we emphasize that we are not making any particular assumption on $A$ and $\cB$ beyond the basic ones to avoid trivial cases such as $\solS=\emptyset$ or $\conv(\solS)=\cK$. Because $\cB$ can be completely arbitrary, the set $\solS$ offers great flexibility, which can be much beyond the relaxations/representations related to MICPs. Specifically, a good understanding of disjunctive conic sets will be particularly relevant to conic complementarity problems as well.

\subsection{Classes of Valid Inequalities and Our Goal} \label{sec:classOfInequalities}

Recall that we are interested in the closed convex hull characterization of the disjunctive conic set 
\[
\solS=\{x\in \cK: ~A x \in\cB \}.
\]
This naturally amounts to the study of valid linear inequalities for $\solS$.
Without loss of generality we assume that all of the linear \emph{valid inequalities} for $\solS$ are of the form
\[
\langle \mu,  x \rangle \geq \eta_0,
\]
where $\mu\in E$ and $\eta_0\in\bR$. 
We  denote the resulting  inequality with $(\mu;\eta_0)$ for short hand notation.  
For any $\mu\in E$, we define 
\begin{equation}\label{mu0}
\vartheta(\mu) := \inf_x\left\{\langle \mu, x \rangle:~ x\in \solS\right\},
\end{equation}
as the best possible right hand side value for an inequality  $(\mu;\eta_0)$ to be valid for $\solS$. 
We say that a valid inequality $(\mu;\eta_0)$ is \emph{tight} if $\eta_0=\vartheta(\mu)$. 
If both $(\mu;\eta_0)$ and $(-\mu;-\eta_0)$ are valid inequalities, then $\langle \mu,  x \rangle =\eta_0$ holds for all $x\in\solS$, and  in this case, we refer to $(\mu;\eta_0)$ as a \emph{valid equation}  for $\solS$. 
{
We let $\viPi\subset E$ be the set of all nonzero vectors $\mu\in E$ such that $\vartheta(\mu)$ is finite. This set $\viPi$ is precisely the subset of  $E$ leading to nontrivial valid inequalities for $\solS$.
}

Let $\coneC  \subset E\times\bR$ denote the convex cone of all valid inequalities given by $(\mu;\eta_0)$. Identifying valid linear inequalities that are necessary in the description of $\clconv(\solS)$ is equivalent to studying $\coneC$ and its generators. 
Because $\coneC$ is a convex cone in $E\times\bR$, it can be written as the sum of a linear subspace $L$ of $E\times\bR$ and a pointed cone $C$, i.e., $\coneC=L+C$. 
Given $L$, the largest linear subspace contained in $\coneC$, let $L^\perp$ denote the orthogonal complement of $L$. Then a unique representation for the pointed cone $C$ in $\coneC=L+C$ is given by  $C=\coneC \cap L^\perp$. 
A \emph{generating set} $(G_L, G_C)$ for a cone $\coneC$ is a minimal set of elements $(\mu;\eta_0)\in \coneC$ such that $G_L \subseteq L$,  $G_C \subseteq C$, and 
\[
\coneC= \left\{\sum_{w\in G_L} \alpha_w w + \sum_{v\in G_C} \lambda_v v: ~\lambda_v \geq 0  \right\}.
\] 
\begin{remark}\label{Remark:GeneratingSet}
From this definition, it is clear that in a generating set $(G_L, G_C)$ of $\coneC$,  without loss of generality, we can assume that each vector from $G_C$ is orthogonal to every vector in $G_L$, and all vectors in $G_L$ are orthogonal to each other.\epr
\end{remark}

Our study of $\coneC$ will be based on characterizing the properties of the elements of its generating sets. We will refer to the vectors in $G_L$ as \emph{generating equalities} and the vectors in $G_C$  as \emph{generating inequalities} of $\coneC$. An inequality $(\mu;\eta_0)\in \coneC$ is called an \emph{extreme inequality}  of $\coneC$, if there exists a generating set for $\coneC$ including $(\mu;\eta_0)$ as a generating inequality either in $G_L$ or in $G_C$. When the cone $\coneC$ is pointed, we have $G_L$ is trivial  and $G_C$ is uniquely defined up to positive scalings. Then, our definition of extreme inequalities based on generating inequalities matches precisely with the usual definition of extreme inequalities stated as ``an inequality is extreme if it cannot be written as  average  of two other distinct valid inequalities." 
Note that  any non-tight valid inequality $(\mu;\eta_0)$ with $\eta_0<\vartheta(\mu)$ does not belong to a generating set of  $\coneC$.  

Clearly, the inequalities in generating set $(G_L, G_C)$ of the cone $\coneC$ are of great importance; they are necessary and sufficient for the description of $\clconv(\solS)$.   
It is easy to note that $G_L$ is finite, as a basis of the subspace $L$ can be taken as $G_L$.  For nonpolyhedral (nonlinear) cones such as $\mathcal{L}^n$ with $n\geq 3$, $G_C$ need not be finite. In fact we provide an example demonstrating this in section \ref{sec:subadditive}.

\subsection{Outline}
The main body of this paper is organized as follows. 
In section \ref{sec:minimal}, we introduce the class of  $\cK$-minimal inequalities and show that under a mild assumption, this class of inequalities together with the constraint $x\in\cK$  is sufficient to describe $\clconv(\solS)$. We follow this by establishing a number of necessary conditions for $\cK$-minimality. 
In particular, we show that $\cK$-minimal inequalities are tight in many cases. Nonetheless, we highlight that depending on the structure of $\solS$, $\cK$-minimality does not necessarily imply tightness of the inequality.  
In addition to this, one of our necessary conditions for $\cK$-minimality leads us to our next class of valid inequalities, $\cK$-sublinear inequalities. We study $\cK$-sublinear inequalities  in section \ref{sec:subadditive} and establish a precise relation between $\cK$-sublinearity and $\cK$-minimality and show that the set of extreme inequalities in the cone of $\cK$-sublinear inequalities contains all of the extreme inequalities from the cone of $\cK$-minimal inequalities.  
In section \ref{sec:support}, we show that every $\cK$-sublinear inequality is  associated with a convex set of particular structure, which we refer to as a \emph{cut generating set}. Moreover, we show that any nonempty cut generating set leads to a valid inequality.
Through this connection with structured convex sets, we provide necessary conditions for $\cK$-sublinearity, as well as sufficient conditions for a valid inequality to be $\cK$-sublinear and $\cK$-minimal. 
In the case of $\cK=\bR^n_+$, our necessary condition and sufficient condition for $\cK$-sublinearity match precisely establishing a strong relation between $\cK$-sublinear inequalities and the support functions of cut generating sets. This relation provides nice connections to the existing  literature, which we highlight in section \ref{sec:ConnectionToLatticeFreeSetsAndCGFs}. 
We close section \ref{sec:support} by examining the conic Mixed Integer Rounding (MIR) inequality from  \cite{Atamturk_Narayanan_10} in our framework.  
We provide some characterizations of the lineality space of $\coneC$ in section \ref{sec:equations}, and finish by stating a few further research questions. 
 
\section{$\cK$-Minimal Inequalities}\label{sec:minimal}

In this section, based on the ordering induced by the regular cone $\cK^*$, we first introduce a domination notion among valid linear inequalities for $\solS$. Based on this domination notion, we identify a relatively small class of valid linear inequalities, \emph{$\cK$-minimal inequalities}, and show that this class is nonempty under a mild technical assumption. This technical assumption is satisfied, for example, when $\conv(\solS)$ is full dimensional. 
Under this assumption, we establish that $\cK$-minimal inequalities     
along with the constraint $x\in\cK$ is sufficient to describe $\clconv(\solS)$. We then study the properties of inequalities from this class. 

We start by pointing out a trivial class of valid linear inequalities for $\solS$, which we refer as \emph{cone-implied inequalities}. These inequalities stem from the observation that $\solS\subseteq \cK$. The definition of  dual cone immediately implies that for any $\delta \in \cK^*$, the inequality $\langle \delta, x \rangle \geq 0$ is valid for $\cK$, and thus, it is also valid for $ \solS$. Therefore, $(\delta; 0) \in \coneC$ for any $\delta \in \cK^*$. Note that all cone-implied inequalities are readily captured by the constraint $x\in\cK$. Hence, they are not of great interest. In particular, unless $\conv(\solS)=\cK$, the family of cone-implied inequalities will not be sufficient to fully describe $\clconv(\solS)$. Because we have already assumed $\clconv(\solS)\neq\cK$, from now on, we focus on the characterization of valid linear inequalities that are non-cone-implied and are needed to obtain a complete description of $\clconv(\solS)$. This leads us to our definition of $\cK$-minimal inequalities.

\begin{definition}\label{minimalvi}
A valid linear inequality $(\mu;\eta_0)$ with $\mu\neq0$ and $\eta_0\in\bR$ is \emph{$\cK$-minimal} (for $\solS$) if for all valid inequalities $(\rho;\rho_0)$ for $\solS$ satisfying $\rho\neq \mu$, and $\rho\preceq_{\cK^*} \mu$, we have $\rho_0<\eta_0$.
\end{definition}

We next observe that the cone $\cK$, indeed, induces a natural dominance relation among the valid linear inequalities, and $\cK$-minimality definition is a result of this dominance relation. Let us consider a valid inequality $(\mu; \eta_0)$ which is not $\cK$-minimal. Thus, there exists another valid inequality $(\rho;\rho_0)$ such that $\rho\neq \mu$, $\rho\preceq_{\cK^*} \mu$, and $\rho_0\geq\eta_0$. But then the inequality $(\rho;\rho_0)$ together with the constraint $x\in\cK$ implies the inequality $(\mu; \eta_0)$ because 
\[
\langle \mu, x \rangle = \langle \rho+(\mu-\rho), x \rangle =  \underbrace{\langle \rho, x \rangle}_{\geq \rho_0} + \underbrace{\langle \mu-\rho, x \rangle}_{\geq 0} \geq \rho_0 \geq \eta_0,
\]
where the first inequality follows from $x\in\cK$ and $\mu-\rho\in\cK^*$. 
The above relation indicates that when the constraint  $x\in\cK$  and the linear inequality $(\rho;\rho_0)$ are included, the non-$\cK$-minimal inequality $(\mu; \eta_0)$ is not necessary in the description of $\clconv(\solS)$. 
The definition of $\cK$-minimality simply requires an inequality not to be dominated in this fashion: a $\cK$-minimal inequality $(\mu;\eta_0)$ cannot be dominated by another inequality, which is the sum of a cone-implied inequality and another valid inequality for $\solS$. 

In general, there are $\cK$-minimal inequalities that are not extreme. In particular, the definition of $\cK$-minimality allows for a $\cK$-minimal inequality to be implied by the sum of two other non-cone-implied valid inequalities. That said, under a technical assumption, we will show that all non-cone-implied extreme inequalities are $\cK$-minimal. Because characterization of extreme inequalities in general is known to be  a much more difficult task, in this paper, we limit our focus on the characterization of $\cK$-minimal inequalities.

\begin{remark}
None of the cone-implied inequalities $(\mu;\eta_0)=(\delta; 0)$ with $\delta \in \cK^*\setminus \{0\}$ is $\cK$-minimal because we can always write them as the sum of a valid inequality $(\rho;\rho_0)=({1\over 2}\delta; 0)$ with $\rho_0=\eta_0$ and a cone-implied inequality $({1\over 2}\delta; 0)$.   
Nevertheless, a cone-implied inequality can be extreme\myfootnote{See section \ref{sec:classOfInequalities} and the definition of extreme inequalities based on generating inequalities.}, and thus, necessary in the description of $\clconv(\solS)$.   
\epr
\end{remark}

\begin{remark}
In the case of MILP, $\cK=\bR^n_+$, a \emph{minimal inequality} is defined as a valid linear inequality $(\mu;\eta_0)$ such that if $\rho\leq \mu$ (where the $\leq$ is interpreted in the component-wise sense) and $\rho\neq \mu$, then $(\rho; \eta_0)$ is not valid, i.e., reducing any $\mu_i$ for $i\in\{1,\ldots,n\}$ will lead to a strict reduction in the right hand side value of the inequality (cf.\@ \cite{Johnson_81}). 
Considering that $\bR^n_+$ is a regular and also self-dual cone, we conclude that  
$\cK$-minimality definition is indeed a natural extension of the \emph{minimality} definition of valid inequalities studied in the context of  MILPs to more general disjunctive conic sets with regular cones $\cK$. 
\epr
\end{remark}

{
\begin{remark}\label{rem:K-minimalEmphasis}
Whether a valid inequality is necessary for the description of $\clconv(\solS)$ depends on $\solS$ and it can very well be independent of the choice of $A,\cB$ and $\cK$.  In particular, if there exists $A',\cB'$ and $\cK'$ such that $\clconv(\cS(A',\cK',\cB'))=\clconv(\solS)$, then the extreme inequalities for these will be the same. Additionally,  as long as the set $\solS$ remains the same $\cK$-minimality definition is independent of $A$ and $\cB$ but depends on $\cK$ explicitly, that is the $\cK$-minimal inequalities for both $\solS$ and $\cS(A',\cK,\cB')$ are the same as long as $\clconv(\cS(A',\cK,\cB'))=\clconv(\solS)$. However, when $\cK'\neq \cK$, $\cK$-minimal inequalities for $\solS$ might differ from $\cK'$-minimal inequalities for $\cS(A',\cK',\cB')$ even when $\clconv(\cS(A',\cK',\cB'))=\clconv(\solS)$.  We comment more on the choice of the cone $\cK$ and its impact on identifying dominance relations and our minimality notion in Remark \ref{rem:HowToDefineMinimality}.
\epr
\end{remark}

In the light of this remark, from now on we will 
emphasize the classification of valid inequalities based on the cone $\cK$ explicitly. 
}

We let $\coneCm$ denote the set of $\cK$-minimal valid inequalities for $\solS$. Note that $\coneCm$ is closed under positive scalar multiplication and is thus a cone (but it is not necessarily a convex cone). 

The following simple example shows a set $\solS$ together with the $\cK$-minimal inequalities describing its convex hull.

\begin{example}\label{ex:basicminimalvi}
Let $\solS$ be defined with $\cK=\mathcal{L}^3=\cK^*$,  $A=[-1, 0, 1]$ \myfootnote{Throughout this paper, we use Matlab notation with brackets $[ \cdot ]$ to denote explicit vectors and matrices.}  and  $\cB=\left\{0,2\right\}$, i.e., 
\[
\solS=\{x\in\cK: ~-x_1+x_3=0\} \bigcup \{x\in\cK: ~-x_1+x_3= 2\}.
\]
Then 
\bse
\Conv(\solS) &= &\{x\in\bR^3:~x\in\cK, ~0\leq -x_1+x_3\leq 2\} \\
&= &\{x\in\bR^3:~ \langle x, \delta \rangle \geq 0~\forall \delta \in \Ext(\cK^*), ~ x_1-x_3\geq -2\},
\ese
is closed, and thus, the cone of valid inequalities is given by
\[
\coneC=\cone\left( \cK^*\times\{0\}, ([1; 0; -1];-2) \right).
\]

The only non-cone-implied extreme inequality in this description is given by $\mu=[1; 0; -1]$ with $\eta_0=-2=\vartheta(\mu)$. It is easy to see that this inequality is valid and also necessary for the description of the convex hull. In order to verify that it is in fact $\cK$-minimal, consider any $\delta\in\cK^*\setminus\{0\}$, and set  $\rho=\mu-\delta$. Then the best possible right hand side value $\rho_0$ for which $\langle \rho,x\rangle \geq \rho_0$ is valid, is given by
\bse
\rho_0 &:=& \inf_x\{\langle \rho, x \rangle:~ x\in\solS \} \\
&\leq&  \inf_x\{\langle \rho, x \rangle:~ x\in\cK, ~ -x_1+x_3 = 2 \}  \\
&=&  \inf_x\{x_1-x_3 - \langle \delta, x \rangle:~ x\in\cK, ~ -x_1+x_3 = 2 \}  \\
&=&  \inf_x\{-2 - \langle \delta, x \rangle:~ x\in\cK, ~ -x_1+x_3 = 2 \}  \\
&=&  -2- \sup_x\{ \langle \delta, x \rangle:~ x\in\cK, ~ -x_1+x_3 = 2 \}  \\
&<& -2 =\vartheta(\mu),
\ese
where the strict inequality follows from the fact that $u=[0;1;2]$ is in the interior of $\cK$ and satisfies $-u_1+u_3 = 2$ (and thus is feasible to the last optimization problem in the above chain), and also for any $\delta\in\cK^*\setminus\{0\}$, $\langle \delta, u \rangle >0$.
Clearly, all of the other inequalities involved in the description of $\conv(\solS)$ are of the form $\langle \delta, x\rangle \geq 0$ with $\delta\in\Ext(\cK^*)$, and hence, are not $\cK$-minimal. \epr
\end{example}

Our goal is to generalize Example \ref{ex:basicminimalvi} and establish that $\cK$-minimal inequalities along with the constraint $x\in\cK$ are sufficient to describe $\clconv(\solS)$. However, we need a structural assumption for this result. This assumption is a result of the important fact that there can be situations where none of the inequalities describing $\clconv(\solS)$  is $\cK$-minimal even when $\clconv(\solS)\subsetneq \cK$. To emphasize  this technical difficulty and motivate our assumption, let us consider a slightly modified version of Example \ref{ex:basicminimalvi} with a different set  $\cB$:

\begin{example}\label{ex:nominimalvi}
Let $\solS$ be defined with $\cK=\mathcal{L}^3$, $A=[-1, 0, 1] \mbox{ and } \cB=\left\{0\right\}$. Then 
\[
\Conv(\solS)=\{x\in\bR^3:~x\in\cK, ~ -x_1+x_3=0\}=\{x\in\bR^3:~ x_1=x_3, ~x_2=0, ~x_1,x_3\geq0\}.
\] 

We claim and prove that none of the inequalities in the description of $\Conv(\solS)$ is $\cK$-minimal. To observe this, let us fix a particular generating set $(G_L,G_C)$ for the cone $\coneC$.  Based on the above representation of $\conv(\solS)$, we can take for example  $G_C=\cL^3\times\{0\}$ and $G_L=(\mu;0)$ where $\mu=[-1; 0; 1]$ with $\eta_0=0=\vartheta(\mu)$.  Note that all of the inequalities in $G_C$ as well as one side of the valid equation given by $(\mu;0)$ are cone-implied (because $\mu\in\cL^3$), and thus are not $\cK$-minimal. Moreover, the inequality given by $(-\mu;0)$, e.g., the other side of the valid equation also cannot be $\cK$-minimal since $\rho=[1.5; 0; -1.5]$ satisfies $\delta=-\mu-\rho=[-0.5; 0; 0.5]\in\Ext(\cK^*)$ and $(\rho;\eta_0)$ is also valid.  In fact, for any valid inequality $(\mu;\eta_0)$ that is in the description of $\Conv(\solS)$, there exists $\tau>0$ such that we can subtract the vector $\delta=\tau[-1;0;1]\in\Ext(\cK^*)$ from $\mu$, and still obtain $(\mu-\delta;\eta_0)$ as a valid inequality. Finally, note that the generators of $\coneC$ are uniquely defined up to shifts by the vector $(\mu;0)$ defining the valid equation. But these shifts do not change the $\cK$-minimality properties of the inequalities. \epr
\end{example}

The peculiar situation of Example \ref{ex:nominimalvi} is a result of the fact that $\solS\subset \{x\in\cK:  -x_1+x_3=0 \}$, i.e., $\solS$ is contained in a subspace defined by a cone-implied valid equation. The next proposition formally states that this is precisely the situation in which none of the valid linear inequalities, including the extreme ones, is $\cK$-minimal. 

\begin{proposition}\label{prop:minimalexistence}
Suppose that there exists $\delta \in \cK^*\setminus\{0\}$ such that $\langle \delta, x\rangle=0$ for all $x\in\solS$, i.e., $(\delta;0)$ is a valid equation. Then $\coneCm=\emptyset$.
\end{proposition}
\myproof{
Let $\delta\in \cK^*\setminus\{0\}$ be such that $(\delta;0)$ is  a valid equation. Consider any valid inequality $(\mu;\eta_0)$.  Because $(-\delta;0)$ is also valid, we get $(\mu-\delta;\eta_0)$ is valid as well. But then $(\mu;\eta_0)$ is not $\cK$-minimal because  $\delta\in \cK^*\setminus\{0\}$. Given that $(\mu;\eta_0)$ was arbitrary, this implies that there is no $\cK$-minimal valid inequality under the hypothesis of the proposition.
}

Based on Proposition \ref{prop:minimalexistence}, in the remainder of this paper, we make the following assumption:
\begin{quote}
\textbf{Assumption 1}: For each $\delta\in\cK^*\setminus \{0\}$, there exists some $x_\delta\in\solS$ such that $\langle \delta, x_\delta \rangle >0$.
\end{quote}
Note that \textbf{Assumption 1} is indeed not very restrictive, and is trivially satisfied, for example, when $\clconv(\solS)\neq\cK$ and is full-dimensional, e.g., when $\ker(A)\cap\intt(\cK) \neq \emptyset$ (see Proposition \ref{prop:MinimalViRhs}).

Our main result in this section shows that 
under \textbf{Assumption 1}, $\cK$-minimal inequalities, along with the constraint $x\in\cK$, are sufficient to describe $\clconv(\solS)$. In particular, we prove that under \textbf{Assumption 1}, all extreme inequalities are $\cK$-minimal. 
Due to the previous discussion on the dominance relation among inequalities and $\cK$-minimality, this result is expected. However, to formalize this, we need the following definition:   
Given two vectors, $u,v \in C$ where $C$ is a cone with lineality space $L$, $u$ is said to be an \emph{$L$-multiple} of $v$ if $u=\tau v+\ell$ for some $\tau>0$, and $\ell\in L$. From this definition, it is clear that if $u$ is an $L$-multiple of $v$, then $v$ is also an $L$-multiple of $u$. Also, we need the following lemma from \cite{Johnson_81}:
\begin{lemma}\label{lem:Johnson81}
Suppose $v$ is in a generating set for cone $C$ and there exist $v^1,v^2\in C$ such that $v=v^1+v^2$, then $v^1,v^2$ are $L$-multiples of $v$. 
\end{lemma}

Let $(G_L,G_C)$ be a generating set for the cone $\coneC$. Note that whenever the lineality space $L$ of the cone $\coneC$ is nontrivial, the generating valid inequalities are only defined uniquely up to the $L$-multiples.  
We define $G^+_C$ to be the vectors from $G_C$ that are not $L$-multiples of any cone-implied inequality $(\delta;0)$ with $\delta\in \cK^*\setminus\{0\}$. Then $G^+_C$ is again only uniquely defined up to $L$-multiples.

The following result is a straightforward extension of the associated result from \cite{Johnson_81} given in the linear case to our conic case.

\begin{proposition}\label{prop:minimalcharacterization}
Let $(G_L,G_C)$ be a generating set for the cone $\coneC$. Under {\bf Assumption 1}, every valid equation in $G_L$ and every generating valid inequality in $G^+_C$ is $\cK$-minimal.
\end{proposition}
\myproof{
Suppose $(\mu;\eta_0)\in G_L\cup G^+_C$ is not $\cK$-minimal. Then there exists a nonzero $\delta\in \cK^*$ such that $(\mu-\delta; \eta_0)\in\coneC$. Note that $(\delta;0)\in\coneC$, therefore, $(\mu+\delta;\eta_0)$ is valid as well. Then Lemma \ref{lem:Johnson81} implies that  $(\delta;0)$ is an $L$-multiple of $(\mu;\eta_0)$. Using the definition of $G^+_C$, we get $(\mu;\eta_0)\in G_L$. Given that $(\delta;0)$ is an $L$-multiple of $(\mu;\eta_0)$ and $G_L$ is uniquely defined up to $L$-multiples, we get that $(\delta;0)\in G_L$. Hence, $\langle \delta,x\rangle = 0$ is a valid equation, which contradicts to {\bf Assumption 1}.
}

Based on Proposition \ref {prop:minimalcharacterization}, \textbf{Assumption 1} ensures that $\coneCm\neq\emptyset$. In particular, Proposition \ref {prop:minimalcharacterization} immediately implies the following result. 
\begin{corollary}\label{corr:minimalcharacterization}
Suppose that \textbf{Assumption 1} holds. Then, for any generating set $(G_L,G_C)$ of $\coneC$, $(G_L,G^+_C)$ generates $\coneCm$. In particular, all non-cone-implied extreme inequalities are $\cK$-minimal. Thus, $\cK$-minimal inequalities along with the original conic constraint $x\in\cK$ are sufficient to describe $\clconv(\solS)$.
\end{corollary} 

Under \textbf{Assumption 1}, in the light of Proposition \ref{prop:minimalcharacterization} and Corollary \ref{corr:minimalcharacterization}, we arrive at 
\bse
\clconv(\solS) & = & \{x\in E:~ x\in\cK,  ~\langle \mu, x\rangle = \eta_0\, \forall (\mu;\eta_0)\in G_L, ~\langle \mu, x\rangle \geq \eta_0\, \forall (\mu;\eta_0)\in G_C^+ \} \\
& = & \{x\in E:~ x\in\cK,  ~\langle \mu, x\rangle \geq \eta_0\, \forall (\mu;\eta_0)\in \coneCm \}. 
\ese
Therefore, under \textbf{Assumption 1}, any valid inequality $(\mu;\eta_0)$  for $\solS$ is dominated by a set of $\cK$-sublinear inequalities $(\mu^i;\eta_0^i)$ where  $i \in I$ is a set indices and a cone-implied inequality $(\delta;0)$ with $\delta\in\cK^*$ (note that the cone of cone-implied inequalities is convex). That is, $\mu=\sum_{i\in I} \mu^i +\delta$ and $\eta_0\leq \sum_{i\in I} \eta_0^i$. When $\coneCm$ is convex, the set of indices $I$ can be taken as a singleton.

Next, we deliberate on the importance of the cone $\cK$ in establishing dominance relations and in our $\cK$-minimality definition.

\begin{remark}\label{rem:HowToDefineMinimality}
Based on Remark \ref{rem:K-minimalEmphasis} and our $\cK$-minimality notion, the structural information encoded in the cone $\cK$ is rather important in identifying smaller classes of valid inequalities  that are sufficient to describe the closed convex hulls of disjunctive conic sets. To emphasize this, let us consider a situation where we are given $A$ and $\cB$, and we have several options for the cone $\cK$ to encode $\solS$. Suppose that we are given two cones $\cK_1\subset \cK_2$ such that $\cS(A,\cK_1,\cB)=\cS(A,\cK_2,\cB)$. Then $C(A,\cK_1,\cB)= C(A,\cK_2,\cB)$. In such a case,  the  smaller cone $\cK_1$ encodes the structural information of the disjunctive conic set $\cS(A,\cK_1,\cB)$ better than $\cK_2$.  
In order to avoid technical difficulties let us assume that $\cS(A,\cK_1,\cB)$ satisfies \textbf{Assumption 1} with respect to $\cK_1$ and $\clconv(\cS(A,\cK_1,\cB))\neq \cK_1$, thus $C_m(A,\cK_1,\cB)$ is nonempty. The definition of $\cK$-minimality together with the relation $\cK_1^*\supset \cK_2^*$  automatically implies that $\cK_1$-minimal inequalities for $\cS(A,\cK_1,\cB)$ are also $\cK_2$-minimal for $\cS(A,\cK_1,\cB)$, but not vice versa because $\cK_1\neq\cK_2$. Therefore, $C_m(A,\cK_1,\cB)\subsetneq C_m(A,\cK_2,\cB)$. Let $(G_L,G_C)$ be a generating set for $C(A,\cK_1,\cB)$. Let us define $G^{1,+}_C$ to be the vectors from $G_C$ that are not $L$-multiples of any cone-implied inequality $(\delta;0)$ with $\delta\in \cK_1^*\setminus\{0\}$, and $G^{2,+}_C$ analogously with respect to $\cK_2^*$. Then by Corollary \ref{corr:minimalcharacterization}, we have $(G_L,G^{1,+}_C)$ generates $C_m(A,\cK_1,\cB)$ and $(G_L,G^{2,+}_C)$ generates $C_m(A,\cK_2,\cB)$. Because $C_m(A,\cK_1,\cB)\subsetneq C_m(A,\cK_2,\cB)$, we conclude that $G^{1,+}_C \subsetneq G^{2,+}_C$. Then,  all extreme $\cK_2$-minimal inequalities are also $\cK_1$-minimal, but some $\cK_2$-minimal inequalities, namely $G^{2,+}_C\setminus G^{1,+}_C$, will not be extreme, and
\bse
\clconv(\solS) & = & \{x\in E:~ x\in\cK_1,  ~\langle \mu, x\rangle = \eta_0~ \forall (\mu;\eta_0)\in G_L, ~\langle \mu, x\rangle \geq \eta_0~\forall (\mu;\eta_0)\in G_C^{1,+} \} \\
& = & \{x\in E:~ x\in\cK_2,  ~\langle \mu, x\rangle = \eta_0~ \forall (\mu;\eta_0)\in G_L, ~\langle \mu, x\rangle \geq \eta_0~ \forall (\mu;\eta_0)\in G_C^{2,+} \}.  
\ese
Hence, we conclude that whenever we have a choice between $\cK_1\subset \cK_2$, minimality defined with respect to the smaller cone $\cK_1$ results in a stronger dominance notion among valid linear inequalities defining $\clconv(\solS)$.

As a consequence of this, we highlight the importance of encoding structural information in $\cK$ as much as possible. For example, among different choices of disjunctive conic representations of the same set suggested in Examples \ref{ex:MICPform2} and \ref{ex:MICPform3}, the one in Example \ref{ex:MICPform3} is superior. This is so, even when the cone $\tilde{\cK}$ is as simple as $\bR^n_+$. Therefore, even in the case of MILPs, whenever such structural information, e.g., a polyhedral relaxation, is present, there is benefit in defining minimality notion based on a regular cone $\cK$ defined in a lifted space as described in Example \ref{ex:MICPform3} as opposed to the usual choice of nonnegative orthant from the MILP literature.   
\epr
\end{remark}

These results motivate us to further study the properties of $\cK$-minimal inequalities in the next section.

\subsection{Necessary Conditions for $\cK$-Minimality}

Our first proposition states that in certain cases, all $\cK$-minimal inequalities are tight. This also gives us our first necessary condition for $\cK$-minimality.

\begin{proposition}\label{prop:minVIrhs}
Let  $(\mu;\eta_0)\in\coneCm$. Then,  whenever 
$\mu\in\cK^*$ or $\mu\in-\cK^*$, 
the valid inequality $(\mu;\eta_0)$ is tight, i.e., $\eta_0=\vartheta(\mu)$ (cf.\@ \eqref{mu0}). Furthermore, $(\mu;\eta_0)\in\coneCm$ and $\mu\in\cK^*$ (respectively $\mu\in-\cK^*$) implies $\vartheta(\mu)>0$ (respectively $\vartheta(\mu)<0$).
\end{proposition}
\myproof{
Consider $(\mu;\eta_0)\in\coneCm$ with $\mu\neq0$. Note that $\mu=0$ leads to trivial valid inequalities which are not of interest. The validity of $(\mu;\eta_0)$ immediately implies $\eta_0\leq\vartheta(\mu)$. Assume for contradiction that $\eta_0<\vartheta(\mu)$. We need to consider only two cases:
\begin{itemize}
\item[{(i)}] $\mu\in\cK^*\setminus\{0\}$: Then $\vartheta(\mu)\geq \eta_0>0$, because otherwise $(\mu;\eta_0)$ is either a cone-implied inequality or is dominated by a cone-implied inequality, both of which are not possible. Let $\beta={\eta_0\over\vartheta(\mu)}$, and consider $\rho=\beta\cdot\mu$. Then $(\rho;\eta_0)$ is a valid inequality because $0<\beta<1$, $(\mu;\vartheta(\mu))\in\coneC$ and $\coneC$ is a cone. But $\mu-\rho = (1-\beta)\mu \in\cK^*\setminus\{0\}$ since $\mu\neq 0$ and $\beta<1$. This is a contradiction, thus, we conclude $\eta_0=\vartheta(\mu)>0$.  

\item[{(ii)}] $-\mu\in\cK^*\setminus\{0\}$: Because $(-\mu;0)$ is trivially valid and we cannot satisfy both $(-\mu;0)$ and $(\mu;\vartheta(\mu))$ when $\vartheta(\mu)>0$ unless  $\solS=\emptyset$. But this is not possible due to our assumptions on $\solS$, thus, we conclude that $\vartheta(\mu)\leq 0$. 
Moreover, if $\vartheta(\mu)=0$, then $\solS\subset \{x\in\cK:~\langle \mu, x\rangle = 0\}$, which contradicts \textbf{Assumption 1}. Hence, we conclude that $\eta_0<\vartheta(\mu)<0$. Once again let $\beta={\eta_0\over\vartheta(\mu)}$, and consider $\rho=\beta\cdot\mu$. Then $(\rho;\eta_0)$ is a valid inequality since $\beta>1$, $(\mu;\vartheta(\mu))\in\coneC$ and $\coneC$ is a cone. But $\mu-\rho = (1-\beta)\mu \in\cK^*\setminus\{0\}$ since $\mu\in-\cK^*\setminus\{0\}$  and $\beta>1$. But, this is a contradiction to the $\cK$-minimality of $(\mu;\eta_0)$. Thus, we conclude that $\eta_0=\vartheta(\mu)<0$.
\vspace{-13pt}
\end{itemize}
}

Clearly, Proposition \ref{prop:minVIrhs} does not cover all possible cases for $\mu$. As a matter of fact, it is possible to have $\mu\not\in\pm \cK^*$ leading to a $\cK$-minimal inequality. While one is naturally inclined to believe that a $\cK$-minimal inequality $(\mu;\eta_0)$ is always tight, i.e., $\eta_0=\vartheta(\mu)$, we have the following counter-example.
\begin{example}\label{ex:MinimalViRhs}
Consider the disjunctive conic set $\solS$ defined with $A=[-1, 1]$, $\cB=\{-2,1\}$ and $\cK=\bR^2_+$. First, note that \textbf{Assumption 1} holds because $\{[0;1],[2;0]\}\in\solS$, and $\conv(\solS)=\clconv(\solS)\neq\bR^2_+$. Thus, $\cK$-minimal inequalities exist, and together with nonnegativity restrictions they are sufficient to describe $\conv(\solS)$. 
In fact, 
\[
\conv(\solS)=\{x\in\bR^2:~ -x_1+x_2 \geq -2, ~ x_1-x_2\geq -1, ~x_1+2x_2 \geq 2, ~ x_1,x_2\geq 0\},
\]
and one can easily show that each of the nontrivial inequalities in this description is in fact $\cK$-minimal. 

\begin{figure}[ht!]
\begin{center}
\begin{tikzpicture}[scale=1.2]
\fill[blue!50!white, thick] (2,0) -- (0,1) -- (1.5,2.5) -- (2.9,2.5) -- (2.9,1) -- cycle;
\draw[blue!70!white, very thick]  (2.9,1) -- (2,0) -- (0,1) -- (1.5,2.5);
\draw[blue!70!black,  very thick]  (-2.4,-0.4) -- node[sloped,above] {\tiny{$x_1-x_2\geq -2$}}  (0,2) -- (0.4,2.4);
\draw[blue!70!black,    very thick, ->]  (-2.4,-0.4) -- (-2.2,-0.6) ;
\draw[blue!70!black,    very thick, ->]  (0.4,2.4) -- (0.6,2.2) ;

\draw[ step = 0.5cm, gray, dashed, very thin] ( -2.9 ,-0.7 ) grid ( 2.9, 2.4 );
\draw[<->] ( 0, 2.5 ) node[ anchor = south ] {$ x_2 $} -- ( 0, -0.8 ); 
\draw[<->] ( -3, 0 ) -- ( 3, 0 ) node[ anchor =  south west ] {$ x_1 $};

\fill ( 0, 0 ) circle ( 0.04 );
\draw ( 0, 0 ) node[anchor = south east] {\tiny{$ ( 0, 0 ) $}};
\fill ( 0, 1 ) circle ( 0.04 );
\draw ( 0, 1 ) node[anchor = south east] {\tiny{$ ( 0, 1 ) $}};
\fill ( 0, 2 ) circle ( 0.04 );
\draw ( 0, 2 ) node[anchor = south east] {\tiny{$ ( 0, 2) $}};
\fill ( -2, 0 ) circle ( 0.04 );
\draw (- 2, 0 ) node[anchor = north] {\tiny{$ ( -2, 0) $}};
\fill ( -1, 0 ) circle ( 0.04 );
\draw (- 1, 0 ) node[anchor = north] {\tiny{$ ( -1, 0) $}};
\fill ( 1, 0 ) circle ( 0.04 );
\draw ( 1, 0 ) node[anchor = north] {\tiny{$ ( 1, 0) $}};
\fill ( 2, 0 ) circle ( 0.04 );
\draw ( 2, 0 ) node[anchor = north] {\tiny{$ ( 2, 0) $}};

\end{tikzpicture}
\end{center}
\caption{Convex hull of $\solS$ for Example \ref{ex:MinimalViRhs}}
\label{fig:ex:MinimalViRhs}
\end{figure}

Now, let us consider the valid inequality given by $(\mu;\eta_0)=([1;-1];-2)$. Note that $\vartheta(\mu)=-1$, therefore, $(\mu;\eta_0)$ is not tight and is  dominated by the valid inequality $x_1-x_2\geq -1$. We will show that $(\mu;\eta_0)$ is $\cK$-minimal regardless of the fact that it is not tight. We note that, in this example, $\cK$-minimality is the same as the usual minimality used in the usual MILP literature.  

Suppose that $(\mu;\eta_0)$ is not $\cK$-minimal, then there exists $\rho=\mu-\delta$ with $0\neq \delta\in\cK^*=\bR^2_+$ such that $(\rho;\eta_0)$ is a valid inequality. This implies  
\bse
-2=\eta_0 &\leq& \inf_{x} \{\langle \rho,x\rangle:~x\in\solS\} = \min_x \{\langle \rho,x\rangle:~x\in\conv(\solS)\} \\ 
& = & \min_x  \{\langle \rho,x\rangle:~ -x_1+x_2 \geq -2, ~ x_1-x_2\geq -1, ~x_1+2x_2 \geq 2, ~ x_1,x_2\geq 0 \} \\
& = & \max_\lambda \{-2\lambda_1 -\lambda_2 +2\lambda_3:~ -\lambda_1+\lambda_2+\lambda_3\leq \rho_1,~ \lambda_1-\lambda_2+2\lambda_3\leq \rho_2,~ \lambda\in\bR^3_+ \} \\
& = & \max_\lambda \{-2\lambda_1 -\lambda_2 +2\lambda_3:~ -\lambda_1+\lambda_2+\lambda_3\leq 1-\delta_1,~ \lambda_1-\lambda_2+2\lambda_3\leq -1-\delta_2,~ \lambda\in\bR^3_+ \},
\ese
where the third equation follows from strong duality (the primal problem is  feasible), and the fourth equation follows from the definition of $\rho=\mu-\delta$. On the other hand, the following system 
\bse
\lambda&\geq& 0\\
\lambda_1-\lambda_2-\lambda_3 &\geq& \delta_1-1\\
-\lambda_1+\lambda_2-2\lambda_3 &\geq& 1+\delta_2,
\ese
implies that $0\geq -3\lambda_3 \geq \delta_1+\delta_2 $. Considering that $\delta\in\bR^2_+$, this leads to $\delta_1=\delta_2=0$, which is a contradiction to $\delta\neq 0$. Therefore, we conclude that $(\mu;\eta_0)=([1;-1];-2)\in\coneCm$ yet $\eta_0\neq \vartheta(\mu)$.
\epr
\end{example}

\begin{remark}\label{rem:non-tightness}
This issue of non-tightness of some $\cK$-minimal inequalities is independent of whether the $\cK$-minimal inequality separates the origin or not. When we consider a variation of Example \ref{ex:MinimalViRhs} given by $A=[-1, 1]$, $\cB=\{-2,-1\}$ and $\cK=\bR^2_+$, we have the valid inequality given by $(\mu;\eta_0)=([1;-1];{1\over 2})$ is $\cK$-minimal due to the same reasoning, yet, it has $\vartheta(\mu)=1$ and hence $(\mu;\eta_0)$ is not tight. Note also that this inequality separates the origin from the closed convex hull.
\epr
\end{remark}

In fact, we can generalize the situation of Example \ref{ex:MinimalViRhs}, and prove the following proposition, which states that under a special condition, $\ker(A)\cap\intt(\cK) \neq \emptyset$, any valid inequality $(\mu;\eta_0)$ with $\mu\in\Im(A^*)$ and $-\eta_0\leq\vartheta(\mu)$ (cf.\@ \eqref{mu0}) is a $\cK$-minimal inequality.  
\begin{proposition}\label{prop:MinimalViRhs}
Suppose $\ker(A)\cap\intt(\cK) \neq \emptyset$. Then, for any $\mu\in\Im(A^*)$ and any $-\infty<\eta_0\leq\vartheta(\mu)$, we have $(\mu;\eta_0)\in\coneCm$.
\end{proposition}
\myproof{
Consider $d \in \ker(A)\cap\intt(\cK) \neq \emptyset$, note that $d\neq 0$. 
For any $b\in\cB$, define the set $\cS_b:=\{x\in E:~Ax=b, ~x\in\cK\}$, and let $\widehat{\cB}:=\{b\in\cB:\,\cS_b\neq\emptyset\}$. Because $\solS\ne\emptyset$, we have $\widehat{\cB}\neq\emptyset$. For any $b\in\widehat{\cB}$, let $x_b\in\cS_b$, then $P_b:=\{x_b+\tau d:~\tau\geq 0\} \subseteq \cS_b$ holds. Moreover, $P_b\cap\intt(\cK)\neq \emptyset$ for any $b\in\widehat{\cB}\neq\emptyset$, and thus, \textbf{Assumption 1} holds here. 

Assume for contradiction that the statement is not true, i.e., there exists $\mu\in\Im(A^*)$ together with $\eta_0\leq\vartheta(\mu)$, such that $(\mu;\eta_0)\not\in\coneCm$. Then there exists $\delta\in\cK^*\setminus \{0 \}$ such that $(\mu-\delta;\eta_0)\in\coneC$, which implies
\bse
-\infty < \eta_0 & \leq & \inf_{x} \{\langle \mu-\delta, x\rangle: ~x\in\solS\} \\
& \leq & \inf_{b\in\cB} \inf_{x} \{\langle \mu-\delta, x\rangle:~ A x=b, ~x\in\cK \} \\
& \leq & \inf_{b\in\widehat{\cB}} \inf_{x} \{\langle \mu-\delta, x\rangle:~ x\in P_b \} \\
& \leq & \inf_{b\in\widehat{\cB}}\left[ \underbrace{\langle \mu-\delta, x_b\rangle}_{\in\bR}  + \inf_{\tau} \{\langle \mu-\delta, \tau d\rangle:~\tau\geq 0 \}  \right] .
\ese
Also, note that $\inf_{\tau} \{\langle \mu-\delta, \tau d\rangle:~\tau\geq 0 \} =-\infty$ when  $\langle \mu-\delta, d\rangle <0$. But  $\langle \mu-\delta, d\rangle <0$ is impossible since it would have implied $-\infty < \eta_0 \leq-\infty$. Therefore, we conclude that $\langle \mu-\delta, d\rangle\geq 0$. 

Finally, because $\mu\in\Im(A^*)$, there exists $\lambda$ such that $\mu=A^*\lambda$. Taking this into account, we arrive at  
\[
0 \leq \langle \mu-\delta, d\rangle = \langle A^*\lambda, d\rangle -  \langle \delta, d\rangle =  \lambda^T\underbrace{(Ad)}_{=0} -  \langle \delta, d\rangle = -  \langle \delta, d\rangle,
\]
where we used the fact that $d\in\ker(A)$. But $d\in\intt(\cK)$ and $\delta\in\cK^*\setminus\{0\}$ implies that $\langle \delta, d\rangle > 0$, which is a contradiction.
}

Example \ref{ex:MinimalViRhs} and Proposition \ref{prop:MinimalViRhs} indicate a weakness of the $\cK$-minimality definition. To address this, we should focus on only \emph{tight $\cK$-minimal} inequalities, that is, $(\mu;\eta_0)\in\coneCm$ where $\eta_0=\vartheta(\mu)$ ({$\eta_0$ cannot be increased without making the current inequality invalid}). While we can include a tightness requirement in our $\cK$-minimality definition, we note that  tightness has a direct characterization through $\vartheta(\mu)$, and also to remain consistent with the original minimality definition for $\cK=\bR^n_+$, we opt to work with our original $\cK$-minimality definition. As will be clear from the rest of the paper, tightness considerations will make minimal change in our analysis. 

We next state a proposition which identifies a key necessary condition for $\cK$-minimality via a certain non-expansiveness property. The following set of linear maps will be of importance for this result.
\[
\cF_{\cK}:=\{(Z:E\rightarrow E):~ Z\mbox{ is linear, and } Z^*v\in\cK ~\forall v\in\cK\},
\]
where $Z^*$ denotes the conjugate linear map of $Z$.\myfootnote{
Here, for a linear map $Z:E\rightarrow E$, we use $Z^*$ to denote its conjugate map given by the identity 
\[
\langle x, Zv\rangle = \langle Z^*x, v\rangle ~~\forall (x\in E, v\in E). 
\]
}

\begin{proposition}\label{matrixform}
Let $(\mu;\eta_0)\in\coneC$ and suppose that there exists a linear map $Z\in \cF_{\cK}$ such that $AZ^*=A$,  and $\mu-Z\mu\in \cK^*\setminus\{0\}$. Then  $(\mu;\eta_0)\not \in \coneCm$.
\end{proposition}
\myproof{
Let $(\mu;\eta_0)\in\coneC$ and $Z$ be a linear map as described in the proposition. Since $Z\in \cF_{\cK}$, for any $x\in\cK$, we have $Z^*x\in\cK$. Moreover, $AZ^*x=Ax$ due to $AZ^*=A$, and thus  for any $x \in \solS$, $AZ^*x=Ax\in\cB$. Therefore, we have $Z^*x\in \solS$ for any $x \in \solS$. 
Now, let $\delta=\mu-Z\mu$, then $\delta\in \cK^*\setminus\{0\}$ by the premise of the proposition. Define $\rho:=\mu-\delta$, then for any $x\in\solS$ we have
\[
\langle \rho, x\rangle = \langle \mu-\delta, x\rangle = \langle Z\mu, x\rangle =   \langle \mu, Z^*x\rangle \geq \eta_0,
\]
where the last inequality follows from the fact that $Z^*x\in\solS$ and  $(\mu;\eta_0)\in\coneC$. Hence, we get
\[
\inf_{x}\{\langle \rho, x\rangle:~ Ax \in \cB, ~x\in\cK \} \geq \eta_0,
\]
which implies that  $(\mu;\eta_0)\not \in \coneCm$ because $(\mu;\eta_0)=(\rho;\eta_0)+(\delta;0)$ with $(\rho;\eta_0)\in\coneC$ and $0\neq \delta\in\cK^*$.
}

Proposition \ref{matrixform} states an involved necessary condition for a valid  inequality to be $\cK$-minimal.  
It states that $(\mu;\eta_0)$ is a $\cK$-minimal inequality only if the following holds:
\[
(\Id-Z)\mu\not \in \cK^*\setminus\{0\}~~\forall Z\in\cF_\cK \mbox{ such that } AZ^* = A.
\]

Based on this result, the set $\cF_\cK$ has certain importance. In fact $\cF_\cK$ is the cone of $\cK^*-\cK^*$ positive maps, which also appear in applications of robust optimization, quantum physics, etc. (see \cite{BenTal_Nemirovski_98}). When $\cK=\bR^n_+$,  $\cF_\cK=\{Z\in\bR^{n\times n}:~ Z_{ij}\geq 0 ~\forall i,j\}$. However, in general, the description of $\cF_\cK$ can be rather nontrivial for different cones $\cK$. In fact, in \cite{BenTal_Nemirovski_98}, it is shown that deciding whether a given linear map takes $\cS^n_+$ to itself is an NP-Hard optimization problem. In another particular case of interest, when $\cK=\mathcal{L}^n$, a quite nontrivial explicit description of $\cF_\cK$ via linear matrix inequalities is given by Hildebrand in \cite{Hildebrand_07,Hildebrand_11}. 
Due to the general difficulty of characterizing $\cF_\cK$, and thus, testing the necessary condition of $\cK$-minimality given in Proposition \ref{matrixform}, in the next section, we study a relaxed version of the condition from Proposition \ref{matrixform}. This leads to a larger class of valid inequalities, namely \emph{sublinear inequalities}, which subsumes the class of $\cK$-minimal inequalities. 

\section{$\cK$-Sublinear Inequalities}\label{sec:subadditive}
\begin{definition}\label{def:sublinearVi}
An inequality $(\mu;\eta_0)$ with $\mu\neq0$  and $\eta_0\in\bR$ is \emph{$\cK$-sublinear} (for $\solS$)
if it satisfies the conditions (\textbf{A.1($\alpha$)}) for all $\alpha\in\Ext(\cK^*)$ and (\textbf{A.2}) where
\bse
\mbox{(\textbf{A.1($\alpha$)})} && 0 \leq \langle \mu, u \rangle  \mbox{ for all } u\in E \mbox{ s.t. } Au=0  \mbox{ and }  \langle \alpha, v\rangle u+ v \in\cK ~~\forall v \in \Ext(\cK),\\ 
\mbox{(\textbf{A.2})} && \mu_0 \leq \langle \mu, x \rangle \mbox{ for all }  x\in \solS.
\ese
\end{definition}
When an inequality satisfies (\textbf{A.1($\alpha$)}) for all $\alpha\in\Ext(\cK^*)$ we say that it satisfies condition (\textbf{A.1}).  

It can be easily verified that the set of $(\mu;\eta_0)$ satisfying conditions (\textbf{A.1})-(\textbf{A.2}) in fact leads to a convex cone in the space $E\times\bR$. We denote this cone of $\cK$-sublinear inequalities with $\coneCs$. 

Condition (\textbf{A.2}) is simply included to ensure the validity of a given inequality, and thus, it is satisfied by every valid inequality. On the other hand, condition (\textbf{A.1}) is not very intuitive. The main role of condition (\textbf{A.1}) is to ensure the necessary non-expansivity condition for $\cK$-minimality established in Proposition \ref{matrixform}. 

There is a particular and simple case of (\textbf{A.1}) that is of interest and deserves a separate treatment: 

Let $(\mu;\eta_0)$ satisfy (\textbf{A.1}), then  $(\mu;\eta_0)$ also satisfies the following condition:
\bse
(\textbf{A.0}) && 0 \leq \langle \mu, u \rangle \mbox{ for all } u\in\cK \mbox{ such that } Au=0.
\ese
In order to see that in fact (\textbf{A.0})  is a special case of (\textbf{A.1}), consider    any $u\in\cK\cap\ker(A)$. Then, for any $\alpha\in\Ext(\cK^*)$, we have $\langle \alpha, v\rangle \geq 0$ for all $v\in\Ext(\cK)$, and  because $u\in\cK$ and $\cK$ is a cone, the requirement of condition (\textbf{A.1}) on $u$, is automatically satisfied for any such $u\in\cK\cap\ker(A)$. 

While condition (\textbf{A.1}) immediately implies (\textbf{A.0}),  treating (\textbf{A.0}) separately seems to be handy as some of our results depend solely on conditions (\textbf{A.0}) and (\textbf{A.2}). 
Note also that condition (\textbf{A.0}) is precisely equivalent to
\bse
(\textbf{A.0}) && \mu\in (\cK\cap \Ker(A))^* = \cK^*+(\Ker(A))^* = \cK^*+\Im(A^*),
\ese
where the last equation follows from the facts that $\Ker(A)^*=\Ker(A)^\perp=\Im(A^*)$ and $\cK^*+\Im(A^*)$ is closed whenever $\cK$ is closed \cite[Corollary 16.4.2]{Rockafellar_70}.

Condition  (\textbf{A.0}) is not as strong as  (\textbf{A.1}). Nevertheless,  condition  (\textbf{A.0}) is necessary for any non-trivial valid inequality, which we prove next. Recall that $\viPi=\{\mu\in E:~\mu\neq 0, ~\vartheta(\mu)\in\bR\}$.  

\begin{proposition}\label{prop:necessityA0}
Suppose $\mu\in \viPi$, then $\mu$ satisfies condition (\textbf{A.0}).
\end{proposition}
\myproof{
Suppose condition (\textbf{A.0}) is violated by some $\mu\in \viPi$.  Then  there exists  $u \in \cK$ such that $Au=0$ and $\langle \mu, u \rangle <0$. Note that for any $\beta>0$ and $x\in \solS$, $x+\beta u \in \cK$ and $A(x+\beta u)=Ax\in \cB$, hence $x+\beta u \in \solS$. On the other hand, the term,  
\[
\langle \mu, x+\beta u \rangle = \langle \mu, x \rangle+  \beta \langle \mu, u \rangle,  
\]
can be made arbitrarily small by increasing $\beta$, which implies $\vartheta(\mu)=-\infty$ where $\vartheta(\mu)$ is as defined in \eqref{mu0}. However, this is a contradiction because we started with $\mu\in \viPi$, and so $\vartheta(\mu)\neq -\infty$.
}
As a consequence of Proposition \ref{prop:necessityA0}, we conclude that in order to obtain $\clconv(\solS)$ one is required to add only an appropriate subset of valid inequalities $(\mu;\vartheta(\mu))$ with $\mu\in \cK^*+\Im(A^*)$. 

Our next theorem states that every $\cK$-minimal inequality is also $\cK$-sublinear. 
\begin{theorem}\label{thm:subadditive}
If $(\mu;\eta_0)\in\coneCm$, then $(\mu;\eta_0)\in\coneCs$.
\end{theorem}
\myproof{
Consider any $\cK$-minimal inequality $(\mu;\eta_0)$. Because $(\mu;\eta_0)\in\coneCm$, $(\mu; \eta_0)$ is valid for $\solS$, and hence, condition (\textbf{A.2}) for $\cK$-sublinearity is automatically satisfied.

Assume for contradiction that $(\mu;\eta_0)$ violates condition (\textbf{A.1}$(\alpha)$) for some $\alpha\in \Ext(\cK^*)$, i.e., there exists  $u$ such that $\langle \mu, u \rangle <0$, $Au=0$, and $\langle \alpha, v\rangle u+ v \in\cK ~~\forall v \in \Ext(\cK)$.  Based on $u$ and $\alpha$, let us define a linear map $Z:E\rightarrow E$ as 
\[
Zx = \langle x,u\rangle \alpha + x~~ \mbox{for any } x\in E. 
\] 
Note that $A:E\rightarrow \bR^m$ and thus its conjugate $A^*: \bR^m \rightarrow E$. We let $ A^*e^i=:A^i\in E$ for $i=1,\ldots,m$, where $e^i$ is the $i^{th}$ unit vector in $\bR^m$. This way,  for all $i=1,\ldots,m$, we have $ZA^*e^i=\langle A^i, u \rangle +A^i=A^i$ because $u\in\ker(A)$ implies $\langle A^i, u \rangle=0$. Therefore, we have $ZA^*=A^*$. Also, since $A:E\rightarrow \bR^m$ and $Z:E\rightarrow E$ are linear maps, we have $ZA^*$ is a linear map and its conjugate is given by $AZ^*=A$ as desired.

Moreover, for all $w\in\cK^*$ and $v\in\Ext(\cK)$, we note that
\[
\langle Zw, v \rangle = \left\langle (\langle w, u\rangle \alpha+ w), v \right\rangle = \langle w, u\rangle \langle \alpha, v \rangle + \langle w, v\rangle = \langle w,  \underbrace{\langle \alpha, v \rangle u + v}_{\in\cK} \rangle \geq 0.
\] 
Because any $v\in\cK$ can be written as a convex combination of points from $\Ext(\cK)$, we conclude that $Z\in\cF_{\cK}$. Finally by recalling that $\alpha \in \cK^*$ and is nonzero, we get 
\[
\mu-Z\mu= - \underbrace{\langle \mu, u \rangle}_{< 0}  \alpha \in \cK^*\setminus \{0\},
\] 
which is a contradiction to the necessary condition for $\cK$-minimality given in Proposition \ref{matrixform}.
}

The proof of Theorem \ref{thm:subadditive} reveals the importance of condition (\textbf{A.1}) and its implications in terms of $\cK$-minimality. Next, we show that condition (\textbf{A.1}) further simplifies in the case of $\cK=\bR^n_+$, and conditions (\textbf{A.0})-(\textbf{A.2}) underlie the definition of \emph{subadditive  inequalities} from \cite{Johnson_81} in the MILP case.

\begin{remark}\label{remark:connectionToJohnson}
When the cone $\cK$ has a simple structure, in particular, when it has finitely many and orthogonal to each other extreme rays, the interesting cases of condition (\textbf{A.1}) that are not covered by condition (\textbf{A.0}) can be simplified. When in addition the cone $\cK$ is assumed to be regular,  we can assume that $\cK=\bR^n_+$ without loss of generality.
 
Suppose $\cK=\bR^n_+$, then the extreme rays of $\cK$ as well as $\cK^*$ are just the unit vectors, $e^i$. Let us consider (\textbf{A.1}$(\alpha)$) for the case of $\alpha=e^i$. Then the vectors $u$ considered in the condition (\textbf{A.1}$(e^i)$) are required to satisfy   
\[
v_i u+ v \in\cK ~~\forall v \in \Ext(\cK)=\{e^1,\ldots,e^n\}.
\]
Because all of the extreme rays of $\cK$ are unit vectors, this requirement affects only the extreme rays $v$ with a nonzero $v_i$ value, which is just  the case of $v=e^i$. Hence, for $i=1,\ldots,n$, we can equivalently rewrite condition (\textbf{A.1}$(e^i)$) as follows: 
\bse
(\textbf{A.1i}) &&0 \leq  \langle \mu, u\rangle  \mbox{ for all } u \mbox{ such that } Au=0 \mbox{ and } u+ e^i \in\bR^n_+ .
\ese
Let $a^i$ denote the $i^{th}$ column of the matrix $A$. By a change of variables, this requirement is equivalent to the following relation for all $i=1,\ldots,n$:  
\bse
(\textbf{A.1i}) && \mu_i \leq \langle \mu, w\rangle \mbox{ for all } w\in\bR^n_+ \mbox{ such that }  ~Aw=a^i .
\ese

When $\cK=\bR^n_+$ and $\cB$ is  a finite set, Johnson \cite{Johnson_81} defines the class of so called \emph{subadditive valid inequalities} precisely as the inequalities that satisfy the collection of conditions (\textbf{A.1i}) for $i=1,\ldots,n$, along with the conditions (\textbf{A.0}) and (\textbf{A.2}).   
In this specific setup, Johnson \cite{Johnson_81} shows further that $\bR^n_+$-sublinearity of an inequality can be verified by checking requirements (\textbf{A.0}),  (\textbf{A.1i}) for $i=1,\ldots,n$, and (\textbf{A.2}) on only a finite set of vectors (those satisfying a minimal linear dependence condition).

Moreover, let us for a moment assume that there exists a function $\sigma(\cdot)$ underlying the $\cK$-sublinear inequality $(\mu;\eta_0)$. That is, given the data associated with variable $x_i$, namely $a^i$, $\sigma(\cdot)$  generates the corresponding coefficient in the valid inequality, i.e., $\mu_i=\sigma(a^i)$ for all $i=1,\ldots,n$. 
Then, condition (\textbf{A.1i}) above  precisely represents the subadditivity property of the function $\sigma(\cdot)$ over the columns of $A$. 
In fact, in section \ref{sec:support} for general disjunctive conic sets $\solS$ with $\cK=\bR^n_+$, we show that for every $\cK$-sublinear inequality $(\mu;\eta_0)$, such a function $\sigma(\cdot)$ generating $\mu$ always exists. In the specific case of $\cK=\bR^n_+$ and a finite set $\cB$, this connection was previously established in \cite{Johnson_81}. We discuss the implications of these with regard to existing MILP literature in detail in section \ref{sec:ConnectionToLatticeFreeSetsAndCGFs}.  
\epr
\end{remark}

Under \textbf{Assumption 1}, there is a precise relation between the generators of the cones of $\cK$-sublinear inequalities and $\cK$-minimal inequalities.  We state this below in Theorem \ref{thm:subadditivegenerators}, which is a generalization of the corresponding result from \cite{Johnson_81} to the conic case. For completeness, we include the following proof, which simultaneously simplifies and generalizes the approach of \cite{Johnson_81}. 

\begin{theorem}\label{thm:subadditivegenerators}
Suppose that \textbf{Assumption 1} holds. Then, any generating set of $\coneCs$ is of form $(G_L, G_s)$ where $G_s\supseteq G^+_C$ and $(G_L, G_C)$ is a generating set of $\coneC$. Moreover, if $(\mu;\eta_0)\in G_s\setminus G^+_C$, then $(\mu;\eta_0)$ is not $\cK$-minimal. 
\end{theorem}
\myproof{
Based on Remark \ref{Remark:GeneratingSet}, let $(G_L, G_C)$ be a generating set of $\coneC$ such that each vector in $G_C$ is orthogonal to every vector in $G_L$, and all vectors in $G_L$ are orthogonal to each other. Let $(G_\ell, G_s)$ be a generating set of $\coneCs$ in which each vector in $G_s$ is orthogonal to every vector in $G_\ell$. 
Note that by Theorem \ref{thm:subadditive}, we have $\coneCm\subseteq\coneCs\subseteq \coneC$. 

Under \textbf{Assumption 1}, using Corollary \ref{corr:minimalcharacterization}, we have $\coneCm$  has a generating set of the form $(G_L,G_C^+)$. Hence, we conclude that the subspace spanned by $G_\ell$ both simultaneously contains, and is contained in, the subspace generated by  $G_L$. Therefore,  we can take $G_\ell=G_L$. 

Let $Q$ be the orthogonal complement to the subspace generated by $G_L$ and define $C'=\coneC\cap Q$, $C'_m=\coneCm\cap Q$ and $C'_s=\coneCs\cap Q$. Then $C'=\cone(G_C)$, and under \textbf{Assumption 1}, $C'_m=\cone(G_C^+)$.
Also, $C', C'_m$ and $C'_s$ are pointed cones and satisfy $C'_m\subseteq C'_s \subseteq C'$. 
Given that the elements of $G_C^+$ are extreme in both $C'$ and $C'_m$, they remain extreme in $C'_s$ as well. Therefore, $G_C^+ \subseteq G_s$. 

Finally, consider any $(\mu;\eta_0)\in G_s\setminus G_C^+$. We need to show that $(\mu;\eta_0)\not \in \coneCm$. Suppose not, then $(\mu;\eta_0)\in \coneCm$ but not in $G_C^+$, which implies that $(\mu;\eta_0)$ is not extreme in $\coneCm$. Noting $\coneCm\subseteq\coneCs$, we conclude that $(\mu;\eta_0)$ is not extreme in $\coneCs$ as well. But this is a contradiction to the fact that $(\mu;\eta_0)\in G_s$ and $(G_L, G_s)$ is a generating set for $\coneCs$. Therefore, for any  $(\mu;\eta_0)\in G_s\setminus G_C^+$, $(\mu;\eta_0)\not \in \coneCm$.
}

Theorem \ref{thm:subadditivegenerators} implicitly describes a way of obtaining all of the nontrivial extreme valid inequalities of $\coneC$: first identify a generating set $(G_L, G_s)$ for $\coneCs$ and then test its elements for $\cK$-minimality to identify $G_C^+$. 
On one hand, this is good news, as we seem to have a better algebraic handle on $\coneCs$ via the conditions given by \textbf{(A.0)}-\textbf{(A.2)}. On the other hand, testing these conditions as stated in \textbf{(A.0)}-\textbf{(A.2)}, is  a nontrivial task. Moreover, we need to establish further algebraic ways of characterizing $\cK$-minimality. Both of these tasks are tackled in the next section.

\section{Relations to Support Functions and Cut Generating Sets}\label{sec:support}

In this section, we first relate $\cK$-sublinear inequalities to the support functions of sets with certain structure. Recall that a \emph{support function} of a nonempty set $D \subseteq \bR^m$ is defined as 
\[
\sigma_D(z) := \sup_\lambda \{  z^T\lambda : ~ \lambda \in D \} \quad\mbox{for any } z\in \bR^m. 
\]

For any nonempty set $D$, it is well known that its support function, $\sigma_D(\cdot)$, satisfies the following properties:
\begin{itemize}
\item[(\textbf{S.1})] $\sigma_D(0)=0$, 
\item[(\textbf{S.2})] $\sigma_D(z^1+z^2) \leq \sigma_D(z^1) + \sigma_D(z^2)$,  (subadditive),
\item[(\textbf{S.3})] $\sigma_D(\beta z) =\beta \sigma_D(z) ~~ \forall \beta>0 \mbox{ and for all } z \in \bR^m$ (positively homogeneous).
\end{itemize}
In particular, support functions are positively homogeneous and subadditive, and thus, sublinear and convex. We refer the reader to  \cite{Hiriart-Urruty_Lemarechal_01,Rockafellar_70} for an extended exposure to the topic.

$\cK$-sublinear inequalities are closely related to support functions of convex sets with certain structure. This connection leads the way to a \emph{cut generating set} point of view as well as a number of necessary conditions for $\cK$-sublinearity. We state this connection in a series of results as follows:

\begin{theorem}\label{thm:supportbasic}
Consider any $\mu\in E$ satisfying condition (\textbf{A.0}), and define 
\begin{equation}\label{eqD}
D_\mu=\{\lambda \in\bR^m: ~A^*\lambda \preceq_{\cK^*} \mu \}.
\end{equation}
Then, $D_\mu\neq \emptyset$, $\sigma_{D_\mu}(0)=0$ and $\sigma_{D_\mu}(Az)\leq \langle \mu, z \rangle$ for all $z\in\cK$.
\end{theorem}
\myproof{ 
Since $\mu$ satisfies condition (\textbf{A.0}), we have $\mu\in\cK^*+\Im(A^*)$, which trivially implies the non-emptiness of $D_\mu$. 
Given that $\sigma_{D_\mu}(\cdot)$ is the support function of $D_\mu$ and $D_\mu\neq \emptyset$, we have $\sigma_{D_\mu}(0)=0$. 

Finally, for any $z\in\cK$, we have
\bse
\sigma_{D_\mu}(Az) &=& \sup_{\lambda} \{ \lambda^T Az: ~ \lambda \in D_\mu \} 
= \sup_{\lambda} \{ \langle z, A^*\lambda \rangle: ~ A^*\lambda \preceq_{\cK^*} \mu \} \\
&\leq&\sup_{\lambda} \{ \langle z, \mu \rangle: ~ A^*\lambda \preceq_{\cK^*} \mu \} = \langle z, \mu \rangle,
\ese
where the last inequality follows from the fact that $z\in\cK$ and for any $\lambda \in D_\mu$, we have $\mu-A^*\lambda\in\cK^*$, implying $\langle \mu - A^*\lambda, z \rangle \geq 0$. Therefore, $\sigma_{D_\mu}(Az)\leq \langle \mu , z \rangle$.
}

We note that every non-trivial valid linear inequality satisfies $\mu\in \textup{Im}(A^*)+\cK^*$ (see Proposition \ref{prop:necessityA0}). Thus,  any $\mu$ such that $\mu\not\in \textup{Im}(A^*)+\cK^*$ is redundant in the description of $\clconv(\solS)$. Furthermore, given a vector $\mu\in \textup{Im}(A^*)+\cK^*$, based on Theorem \ref{thm:supportbasic},  we can use the support function of the corresponding set $D_\mu$, and easily establish a condition on the right hand side value, $\eta_0$ that will ensure the validity of the inequality $(\mu;\eta_0)$. We state this result next.

\begin{proposition}\label{prop:cutgen}
Suppose  $\mu\in E$ satisfies condition (\textbf{A.0}). 
Then, $\inf_{b\in\cB} \sigma_{D_\mu}(b)\leq \vartheta(\mu)$, and thus, the inequality given by $(\mu;\eta_0)$ with $\eta_0\leq\inf_{b\in\cB} \sigma_{D_\mu}(b)$ is valid for $\solS$.
\end{proposition}
\myproof{
From the condition on $\mu$ and by Theorem \ref{thm:supportbasic}, we immediately have $D_\mu\neq \emptyset$ and $\sigma_{D_\mu}(Az)\leq \langle \mu , z \rangle$. Let $\widehat{\cB}:=\{b\in\cB:~\exists x \mbox{ s.t. } Ax=b,\, x\in\cK\}$. Then 
\bse
\eta_0 \leq \inf_{b\in\cB} \sigma_{D_\mu}(b) \leq \inf_{b\in\widehat{\cB}} \sigma_{D_\mu}(b) &=& \inf_{b\in\bR^m,x\in E}  \left\{\sigma_{D_\mu}(Ax):~ Ax=b,~b\in\widehat{\cB}  \right\} \\
&\leq& \inf_{x} \left\{\sigma_{D_\mu}(Ax):~ ~ x\in\cK, ~Ax\in\widehat{\cB} \right\} \\
& \leq & \inf_{x} \left\{ \langle \mu, x \rangle:~ x\in\cK, ~Ax\in\widehat{\cB} \right\} \\
& = & \inf_{x} \left\{ \langle \mu, x \rangle:~ x\in\cK, ~Ax\in\cB \right\} =\vartheta(\mu),
\ese
where the 
last inequality follows from the fact that for all $z\in\cK$, we have $\sigma_{D_\mu}(Az)\leq \langle \mu, z\rangle$, and the last two equations follow from $\widehat{\cB} \subseteq B$ and the definition of $\vartheta(\mu)$ (cf.\@ \eqref{mu0}). Then the validity of the inequality $(\mu;\eta_0)$  with $\eta_0 \leq \inf_{b\in\cB} \sigma_{D_\mu}(b)$ follows right away because $\inf_{b\in\cB} \sigma_{D_\mu}(b) \leq \vartheta(\mu)$. 
}

In addition to this, under a structural assumption on $\solS$,  we show that $\vartheta(\mu)$  for any $\mu\in E$ leading to a nontrivial inequality is exactly equal to $\inf_{b\in\cB} \sigma_{D_\mu}(b)$.

\begin{corollary}\label{cor:supportRHS}
Suppose $\Ker(A)\cap\intt(\cK)\neq\emptyset$. Then, for any  $\mu\in E$ satisfing condition (\textbf{A.0}), we have
 $\vartheta(\mu)=\inf_{b\in\cB} \sigma_{D_\mu}(b)$. 
\end{corollary}
\myproof{
By Proposition \ref{prop:cutgen}, we already have $\inf_{b\in\cB} \sigma_{D_\mu}(b) \leq \vartheta(\mu)$. Moreover, 
\bse
\inf_{b\in\cB} \sigma_{D_\mu}(b) &=& \inf_{b\in\cB} \sup_{\lambda\in\bR^m}\{b^T\lambda:~A^*\lambda \preceq_{\cK^*} \mu\}\\
&=& \inf_{b\in\cB} \underbrace{\inf_x \{ \langle \mu, x\rangle:~x\in\cK, ~Ax=b\}}_{\geq \vartheta(\mu)} \geq \vartheta(\mu), 
\ese
where the last equation follows from strong conic duality due to $\Ker(A)\cap\intt(\cK)\neq\emptyset$, and the last inequality follows from $b\in\cB$, and the definition of $\vartheta(\mu)$ in \eqref{mu0}. 
Thus, we obtain $\inf_{b\in\cB} \sigma_{D_\mu}(b) =\vartheta(\mu)$.
}

{
Given $\mu$ satisfying condition (\textbf{A.0}), there is a unique $D_\mu$ set associated with it, and Proposition \ref{prop:cutgen} highlights that one can use the support functions $\sigma_{D_\mu}(\cdot)$ of these sets $D_\mu$ to obtain a valid inequality based on $\mu$. 
Note that it is possible to have two distinct vectors $\mu'\neq\mu$ such that $D_\mu=D_{\mu'}$ (cf.\@ Example \ref{ex:SupportAtLeastOne}). 
These sets $D_\mu$ have a particular importance in our discussion in section \ref{sec:ConnectionToLatticeFreeSetsAndCGFs}. Due to the common structure of these sets $D_\mu$, we refer to the sets of this form as \emph{cut generating sets}.  
}

\subsection{Necessary Conditions for $\cK$-Sublinearity}

We next establish a number of necessary conditions for $\cK$-sublinearity via cut generating sets and their support functions.

\begin{proposition}\label{supportmain}
Suppose $\mu\in E$ satisfies condition (\textbf{A.0}). For any given $z\in\cK$, define 
\begin{equation}\label{eq:perp}
\perp_z := \{\gamma\in\cK^*: ~ \langle \gamma, z\rangle =0 \}.
\end{equation}
Then, for all $z\in\cK$ such that $\perp_z\cap(\mu-\textup{Im}(A^*))\neq \emptyset$, we have $\sigma_{D_\mu}(Az)= \langle \mu, z \rangle$ where $D_\mu$ is defined by \eqref{eqD}. 
\end{proposition}
\myproof{
Consider any $z\in\cK$, then we have
\bse
\sigma_{D_\mu}(Az) &=& \sup_{\lambda \in \bR^m} \{ \lambda^T Az : ~ \lambda \in D_\mu \} \\
&=& \sup_{\gamma \in E,\ \lambda \in \bR^m}\{ \langle z, A^*\lambda \rangle:  ~A^*\lambda = \mu-\gamma, ~ \gamma\in\cK^* \} \\
&=& \langle z, \mu \rangle - \inf_{\gamma \in E}\{ \langle z, \gamma \rangle:  ~\gamma\in \mu-\mbox{Im}(A^*), ~ \gamma\in\cK^* \} 
= \langle z, \mu \rangle
\ese
where the last equation follows from the fact that $ \langle z, \gamma \rangle\geq 0$ for all $z\in\cK$ and $\gamma\in\cK^*$, and there exists $\bar{\gamma}\in\perp_z\cap(\mu-\Im(A^*))$, i.e., $\bar{\gamma}\in\cK^* \cap (\mu -\textup{Im}(A^*))$ and $\langle \mu, \bar{\gamma}\rangle=0$.
}

Note that for $\mu\in\partial(\cK^*)+\Im(A^*)$, we immediately have $\partial(\cK^*)\cap(\mu-\Im(A^*))\neq\emptyset$, and thus, there exists $z\in\partial\cK$ such that $\perp_z\cap(\mu-\Im(A^*))\neq\emptyset$. In particular, for $\mu\in\Im(A^*)$, we have $0\in\cK^*\cap(\mu-\Im(A^*))$. Therefore, taking into account condition (\textbf{A.0}) and Theorem \ref{thm:supportbasic},  we have the following corollary of Proposition \ref{supportmain}: 
\begin{corollary}\label{cor:supportmain}
For any $\mu\in\partial(\cK^*)+\Im(A^*)$, we have $D_\mu\neq\emptyset$ and $\sigma_{D_\mu}(Az)= \langle \mu, z \rangle$ holds for at least one $z\in\Ext(\cK)$ where  $D_\mu$ is defined as in \eqref{eqD}. 
Moreover, for any $\mu\in\Im(A^*)$, we have $\sigma_{D_\mu}(Az)= \langle \mu, z \rangle$ for all $z\in\cK$. 
\end{corollary}

In the case of $\cK=\bR^n_+$,  using Remark \ref{remark:connectionToJohnson} the relationship between $\cK$-sublinearity and the support functions of cut generating sets can be further enhanced. 
\begin{proposition}\label{prop:RnSupportMain}
Consider a disjunctive conic set $\solS$ where $\cK=\bR^n_+$, and a $\cK$-sublinear inequality $(\mu;\eta_0)$ for it. 
Then, $\perp_{e^i}\cap(\mu-\textup{Im}(A^*))\neq \emptyset$, and thus, $\sigma_{D_\mu}(a^i)=\mu_i$ for all $i=1,\ldots,n$ where $a^i$ is the $i^{th}$ column of the matrix $A$.
Moreover, $\,\inf_{b\in\cB}\sigma_{D_\mu}(b)=\vartheta(\mu)$. 
\end{proposition}
\myproof{
Because $(\mu;\eta_0)$ is $\cK$-sublinear where $\cK=\bR^n_+$, $\mu\in E=\bR^n$ satisfies  conditions \textbf{(A.0)}-\textbf{(A.1i)} for all $i=1,\ldots,n$, and $\eta_0\leq \vartheta(\mu)$.  
Assume for contradiction that the statement is not true. Then there exist $i$ such that $\perp_{e^i}\cap(\mu-\textup{Im}(A^*))= \emptyset$. Note that $\perp_{e^i}=\{\gamma\in\bR^n_+:~ \gamma_i=0\}=\cone\{e^1,\ldots,e^{i-1},e^{i+1},\ldots,e^n\}$. Therefore, we arrive at the following system of linear inequalities in $\gamma,\lambda$ being infeasible:
\bse
\gamma+A^*\lambda = \mu, \\
\gamma_j \geq 0\ &&~~\forall j\neq i, \\
\gamma_i = 0.
\ese
Using Farkas' Lemma, we conclude that $\exists u,v$ such that $u+v=0$, $v_j\geq 0$ for all $j\neq i$, $Au=0$ and $\langle u, \mu\rangle \geq 1$. By eliminating $u$, this implies that $\exists v$ such that $v_j\geq 0$ for all $j\neq i$, $Av=0$ and $\langle v, \mu\rangle \leq -1$. Hence, if $v_i <-1$, we can scale $v$ so that $v_i \geq -1$, and arrive at the conclusion that there exists $v$ such that $v+e^i\in \bR^n_+=\cK$, $Av=0$ and $\langle v, \mu\rangle <0$, which is a contradiction to the condition \textbf{(A.1i)}. 

Because the conditions \textbf{(A.0)}-\textbf{(A.1i)} are necessary for the $\cK$-sublinearity (and also $\cK$-minimality) of $(\mu;\eta_0)$, using Proposition \ref{supportmain}, we conclude that for all $i=1,\ldots,n$, we have $\perp_{e^i}\cap(\mu-\textup{Im}(A^*)) \neq \emptyset$, and $\sigma_{D_\mu}(a^i)=\mu_i$. 

Finally, note that 
\bse
 \inf_{b\in\cB} \sigma_{D_\mu}(b) &=& \inf_{b\in\cB}  \sup_{\lambda\in\bR^m} \left\{b^T\lambda:~ A^*\lambda \leq \mu  \right\} \\
&=& \inf_{b\in\cB}  \inf_{x} \left\{ \mu^T x: ~Ax=b, ~ x\in\bR^n_+ \right\}  
=\vartheta(\mu),
\ese
where the second equation follows because condition \textbf{(A.0)} implies $\mu\in\bR^n_+ + \Im(A^*)$, and thus, the inner linear optimization problem is feasible, and so strong linear programming  duality holds. 
Thus, we have equality relations throughout implying $\inf_{b\in\cB} \sigma_{D_\mu}(b)=\vartheta(\mu)$. 
}

Proposition \ref{prop:RnSupportMain} has an important consequence that we point out next.
\begin{remark}\label{rem:SupportSummaryRnp}
Let $(\mu;\eta_0)$ be a $\bR^n_+$-sublinear inequality for $\cS(A,\bR^n_+,\cB)$. Given linear map $A\in\bR^{m\times n}$, let $a^i$ denote the $i^{th}$ column of $A$. Then Proposition \ref{prop:RnSupportMain} guarantees that for all $i=1,\ldots,n$ the support function $\sigma_{D_\mu}(\cdot)$ evaluated at the vector $a^i$, namely the data corresponding to variable $x_i$, precisely matches with the corresponding coefficient of $x_i$ in the inequality $(\mu;\eta_0)$,  i.e., $\mu_i=\sigma_{D_\mu}(a^i)$. Besides, $\sigma_{D_\mu}(\cdot)$ generates the tightest possible right hand side value, $\inf_{b\in\cB}\sigma_{D_\mu}(b)=\vartheta(\mu) \geq \eta_0$. Another way to state this is that \emph{every tight $\bR^n_+$-sublinear inequality (its coefficient vector, and the corresponding best possible right hand side value) is  generated by the support function $\sigma_{D_\mu}(\cdot)$}, a very specific sublinear function.
Furthermore, in the case of $\cK=\bR^n_+$,  the cut generating sets $D_\mu$ defined in \eqref{eqD} are polyhedral. Precisely, they are of the form   
\[
D_\mu=\{\lambda\in\bR^m:~A^*\lambda \leq \mu\}. 
\]
Thus, the support functions of these sets are automatically sublinear (subadditive and positively homogeneous), and in fact piecewise linear and convex. This relates nicely with the lattice-free sets, and cut generating functions.  
We discuss these in detail in section \ref{sec:ConnectionToLatticeFreeSetsAndCGFs}. 

Moreover, given any valid inequality $(\mu;\eta_0)$ for $\cS(A,\bR^n_+,\cB)$, if it is not $\bR^n_+$-sublinear, using the support function $\sigma_{D_\mu}(\cdot)$, one can immediately obtain an $\bR^n_+$-sublinear inequality dominating it (cf.\@ \cite[Proposition 2]{KKS14}.  
\epr
\end{remark}

Motivated by the positive result of Proposition \ref{prop:RnSupportMain} given  in the specific case of $\cK=\bR^n_+$, one is inclined to think that a similar result will hold for general regular cones $\cK$. We address this question in Proposition \ref{prop:SupportAtLeastOne}, and prove that in the case of general regular cones $\cK$,  for any $\cK$-sublinear inequality $(\mu;\eta_0)$, there exists at least one $z\in\Ext(\cK)$ such that $\sigma_{D_\mu}(Az)= \langle \mu, z \rangle$. 
Unfortunately, in the case of general regular cones $\cK$, the result of Proposition \ref{prop:SupportAtLeastOne} is not as strong as that of  Proposition \ref{prop:RnSupportMain}. 
Before we proceed with Proposition \ref{prop:SupportAtLeastOne}, we need a few technical lemmas. 

\begin{lemma}\label{lem:InfSwitch}
For any two sets $U$ and $V$ that are independent of each other, we have
\[\inf_{u\in U} \inf_{v\in V} \langle u, v \rangle = \inf_{v\in V} \inf_{u\in U} \langle u, v \rangle.\]
\end{lemma}
\myproof{
Let us  consider a given $\bar u\in U$. Then for any  $v\in V$, we have $\inf_{u\in U}\langle u, v\rangle \leq \langle \bar u, v\rangle$, and by taking the infimum of both sides of this inequality over $v\in V$, we obtain $\inf_{v\in V} \inf_{u\in U}\langle u, v\rangle \leq \inf_{v\in V} \langle \bar u, v\rangle$ holds for any $\bar u\in U$. Now, by taking the infimum of this inequality over $\bar u\in U$, and noting that the left hand side is simply a constant, we arrive at $\inf_{v\in V} \inf_{u\in U}\langle u, v\rangle \leq \inf_{\bar u \in U}\inf_{v\in V} \langle \bar u, v\rangle=\inf_{u \in U}\inf_{v\in V} \langle u, v\rangle$. To see that the reverse inequality also holds, we can start by considering a given $\bar v\in V$, and repeat the same reasoning by interchanging roles of $u$ and $v$.
}

\begin{lemma}\label{lem:SupportAtLeastOne}
Suppose that $\mu\in E$ satisfies condition (\textbf{A.0}), 
and $\perp_{z}\cap(\mu-\textup{Im}(A^*))= \emptyset$ holds for all $z\in\Ext(\cK)$ where $\perp_z$ is as defined by \eqref{eq:perp}. Then, there exists $\bar{\gamma}\in\intt(\cK^*)\cap (\mu-\textup{Im}(A^*))$. Moreover, $\inf_{b\in\cB}\sigma_{D_\mu}(b)=\vartheta(\mu)$. \end{lemma}
\myproof{
First, note that because $\mu$ satisfies condition (\textbf{A.0}), 
by Theorem \ref{thm:supportbasic}, $D_\mu\neq\emptyset$, which implies that $\{\gamma\in E: ~\exists\lambda\in\bR^m \mbox{ s.t. }  \gamma+A^*\lambda = \mu, ~ \gamma\in\cK^* \}\neq \emptyset$.
\par
In addition to this, $0 \in \bigcap_{z\in \textup{Ext}(\cK)} \perp_{z}$, and therefore,  together with the premise of the lemma that $\perp_{z}\cap(\mu-\textup{Im}(A^*))= \emptyset$, we conclude $0\not\in \mu-\textup{Im}(A^*)$. Moreover, by rephrasing the statement of lemma and definition of $\perp_z$, we get
\bse
0 &<& \inf_{z\in\textup{Ext}(\cK)} \inf_{\gamma\in E,\ \lambda \in \bR^m} \{\langle\gamma, z\rangle:~ \gamma+A^*\lambda = \mu, ~ \gamma\in\cK^* \} \\
&=& \inf_{\gamma\in E,\ \lambda \in \bR^m}\left\{ \inf_{z} \{\langle\gamma, z\rangle:~ z\in\Ext(\cK)\} :~ \gamma+A^*\lambda = \mu, ~ \gamma\in\cK^* \right\},
\ese
in which the last equation follows from Lemma \ref{lem:InfSwitch} where we take $U=\Ext(\cK)\times 0 \subseteq E \times \bR^m$ and $V=\{(\gamma,\lambda)\in E\times \bR^m:~  \gamma+A^*\lambda = \mu, ~ \gamma\in\cK^*\}$.
 
Now assume for contradiction that the set $\{\gamma: ~\exists\lambda\in\bR^m \mbox{ s.t. }  \gamma+A^*\lambda = \mu, ~ \gamma\in\cK^* \} \subseteq \partial\cK^*$. This together with the above inequality implies that there exists $\bar{\gamma}\in\partial\cK^*$ such that $\langle \bar{\gamma}, z\rangle >0$ for all $z\in\Ext(\cK)$. Hence, $\langle \bar{\gamma}, z\rangle >0$ for all $z\in\cK\setminus\{0\}$. Since $\cK^*$ is a closed convex cone, $\langle \bar{\gamma}, z\rangle >0$ for all $z\in\cK\setminus\{0\}$ implies that $\bar{\gamma}\in\intt(\cK^*)$, which is a contradiction.
Thus, we conclude that there exists $\bar{\gamma}\neq 0$ such that $\bar{\gamma}\in \intt(\cK^*)\cap(\mu-\textup{Im}(A^*))$.

To finish the proof note that 
\bse
 \vartheta(\mu) &:=& \inf_x \{ \langle \mu, x\rangle:~ x\in\solS \} \\
&=&\inf_{b\in\cB}\inf_x \{ \langle \mu, x\rangle:~ Ax = b, ~x\in\cK \}  \\
 &=&  \inf_{b\in\cB}\sup_{\lambda\in \bR^m,\gamma\in E} \{ b^T\lambda:~ A^*\lambda +\gamma = \mu, ~\gamma\in\cK^* \} = \inf_{b\in\cB}\sigma_{D_\mu}(b),
\ese
where the third equality follows from strong conic duality, which holds due to the existence of a strictly feasible solution $\bar{\gamma}\in\intt(\cK^*)$.  Therefore, we have $\vartheta(\mu)= \inf_{b\in\cB} \sigma_{D_\mu}(b)$.
}

We are now ready to state and prove Proposition \ref{prop:SupportAtLeastOne}.
\begin{proposition}\label{prop:SupportAtLeastOne}
Suppose that $\mu\in\viPi$, and $\perp_{z}\cap(\mu-\textup{Im}(A^*))= \emptyset$ holds for all $z\in\Ext(\cK)$ where $\perp_z$ is as defined by \eqref{eq:perp}. 
Then, there exists at least one $z\in\Ext(\cK)$ such that $\sigma_{D_\mu}(Az)= \langle \mu, z \rangle$. 
\end{proposition}
\myproof{
Assume for contradiction that $\sigma_{D_\mu}(Az)< \langle \mu, z \rangle$ for all $z\in\Ext(\cK)$. 
Then by Lemma \ref{lem:SupportAtLeastOne}, there exists $\bar{\gamma}\in\intt(\cK^*)\cap (\mu-\textup{Im}(A^*))$ and $\inf_{b\in\cB}\sigma_{D_\mu}(b) = \vartheta(\mu)$. 
Note that due to weak conic duality and the existence of such $\bar \gamma$,  we have for all $b$
\[
 \inf_{x} \{\langle \mu, x\rangle:~ Ax=b, ~x\in\cK\} \geq \sigma_{D_\mu}(b) = \sup_{\lambda\in\bR^m,\ \gamma\in E} \{ b^T \lambda:~ A^*\lambda+\gamma=\mu, ~\gamma\in\cK^*\} > -\infty .
\] 

For any $b\in\cB$, define $\cS_b:=\{x\in\cK:\, Ax=b\}$, and let $\widehat B:=\{b\in\cB:\, \cS_b\neq \emptyset\}$. Because $\solS\neq\emptyset$, $\widehat \cB\neq\emptyset$. Then for any $b\in\widehat \cB$, $x_b\in\cS_b$ leads to an upper bound on $\sigma_{D_\mu}(b)$, i.e., $\sigma_{D_\mu}(b)\leq \langle \mu, x_b \rangle$. Therefore, for any $b\in\widehat \cB$, the conic optimization problem defining $\sigma_{D_\mu}(b)$ is bounded above and is strictly feasible, and so we have strong conic duality and the dual problem given by the $\inf_x$ above is solvable. Consider any $b\in\widehat \cB$, let $\bar{x}_b$ be the corresponding optimal dual solution, i.e., $\bar{x}_b\in\cS_b$ and $\langle \mu, \bar{x}_b\rangle = \sigma_{D_\mu} (b)$.
Because $\bar{x}_b\in\cK$, there exists $z^1,\ldots,z^\ell\in\Ext(\cK)$ with $\ell\leq n$ such that $\bar{x}_b=\sum_{i=1}^\ell z^i$, which leads to
\bse
\langle \mu, \bar{x}_b \rangle = \sigma_{D_\mu} (b) = \sigma_{D_\mu}(A\bar{x}_b) \underbrace{\leq}_{(*)} \sum_{i=1}^\ell \sigma_{D_\mu}(Az^i) \underbrace{<}_{(**)} \sum_{i=1}^\ell \langle \mu, z^i \rangle = \langle \mu, \bar{x}_b \rangle,
\ese 
where the inequality $(*)$ follows because $\sigma_{D_\mu}(\cdot)$ is a support function, and thus is subadditive, and $(**)$ follows from the assumption that $\sigma_{D_\mu}(Az)< \langle \mu, z\rangle$ for all $z\in\Ext(\cK)$. But this is a contradiction. Thus, there exists  $z\in\Ext(\cK)$ such that $\sigma_{D_\mu}(Az)= \langle \mu, z\rangle$.
}

To summarize whenever $\mu\in\viPi$, Propositions \ref{supportmain} and \ref{prop:SupportAtLeastOne} together  cover all possible cases and indicate that for a $\cK$-sublinear inequality, there exists at least one $z\in\Ext(\cK)$ such that $\sigma_{D_\mu}(Az)=\langle \mu, z\rangle$. 

We illustrate the necessary conditions for $\cK$-sublinearity established so far via the following example. 

\begin{example}\label{ex:SupportAtLeastOne}
Consider the set $\solS$ with $\cK=\mathcal{L}^3$,  
$
~A=[1, 0, 0] \mbox{ and } \cB=\left\{-1,1\right\}.
$
In this case, $\Conv(\solS)=\{x\in\bR^3:~x\in\cK,~x_3 \geq \sqrt{1+x_2^2},~ -1\leq x_1\leq 1\}$ (see Figure \ref{fig:SupportAtLeastOne}). 
\begin{figure}[ht!]
\begin{center}
\begin{tikzpicture}[scale=1.2]
	\draw[dashed] (0,0) -- (1.945,1.88); 
	\draw[dashed] (0,0) -- (-1.945,1.88); 
	\draw[dashed] (0,2) ellipse (2 and 0.5);
	\draw[dashed,color=gray] (-1,1) arc (180:0:1 and 0.25);
	\draw[dashed] (-1,1) arc (-180:0:1 and 0.25);
	
	\fill[fill=blue!40!magenta!80!black, opacity=0.45] 
		(-1,1) parabola (-1.15,1.59) -- (-0.75,2.46)
		(-1,1) parabola (-0.75,2.46);
	\fill[fill=blue!40!magenta, opacity=0.9] 
		(1,1) -- (-1,1) parabola (-1.15,1.59) 
		 -- (0.85,1.55) 
		parabola (1,1);
		 
	\fill[fill=blue!40!white, opacity=0.85] 
		(-1.15,1.59)  --  (-0.75,2.46)  -- (1.25,2.39) -- (0.85,1.55) -- cycle;

	\fill[fill=blue!40!magenta!80!black, opacity=0.85] 
		(1,1) parabola (1.25,2.39) -- (0.85,1.55)
		(1,1) parabola (0.85,1.55);

	\draw[very thick,color=blue!70!black]  (1,1) parabola (1.25,2.39); 
	\draw[very thick,color=blue!70!black]  (1,1) parabola (0.85,1.55); 
	\draw[very thick,color=blue!70!black]  (0.85,1.55) -- (1.25,2.39); 
	\draw[very thick, dashed, color=blue!70!black,  opacity=0.5]  (-1,1) parabola (-0.75,2.46); 
	\draw[very thick,color=blue!70!black]  (-1,1) parabola (-1.15,1.59); 
	\draw[very thick,color=blue!70!black]  (-1.15,1.59) -- (-0.75,2.46); 
	\draw[very thick,color=blue!70!black] (-1,1) -- (1,1);	
	\draw[very thick,color=blue!70!black]  (-1.15,1.59) -- (0.85,1.55); 
	\draw[very thick,color=blue!70!black]  (-0.75,2.46) -- (1.25,2.39); 
	
	\draw[dashed] (0,2) ellipse (2 and 0.5);
	\draw[thick,->] ( 0, 0 ) -- ( 2.5, 0 ) node[ anchor = south east ] {$ x_1 $};
	\draw[thick,->] ( 0, 0 ) node[ anchor = north west ] {$(0,0)$} -- (-1.75,-1.25) node[ anchor = south east ] {$ x_2 $};
	\draw[thick,-] ( 0, 0 ) -- ( 0, 1 );
	\draw[thick,-,dotted, opacity=0.5] ( 0, 1 ) -- ( 0, 2 );
	\draw[thick,->] ( 0, 2 ) -- ( 0, 3.1 ) node[ anchor = north east ] {$ x_3 $};
	
\end{tikzpicture}

\end{center}
\caption{Convex hull of $\solS$ corresponding to Example \ref{ex:SupportAtLeastOne} }\label{fig:SupportAtLeastOne}
\end{figure}

Note that this description of the convex hull of $\solS$ involves the following inequalities:
\begin{itemize}
\item[(a)] $\mu^{(+)}=[1; 0 ;0]$ with $\eta^{(+)}_0=-1$ and $\mu^{(-)}=[-1; 0 ;0]$ with $\eta^{(-)}_0=-1$;
\item[(b)] $\mu^{(t)}=[0;t;\sqrt{t^2+1}]$ with $\eta^{(t)}_0=1$ for all $t\in\bR$.
\end{itemize}
Here, we show that these inequalities satisfy the necessary conditions for $\cK$-sublinear inequalities; later on we will in fact show that all of these inequalities are $\cK$-minimal.

In case (a), it is easily seen that the corresponding sets  associated with these inequalities $\mu^{(+)},\mu^{(-)}$ are given by
\bse
D_{\mu^{(+)}} &=& \{\lambda:~\exists\gamma\in\cK^*\mbox{ s.t. } \lambda+\gamma_1 = 1; ~\gamma_2 =0;  ~\gamma_3=0 \} =  \{\lambda:~ \lambda= 1\}, \\
D_{\mu^{(-)}} &=& \{\lambda:~ \lambda = -1\}.
\ese
Also, both $\mu^{(+)},\mu^{(-)}\in \textup{Im}(A^*)$, and thus, Corollary \ref{cor:supportmain} implies,  $\sigma_{D_{\mu^{(i)}}}(Az)=\sigma_{D_{\mu^{(i)}}}(z_1)=\langle \mu^{(i)}, z\rangle$ for all $z\in\cK$ for $i\in\{+,-\}$. In addition to this,  $\inf_{b\in \cB}\sigma_{D_{\mu^{(i)}}}(b) =-1= \eta^{(i)}_0$ for $i\in\{+,-\}$.

In case (b), for any given $t\in\bR$, we have the associated sets $D_{\mu^{(t)}}$ given by
\[
D_{\mu^{(t)}} = \{\lambda:~\exists\gamma\in\cK^*\mbox{ s.t. } \lambda+\gamma_1 = 0; ~\gamma_2 =t;  ~\gamma_3=\sqrt{t^2+1} \} =  \{\lambda:~ -1\leq \lambda\leq 1\}.
\]
Moreover, for  all $t$, by considering $z^t\in\{[1;-t;\sqrt{t^2+1}],[-1;-t;\sqrt{t^2+1}]\}\subset\Ext(\cK)$, we have $\langle \mu^{(t)}, z^t \rangle = 1$ and $\sigma_{D_{\mu^{(t)}}}(Az^t)=\sigma_{D_{\mu^{(t)}}}(z^t_1)=\sigma_{D_{\mu^{(t)}}}(1)=1$, proving $\langle \mu^{(t)}, z^t \rangle = \sigma_{D_{\mu^{(t)}}}(Az^t)$. Additionally, $\sigma_{D_{\mu^{(t)}}}(1)=1=\sigma_{D_{\mu^{(t)}}}(-1)$ implying $\inf_{b\in\cB}\sigma_{D_{\mu^{(t)}}}(b) = 1 = \eta^{(t)}_0$ for all $t$.

We highlight that $D_{\mu^{(t)}}$ is common for all distinct vectors $\mu^{(t)}$ corresponding to the valid inequalities $(\mu^{(t)};1)$. Nevertheless, each of these inequalities $(\mu^{(t)};1)$ are required for the description of $\clconv(\solS)$.

Let us also consider another valid inequality $(\nu;\nu_0)$ given by $\nu=[0;1;2]$ and $\nu_0=1$. Note that the associated $D_\nu$ set is given by
\[D_\nu=\left\{\lambda:~ -\sqrt{3}\leq \lambda \leq \sqrt{3}\right\}. 
\]
Furthermore, for any $z_\nu\in\left\{[{1\over \sqrt{3}};-{1\over 3}; {2\over 3}], [-{1\over \sqrt{3}};-{1\over 3}; {2\over 3}]\right\}\subset\Ext(\cK)$ we have $\sigma_{D_\nu}(Az_\nu)=\sigma_{D_\nu}(\pm{1\over\sqrt{3}}) = 1 =\langle\nu,z_\nu\rangle$. Also, $\inf_{b\in\cB}\sigma_{D_\nu}(b) = \sqrt{3} >1=\nu_0$. Therefore, in terms of the necessary conditions established so far  for $\cK$-sublinearity, there seems to be no difference between $(\nu;\nu_0)$ and the previous inequalities from above. When we revisit this example in the next section, we will show that while $(\nu;\nu_0)$ is $\cK$-sublinear,  $(\nu;\nu_0)$ is not $\cK$-minimal. 
In fact, we can easily show that $(\nu;\nu_0)$ is dominated by $\mu^{(1)}=[0;1;\sqrt{2}]$ and $\eta^{(1)}=1$. Because $\delta=\nu-\mu^{(1)}=[0;0;2-\sqrt{2}]\in\cK^*\setminus\{0\}$, we conclude that $(\nu;\nu_0)\not\in\coneCm$. 
\epr
\end{example}

\subsection{Sufficient Conditions for $\cK$-Sublinearity and $\cK$-Minimality}

Given any valid inequality $(\mu;\eta_0)$ satisfying condition (\textbf{A.0}), we can easily test $(\mu;\eta_0)$ for $\cK$-sublinearity with the help of the following proposition.
\begin{proposition}\label{prop:cutgenSubadditiveRel}
Let $(\mu;\eta_0)$ be such that $\mu$ satisfies condition (\textbf{A.0}) and $\eta_0\leq \inf_{b\in\cB}\sigma_{D_\mu}(b)$ (or it is known that $(\mu;\eta_0)\in\coneC$). Then, whenever there exists $x^i\in\Ext(\cK)$ such that $\sigma_{D_\mu}(Ax^i)=\langle \mu, x^i\rangle$ for all $i\in I$ and $\sum_{i\in I} x^i \in \intt(\cK)$, then the inequality $(\mu;\eta_0)$ is $\cK$-sublinear.
\end{proposition}
\myproof{
If we are given that $\eta_0\leq \inf_{b\in\cB}\sigma_{D_\mu}(b)$, then using Proposition \ref{prop:cutgen}, we have $(\mu;\eta_0)\in\coneC$, which automatically implies that condition (\textbf{A.2}) is satisfied. 
 
Next, given any $\alpha\in\Ext(\cK^*)$, we will verify condition (\textbf{A.1}$(\alpha)$). Consider any $u$ such that $Au=0$ and $\langle \alpha, v\rangle u+ v \in\cK ~\forall v \in \Ext(\cK)$.  Let $\cV_\alpha=\{v\in\Ext(\cK):~ \langle \alpha, v\rangle=1\}$, it is clear that $\langle u+v, \gamma \rangle \geq 0$ holds for all $v\in\cV_\alpha$ and $\gamma\in\cK^*$. Also, there exists $\bar{\lambda}$ and $\bar{\gamma}\in\cK^*$ satisfying $A^*\bar{\lambda}+\bar{\gamma}=\mu$ because $\mu$ satisfies condition (\textbf{A.0}), and hence, $\mu\in\cK^*+\Im(A^*)$. In fact, for any such $\bar{\lambda},\bar{\gamma}$, we have
\bse
\langle \mu, u\rangle &=& \langle A^*\bar{\lambda}+\bar{\gamma}, u\rangle 
= \langle \bar{\lambda}, \underbrace{Au}_{=0}\rangle + \langle\bar{\gamma}, u\rangle 
\,\geq\,  \langle \bar{\gamma}, -v \rangle~\quad\forall v \in \cV_\alpha .
\ese
Note that $\langle \gamma, -v \rangle\leq 0$ for all $\gamma \in \cK^*$ and $v \in \cV_\alpha\subset \cK$. 
In order to finish the proof, all we need to show is that there exists $\bar{v}\in\cV_\alpha$ such that $\langle \bar{\gamma}, \bar{v} \rangle=0$. Clearly, when $\mu\in\textup{Im}(A^*)$, we can take $\bar{\gamma}=0$, and hence conclude that $\langle \mu, u\rangle \geq -\langle \bar{\gamma}, \bar{v} \rangle=0$ holds for all such $u$. In the more general case, we have
\bse
&&\inf_{\gamma,\lambda} \left\{\inf_v\{\langle \gamma, v \rangle:~v\in\cV_\alpha\}:~A^*\lambda+\gamma=\mu, ~\gamma\in\cK^* \right\}\\
&=& \inf_v\left\{\inf_{\gamma,\lambda}\{\langle \mu-A^*\lambda, v \rangle: ~A^*\lambda+\gamma=\mu, ~\gamma\in\cK^*\}:~v\in\cV_\alpha  \right\}\\
&=& \inf_v\left\{\langle \mu, v \rangle - \underbrace{\sup_{\gamma,\lambda}\{ \lambda^T (Av) : ~A^*\lambda+\gamma=\mu, ~\gamma\in\cK^*\}}_{=\sigma_{D_\mu}(Av)}:~v\in\cV_\alpha  \right\}
\ese
Because there exists $x^i\in\Ext(\cK)$ such that $\sigma_{D_\mu}(Ax^i)=\langle \mu, x^i\rangle$ for all $i\in I$ and $\sum_{i\in I} x^i \in \intt(\cK)$, for any $\alpha\in\Ext(\cK^*)$, at least one of these $x^i$'s will be in $\cV_\alpha$. Otherwise, we have $\langle \alpha, x^i\rangle=0$ for all $i\in I$, and thus $\langle \alpha, \sum_{i\in I}x^i\rangle=0$, which is not possible since $\sum_{i\in I}x^i\in\intt(\cK)$ and $\alpha\in\Ext(\cK^*)$.  Thus, we conclude that the above infimum is zero. 
This gives us the desired conclusion that $\langle \mu, u\rangle\geq 0$, which proves that condition (\textbf{A.1}$(\alpha)$) is satisfied for any $\alpha\in\Ext(\cK^*)$. Hence, condition (\textbf{A.1}) is satisfied.
}

When $\cK=\bR^n_+$, Proposition \ref{prop:RnSupportMain} together with Theorem \ref{thm:supportbasic} implies that the conditions stated in Proposition \ref{prop:cutgenSubadditiveRel} are necessary and sufficient for $\cK$-sublinearity. 
For general regular cones $\cK$, based on the results from Theorem \ref{thm:supportbasic}, and Propositions \ref{supportmain} and \ref{prop:SupportAtLeastOne}, we conclude that the conditions stated in Proposition \ref{prop:cutgenSubadditiveRel} are almost necessary. This is up to the fact that for $\cK$-sublinear inequalitis $(\mu;\eta_0)$, we can prove the existence of at least one $x\in\Ext(\cK)$ satisfying $\sigma_{D_\mu}(Ax)=\langle \mu, x\rangle$, yet the sufficient condition in Proposition \ref{prop:cutgenSubadditiveRel} requires a number of such extreme rays summing up to an interior point of $\cK$. 
We next provide an example highlighting that  for general regular cones $\cK$ other than the nonnegative orthant, we cannot close this gap between the sufficient condition and the necessary conditions, i.e., there exists $\cK$-sublinear inequalities that satisfy only the necessary conditions from Propositions \ref{supportmain} and \ref{prop:SupportAtLeastOne} but not the sufficient condition of Proposition \ref{prop:cutgenSubadditiveRel}.

\begin{example}
Consider disjunctive conic set $\solS$ with $\cK=\mathcal{L}^3$,  
$
~A=[0, 1, 1] \mbox{ and } \cB=\left\{-1,1\right\}.
$
In this case, $\Conv(\solS)=\{x\in\cL^3:~ x_2+x_3= 1\}$.  
Let us examine the valid inequality $(\mu;\eta_0)$ given by $\mu=[0 ;0; 1]$ and $\eta_0=\vartheta(\mu)={1\over 2}$.  
Here, we first show that there is precisely a single ray $\bar{x}\in\Ext(\cK)$ such that $\sigma_{D_\mu}(A\bar{x})=\langle \mu, \bar{x} \rangle$, yet the inequality $(\mu;\eta_0)$ is a $\cK$-sublinear inequality.

The cut generating set associated with $\mu$ is $D_{\mu} = \{\lambda\in\bR:~ |\lambda | + \lambda \leq 1  \}$. Consider any $z\in\Ext(\cK)=\Ext(\cL^3)$, without loss of generality let us assume that $z$ is normalized to have $z_3=1$. Then 
\[
\langle \mu, z \rangle = \sigma_{D_\mu}(A z) 
\quad \Leftrightarrow \quad  z_3 = \sup_{\lambda\in\bR}\{ \underbrace{(z_2+z_3)}_{\geq 0 \text{ since } z\in \cL^3}\lambda:~ |\lambda | + \lambda \leq 1 \} 
\quad \Leftrightarrow \quad  z_3 = {1\over 2} (z_2+z_3).
\]
Therefore, $z_2=z_3=1$, and by noting $z\in\Ext(\cL^3)$, we get $z_1=0$.  Thus, we conclude that there is a unique extreme ray of $\cL^3$, in particular $z=[0;1;1]$  that satisfies $\langle \mu, z \rangle = \sigma_{D_\mu}(A z)$.

Let us now prove that $(\mu;\eta_0)$ is indeed $\cK$-sublinear. The conditions (\textbf{A.0}) and (\textbf{A.2}) are easily verified. In order to verify condition (\textbf{A.1}), we need to verify that for any $\alpha\in\Ext(\cK^*)$,
\[
0 \leq \langle \mu, u \rangle  \mbox{ for all } u\in E \mbox{ such that } Au=0  \mbox{ and }  \langle \alpha, v\rangle u+ v \in\cK ~~\forall v \in \Ext(\cK), 
\]
holds. Let $\alpha\in\Ext(\cK^*)$ be given. For any $v \in \Ext(\cK)$ if $\langle \alpha, v\rangle=0$, then we automatically have $\langle \alpha, v\rangle u+ v \in\cK$. And if $\langle \alpha, v\rangle\geq 0$, then we can normalize $v$ to assume that $\langle \alpha, v\rangle=1$. So, by defining $\cV_\alpha:=\{v\in\Ext(\cK):\, \langle \alpha, v \rangle =1\}$, we can state the above requirement as 
\[
0 \leq \langle \mu, u \rangle  \mbox{ for all } u\in E \mbox{ such that } Au=0  \mbox{ and }  u+ v \in\cK ~~\forall v \in \cV_\alpha, 
\]
which, in our particular case, becomes 
\[
0 \leq u_3  \mbox{ for all } u\in\bR^3 \mbox{ such that } u_3=-u_2  \mbox{ and }  u+ v \in\cL^3 ~~\forall v \in \cV_\alpha. 
\]
Now notice that $u_3=-u_2$, $u+ v \in\cL^3$ and $v\in\Ext(\cL^3)$ implies $u_3+v_3\geq 0$ and $u_3^2 + (v_1^2+v_2^2) + 2u_3v_3 \geq u_3^2 + v_2^2 - 2u_3v_2 + u_1^2 + v_1^2 + 2u_1v_1$, which is equivalent to $ 2u_3(v_2+v_3) \geq   u_1^2  + 2u_1v_1$. Now suppose that $\alpha_1=0$, then $\bar{v}=[{1\over \alpha_3};0;{1\over \alpha_3}]\in\cV_\alpha$ and $\tilde{v}=[{-1\over \alpha_3};0;{1\over \alpha_3}]\in\cV_\alpha$. In this case, using these particular $\bar{v}$ and $\tilde{v}$, we conclude $u_3 \geq   \max\left\{ {u_1^2  + 2u_1\bar{v}_1 \over 2(\bar{v}_2+\bar{v}_3)} ,  {u_1^2  + 2u_1\tilde{v}_1 \over 2(\tilde{v}_2+\tilde{v}_3)} \right\}={u_1^2 +2|u_1\bar{v}_1|\over 2\bar{v}_3}\geq0$. Moreover, when $\alpha_1\neq 0$, we have $\alpha_2+\alpha_3>0$ (since $\alpha\in\Ext(\cL^3)$) and considering $\hat{v}=\left[0;{1\over 2(\alpha_2+\alpha_3)};{1\over 2(\alpha_2+\alpha_3)}\right]\in\cV_\alpha$, we once again conclude that $u_3\geq0$.
Note that this is precisely what was needed to prove that $(\mu;\eta_0)$ is $\cK$-sublinear.
\epr
\end{example}

In addition to Proposition \ref{prop:cutgenSubadditiveRel}, under  \textbf{Assumption 1}, we can state a sufficient condition for $\cK$-minimality as follows:
\begin{proposition}\label{prop:subadditiveMinimalRel}
Suppose that \textbf{Assumption 1} holds and we are given a $\cK$-sublinear inequality $(\mu;\eta_0)$ satisfying $-\infty<\eta_0=\inf_{b\in\cB}\sigma_{D_\mu}(b)$. Let $\widehat{\cB}=\{b\in\cB:~ \sigma_{D_\mu}(b)\leq \eta_0 \}$. 
Then, if there exists $b^i\in\widehat{\cB}$ and $x^i\in\cK$ such that $\sum_i x^i \in\intt(\cK)$,  $Ax^i=b^i$ and $\langle \mu, x^i \rangle = \eta_0$, then $(\mu;\eta_0)$ is $\cK$-minimal.
\end{proposition}
\myproof{
Consider any  $(\mu;\eta_0)\in\coneCs$ satisfying $\eta_0=\inf_{b\in\cB}\sigma_{D_\mu}(b)$. Assume for contradiction that $(\mu;\eta_0)\not\in\coneCm$, i.e., $\exists \delta\in\cK^*\setminus \{0\}$ such that $(\mu-\delta;\eta_0)\in\coneC$.

Suppose the premise of the proposition holds for some $b^i\in\widehat{\cB}$ and $x^i\in\cK$ such that $\sum_i x^i \in\intt(\cK)$,  $Ax^i=b^i$ and $\langle \mu, x^i \rangle = \eta_0$. Note that for $\beta_i>0$ with $\sum_i\beta_i = 1$, we have $\bar{x}:=\sum_i\beta_i x^i \in\intt(\cK)$ and moreover, by definition, $\bar{x}\in\conv(\solS)$, and $\langle \mu, \bar{x} \rangle = \eta_0$. Because  any valid inequality for $\solS$, in particular $(\mu-\delta;\eta_0)$, is valid for $\conv(\solS)$ as well, we arrive at the contradiction 
\[
\eta_0  \leq  \langle \mu-\delta, \bar{x} \rangle < \eta_0,  
\]
where the last inequality follows from $\bar{x}\in\intt(\cK)$ and $\delta\in\cK^*\setminus\{0\}$ implying $\langle \delta, \bar{x} \rangle >0$ together with $\langle \mu, \bar{x} \rangle = \eta_0$.
}

Proposition \ref{prop:subadditiveMinimalRel}, in particular, states that a $\cK$-sublinear inequality is also $\cK$-minimal whenever the inequality is tight at a point at the intersection of $\intt(\cK)$ and $\conv(\solS)$. In the MILP case, this resembles a sufficient condition for an inequality to be facet defining. Nonetheless, our minimality notion in general is much weaker. In the MILP case, all of the facets are necessary and sufficient for the description of $\clconv(\solS)$; yet in general, one does not need all of the $\cK$-minimal inequalities, only a generating set for $\coneCm$ together along with the constraint $x\in\cK$ is needed. 

Moreover, an immediate implication of Proposition \ref{prop:cutgenSubadditiveRel} and Corollary \ref{cor:supportmain} is as follows:
\begin{corollary}\label{cor:ImAsublinearity}
For any $\mu\in\textup{Im}(A^*)$ and $\eta_0\leq\vartheta(\mu)$, the inequality  $(\mu;\eta_0)$ is $\cK$-sublinear. 
\end{corollary}

We have  already seen in Proposition \ref{prop:MinimalViRhs} that when $\ker(A)\cap\intt(K) \neq \emptyset$, then any $\mu\in\Im(A^*)$ and any $-\infty<\eta_0\leq\vartheta(\mu)$ leads to a $\cK$-minimal inequality $(\mu;\eta_0)$. Corollary \ref{cor:ImAsublinearity} complements this result by showing that valid inequalities $(\mu;\eta_0)$ with $\mu\in\Im(A^*)$ are always $\cK$-sublinear regardless of the requirement $\ker(A)\cap\intt(K) \neq \emptyset$. Indeed, when $\ker(A)\cap\intt(K) \neq \emptyset$, it is easy to see that the additional $\cK$-minimality requirements of Proposition \ref{prop:subadditiveMinimalRel} are trivially satisfied by $(\mu;\vartheta(\mu))$.

\begin{proposition}\label{prop:supportbasicrhs}
Let $(\mu;\eta_0)$ be a $\cK$-minimal inequality such that $\mu\in\intt(\cK^*)$. Then $ \eta_0=\vartheta(\mu)=\inf_{b\in\cB}\sigma_{D_\mu}(b)$.
\end{proposition}
\myproof{
Because $(\mu;\eta_0)\in\coneCm$ and $\mu\in\cK^*$, by Proposition \ref{prop:minVIrhs}, we have 
\[
\eta_0=\vartheta(\mu)=\inf_x\{\langle \mu, x\rangle:~ x\in\solS\}.
\]
Moreover, because $(\mu;\eta_0)$ is $\cK$-minimal, it is also $\cK$-sublinear, and therefore $D_\mu$ as defined in \eqref{eqD} is nonempty. Besides, by Proposition \ref{prop:cutgen}, $\mu\in\cK^*$ implies that $\inf_{b\in\cB}\sigma_{D_\mu}(b) \leq \vartheta(\mu)$.
Assume for contradiction that $\vartheta(\mu) > \inf_{b\in\cB}\sigma_{D_\mu}(b)$, which implies
\bse
\vartheta(\mu)&>& \inf_{b\in\cB} \sigma_{D_\mu}(b) =\inf_{b\in\cB}\sup_{\lambda,\gamma} \{ b^T \lambda:~ A^*\lambda +\gamma = \mu, ~\gamma\in\cK^* \} \\  
& = & \inf_{b\in\cB}\underbrace{\inf_x \{ \langle \mu, x\rangle:~ Ax = b, ~x\in\cK \} }_{\geq\vartheta(\mu)~\textup{since } b\in\cB} \\
&\geq & \vartheta(\mu), 
\ese
where the last equality follows from strong conic duality, which holds due to the  fact that $\mu\in\intt(\cK^*)$, and the last inequality follows from the definition of $\vartheta(\mu)$ and the fact that infimum is over $b\in\cB$. But, this is a contradiction. Therefore, $\eta_0=\vartheta(\mu)=\inf_{b\in\cB}\sigma_{D_\mu}(b)$.
}

To demonstrate the proper uses of Propositions \ref{prop:cutgenSubadditiveRel},  \ref{prop:subadditiveMinimalRel} and \ref{prop:supportbasicrhs}, let us return to our previous example.

\vspace{10pt}
\noindent\textbf{Example \ref{ex:SupportAtLeastOne}  (cont.)}
First note that the convex hull of $\solS$ is full dimensional. To see this, one can demonstrate the existence of $n+1$ affinely independent points from $\solS\subseteq \bR^n$ where $n=3$. Thus, there is no valid equation for $\solS$ implying that the lineality space of $\coneC$ is just  the zero vector. 
Moreover, $\hat{z}=[1;0;2]\in\intt(\cK)\cap\solS$ and hence \textbf{Assumption 1} is satisfied.

We claim that 
\begin{itemize}
\item[(a)] $\mu^{(+)}=[1; 0 ;0]$ with $\eta^{(+)}_0=-1$ and $\mu^{(-)}=[-1; 0 ;0]$ with $\eta^{(-)}_0=-1$;
\item[(b)] $\mu^{(t)}=[0;t;\sqrt{t^2+1}]$ with $\eta^{(t)}_0=1$ for all $t\in\bR$.
\end{itemize}
are all $\cK$-minimal inequalities. We have already seen that the associated sets $D_{\mu^{(i)}}$ are nonempty, $\inf_{b\in\cB} \sigma_{D_{\mu^{(i)}}}(b) = \eta^{(i)}_0$ holds and there are tight extreme points, i.e., $\sigma_{D_{\mu^{(i)}}}(Az^{(i)})=\langle \mu^{(i)}, z^{(i)} \rangle$ satisfying the requirement of Proposition \ref{prop:cutgenSubadditiveRel}, and hence, all of them are in $\coneCs$ by Proposition \ref{prop:cutgenSubadditiveRel}. Moreover, in case $(a)$, by considering the points $z^{(+)}=[1;0;2]\in\intt(\cK)\cap\solS$ and $z^{(-)}=[-1;0;2]\in\intt(\cK)\cap\solS$, we get $\langle \mu^{(i)}, z^{(i)} \rangle = \eta^{(i)}_0$ holds for all $i\in\{+,-\}$. Therefore, using Proposition \ref{prop:subadditiveMinimalRel}, we conclude that these inequalities are also $\cK$-minimal. In case $(b)$, for any $t\in\bR$, consider  $z^{(t)}_+=[1;-t;\sqrt{t^2+1}]\in\cK\cap\solS$ and $z^{(t)}_-=[-1;-t;\sqrt{t^2+1}]\in\cK\cap\solS$. Note that we have $\langle \mu^{(t)}, z^{(t)}_+ \rangle = \eta^{(t)}_0=\langle \mu^{(t)}, z^{(t)}_- \rangle $ for all $t\in\bR$, and hence $z^{(t)}:={1\over 2}( z^{(t)}_+  +  z^{(t)}_- )=[0;-t;\sqrt{t^2+1}]\in\intt(\cK)\cap\conv(\solS)$. Thus, by Proposition \ref{prop:subadditiveMinimalRel}, we conclude that $(\mu^{(t)};\eta^{(t)}_0)\in\coneCm$ for all $t\in\bR$. 

We proceed by showing that the system of infinitely many linear inequalities corresponding to $(\mu^{(t)};\eta^{(t)}_0)=([0;t;\sqrt{t^2+1}];1)$ for all $t\in\bR$ indeed has a compact conic representation as follows: For all $x\in\solS$, we have
\bse
&& 1\leq 0x_1 + t x_2 + \sqrt{t^2+1} x_3 ~~~~~\forall t\in\bR \\
&\iff& 1\leq \inf_{t} \{0x_1 + t x_2 + \sqrt{t^2+1} x_3: ~ t\in\bR\} \\
&\iff& 1\leq \inf_{t,\tau} \{ t x_2 + \tau x_3: ~ t\in\bR, ~\tau\geq \sqrt{t^2+1} \} \\
&\iff&1\leq \inf_{t,\tau} \{ t x_2 + \tau x_3: ~ t\in\bR, ~\tau\geq \sqrt{t^2+1} \} \\
&\iff&1\leq \inf_{t,\tau} \{ t x_2 + \tau x_3: ~ t\in\bR, ~(1;t;\tau)\in\mathcal{L}^3 \} \\
&\iff&1\leq \sup_{\alpha} \{-\alpha_1: ~ \alpha_2=x_2, ~ \alpha_3=x_3, ~[\alpha_1;\alpha_2;\alpha_3]\in\mathcal{L}^3 \} ~~~~~\mbox{ due to $(*)$}\\
&\iff& [-1; x_2; x_3]\in\mathcal{L}^3, 
\ese
where $(*)$ is due to the fact that the primal conic optimization problem is strictly feasible, and hence, strong duality applies here. Note that we have arrived at the constraint $x_3 \geq \sqrt{1+x_2^2}$, which a cylinder in $\bR^3$, hence a particular conic quadratic inequality $[1;x_2;x_3]\in\cL^3$.  The validity of $x_3 \geq \sqrt{1+x_2^2}$ for all $x\in\solS$ follows from its derivation. Moreover, this conic quadratic inequality exactly implies all of the $\cK$-minimal inequalities $(\mu^{(t)};\eta^{(t)}_0)$ for all $t\in\bR$. Thus, in this example, the constraint $x_3 \geq \sqrt{1+x_2^2}$ along with the constraint $x\in\cL^3$, completely describes $\conv(\solS)$.

Finally, recall that we have seen the valid inequality $(\nu;\nu_0)$ given by $\nu=[0;1;2]$ and $\nu_0=1$ has an associated $D_\nu$ set which is nonempty and there are tight extreme points, i.e., $\sigma_{D_\nu}(Az^{(i)})=\langle \nu, z^{(i)} \rangle$ satisfying the requirement of Proposition \ref{prop:cutgenSubadditiveRel} and $\nu_0=1<\sqrt{3}= \inf_{b\in\cB} \sigma_{D_\nu}(b)$,  hence by Proposition \ref{prop:cutgenSubadditiveRel} $(\nu;\nu_0)\in\coneCs$. While $\sigma_{D_\nu}(Az_\nu)=\langle\nu,z_\nu\rangle=\nu_0=1$ holds for any (and only) $z_\nu\in\left\{[{1\over \sqrt{3}};-{1\over 3}; {2\over 3}], [-{1\over \sqrt{3}};-{1\over 3}; {2\over 3}]\right\}\subset\Ext(\cK)$ and the mid point of these two points is in the interior of $\cK$, this mid point is not in $\conv(\solS)$, i.e.,  the sufficiency condition for $\cK$-minimality stated in Proposition \ref{prop:subadditiveMinimalRel} fails. In fact,   $\nu\in\intt(\cK^*)$  and $(\nu;\nu_0)$ fails the necessary condition for $\cK$-minimality given in Proposition \ref{prop:supportbasicrhs}, that is, $\inf_{b\in\cB}\sigma_{D_\nu}(b)=\sigma_{D_\nu}(1)=\sigma_{D_\nu}(-1)=\sqrt{3}>1=\nu_0$. Hence, we conclude that $(\nu;\nu_0)$ is not $\cK$-minimal.
\epr

This example also suggests a technique to derive closed form expressions for convex valid inequalities by grouping all of the tight $\cK$-minimal inequalities associated with the same cut generating set. This approach is further exploited in \cite{KKY14_ipco,KKY14} in analyzing specific disjunctive conic sets obtained from a two-term disjunction on a regular cone $\cK$. 
In particular, in \cite{KKY14_ipco,KKY14} a characterization of tight $\cK$-minimal inequalities for this specific disjunctive conic set is given, and in the case of $\cK=\cL^n$, using conic duality, it is shown that these tight $\cK$-minimal inequalities can be grouped appropriately leading to a class of convex inequalities.

\subsection{Connections to Lattice-free Sets and Cut Generating Functions}\label{sec:ConnectionToLatticeFreeSetsAndCGFs}

In this section, we relate our results to the existing literature on lattice-free sets and cut-generating functions in the case of $\cK=\bR^n_+$ and discuss some implications for general cones $\cK$.  

In the case of $\cK=\bR^n_+$, Proposition \ref{prop:RnSupportMain} and Remark \ref{rem:SupportSummaryRnp} together with the basic facts on support functions conjoin nicely with  the views based on cut generating functions and lattice-free sets. 
To summarize, we have shown that in the case of disjunctive conic sets $\cS(A,\bR^n_+,\cB)$, all tight  $\bR^n_+$-sublinear inequalities $(\mu;\vartheta(\mu))$ 
are generated by the support functions $\sigma_{D_\mu}(\cdot)$  of cut generating sets $D_\mu=\{\lambda\in\bR^m:~A^*\lambda \leq \mu\}$.
That is, $\sigma_{D_\mu}(\cdot)$  take as input $a^i$, the $i^{th}$ column of the linear map $A$, compute the corresponding cut coefficient  of the variable $x_i$, $\mu_i=\sigma_{D_\mu}(a^i)$ for all $i=1,\ldots,n$, and the best possible right hand side value $\vartheta(\mu)=\inf_{b\in\cB}\sigma_{D_\mu}(b)$. Note that  these  support functions are automatically sublinear (subadditive and positively homogeneous), and in fact piecewise linear and convex. 
Moreover, under \textbf{Assumption 1}, using the sufficiency of $\cK$-minimal inequalities (Proposition \ref{prop:minimalcharacterization}) and Theorem \ref{thm:subadditive}, we conclude that all non-cone-implied inequalities for disjunctive conic sets $\cS(A,\bR^n_+,\cB)$  are generated by piecewise-linear, subadditive, and convex functions. 
In addition to this, recently, it is shown in \cite[Proposition 4.1]{KKS14} that without making any assumptions such as \textbf{Assumption 1}, $\bR^n_+$-sublinear inequalities always exist, and along with the nonnegativity restrictions $x\in\bR^n_+$, they are always sufficient to describe $\clconv(\cS(A,\bR^n_+,\cB))$.

These observations on the structure and sufficiency of $\bR^n_+$-sublinear inequalities for $\cS(A,\bR^n_+,\cB)$ provide a simple and intuitive explanation of the well-known strong functional dual for MILPs, e.g., all cutting planes for MILPs are generated by nondecreasing subadditive convex functions (cf.\@ \cite{Nemhauser_Wolsey_88}).  

{
The literature on cutting plane theory for MILP is extensive, we refer the reader to the recent survey \cite{Cornuejols_07}. 
A particular stream of research initiated by Gomory and Johnson \cite{Gomory_Johnson_72_1, Gomory_Johnson_72_2} and followed up by Johnson \cite{Johnson_74} studies  an infinite relaxation of an MILP, i.e., the \emph{mixed-integer group problem} of \cite{Gomory_Johnson_72_1}, and introduces \emph{cut generating functions}, that is, functions $\psi:\bR^m\rightarrow\bR$ such that the inequality
\[
\sum_{i=1}^{n} \psi(a^i)x_i \geq 1
\]
holds for all feasible solutions $x\in\bR^n_+$ for any possible number of  variables $n$ and any choice of columns,  $a^i$, corresponding to these variables and a fixed set $S=\Z^m$ and a point $f\notin S$ (leading to $\cB=-f+S=-f+\Z^m$ in our context).\myfootnote{Note that the cut generating functions studied in these infinite models are independent of the problem data $a^i$, and thus, they work for all problem instances of arbitrary dimension $n$ and problem data $A$ but for a given set $\cB$. We refer the reader to the survey \cite{Conforti_Cornuejols_Zambelli_11} and references therein.}   
The interest in these infinite models originates from deriving cuts from multiple rows of a simplex tableau and its various relaxations that are  obtained by imposing further structural restrictions on the set $S$, and thus on $\cB$. These models have been  investigated extensively 
(see \cite{Conforti_Cornuejols_Zambelli_11} for a recent survey in this area). 
In this framework, \emph{extreme functions} and \emph{minimal functions} are used as convenient ways of creating a hierarchy of functions that are sufficient to generate all cuts.  A valid function $\psi$ is said to be \emph{extreme} if there are no two distinct valid functions $\psi_1$, $\psi_2$ such that $\psi={1\over 2} \psi_1 + {1\over 2} \psi_2$. Extreme functions are sufficient to generate all valid inequalities. Furthermore, all extreme functions are minimal. A valid function $\psi$ is  \emph{minimal} if there is no valid function $\psi'$ distinct from $\psi$ such that $\psi'\leq \psi$ (the inequality relation between functions is stated as a pointwise relation). 

This literature is closely connected to the \emph{$S$-free (lattice-free) cutting plane} theory for MILPs. An \emph{$S$-free} convex set is a convex set that does not contain any point from the given set $S$ in its interior. When $S=\Z^m$ an $S$-free set is called a \emph{lattice-free} set. Usually, one is interested in finding an $S$-free set to generate a valid inequality that cuts off a given point $f\not\in S$. Thus, one seeks an $S$-free convex set that contains $f$ in its interior. 
These results are particularly related to the \emph{intersection cuts} of Balas \cite{Balas_71,Balas_72}. 
In his seminal work \cite{Balas_71,Balas_72} Balas initiated the use of \emph{gauge functions of lattice-free sets} to generate cuts. 
This view continues to attract a lot of attention in the MILP context because the gauge functions have the advantage that they can be evaluated using simpler formulas in comparison to cut generating functions from Gomory-Johnson's infinite group
problem. 
Several papers in this literature \cite{Andersen_etal_07,Borozan_Cornuejols_09,Conforti_Cornuejols_Zambelli_10,Conforti_Cornuejols_Zambelli_11} establish an intimate connection between minimal functions and maximal (with respect to inclusion) $S$-free convex sets for various different models of $S$. 
For example, Borozan and Cornu\'ejols \cite{Borozan_Cornuejols_09} showed that minimal valid inequalities for the infinite relaxation with $\cB=-f+\Z^m$ correspond to maximal lattice-free convex sets, and thus, they arise from nonnegative, piecewise linear, positively homogeneous, convex functions. 
In many cases, e.g., when the sufficiency of nonnegative cut generating functions is known, it is known that every  minimal cut generating function $\psi(\cdot)$, the corresponding set $\{r\in\bR^m:~\psi(r)\leq 1\}$ is a maximal lattice-free set, and vice-versa.  
We refer the interested reader to 
\cite{Conforti_Cornuejols_etal_13,Conforti_Cornuejols_Zambelli_10,Conforti_Cornuejols_Zambelli_11} for further details and recent results. 
}

For finite dimensional problem instances $\cS(A,\bR^n_+,\cB)$  our study provides an alternative view  on the same topic based on support functions. We underline that the finite dimensional setup is indeed more relevant in obtaining strong cuts from the simplex tableau because it does not further relaxes the problem to an infinite model. Besides, it is well known that not all extreme inequalities in an infinite model remain extreme in the underlying finite dimensional model (cf.\@ \cite{Cornuejols_Margot_09}). 

{
Let us consider a  given $\bR^n_+$-sublinear inequality $(\mu;\eta_0)$.  Without loss of generality we assume that $\eta_0\in\{0,\pm1\}$.  First, note that the sets underlying gauge functions and support functions are nicely related via polarity.  
To observe this, let us consider the polar set of $D_\mu$ given by 
\[
D_\mu^o := \{r\in\bR^m:~ \lambda^T r\leq 1~\forall \lambda\in D_\mu \}.
\]
Clearly, $D_\mu^o$ is a closed convex set containing the origin, and the (Minkowski) gauge function of $D_\mu^o$,  $\gamma_{D_\mu^o}(\cdot)$, is given by
\[
\gamma_{D_\mu^o}(r)=\inf_t\{t>0:~r\in t\,D_\mu^o\}.
\]
Note that $\gamma_{D_\mu^o}(\cdot)$ is a nonnegative, closed, and sublinear function, and when $0\not\in\intt(D_\mu^o)$, $\gamma_{D_\mu^o}(\cdot)$ can take the value of $+\infty$. Moreover, by \cite[Theorem C.1.2.5]{Hiriart-Urruty_Lemarechal_01}, we have $D_\mu^o=\{r\in\bR^m:~\gamma_{D_\mu^o}(r) \leq 1\}$, i.e., the gauge function $\gamma_{D_\mu^o}(\cdot)$ represents the set $D_\mu^o$. For a given sublinear function there is  a unique set associated with it in this manner. However, there can be other sublinear functions $\psi(\cdot)$ representing the same set $D_\mu^o$, i.e., $D_\mu^o=\{r\in\bR^m:~\psi(r) \leq 1\}$. 
Because sublinear functions are positively homogeneous, for any sublinear function $\psi(\cdot)$ such that $D_\mu^o=\{r\in\bR^m:~\psi(r)\leq 1 \}$, we have $\gamma_{D_\mu^o}(r)=\psi(r)$ for every $r$ satisfying $\psi(r)>0$.  
In order to obtain strong valid inequalities, one is interested in the smallest possible such sublinear function $\psi(\cdot)$ representing $D_\mu^o$.   
It is also well-known \cite[Corollary C.3.2.5]{Hiriart-Urruty_Lemarechal_01} that whenever $Q$ is a closed convex set containing the origin, the support function of $Q$ is precisely the gauge function $\gamma_{Q^o}$. 
For any $\mu\in\viPi$, the set $D_\mu$ is always closed and convex, yet, we are not always guaranteed to have $0\in D_\mu$. That said, when $\mu\in\cK^*$, we always have $0\in D_\mu$. 
Furthermore, whenever $0\in D_\mu$, we conclude 
the support function of $D_\mu$ studied here is precisely the gauge function of the polar set $D_\mu^o$, that is $\sigma_{D_\mu}=\gamma_{D_\mu^o}$. Next, we make this connection more explicit and comment on when $D_\mu^o$ is  $\cB$-free.

Based on the given $\bR^n_+$-sublinear inequality $(\mu;\eta_0)$, let us also define the set 
\[
V_\mu:=\{r\in\bR^m:\, \sigma_{D_\mu}(r) \leq \eta_0\}.
\]
Note that $V_\mu$ is a closed convex set since $\sigma_{D_\mu}(\cdot)$ is a sublinear function. When $\cK=\bR^n_+$, Proposition \ref{prop:RnSupportMain} implies $\vartheta(\mu)=\inf_{b\in\cB}\sigma_{D_\mu}(b) \geq \eta_0$, and thus, $ \cB\cap\intt(V_\mu)=\emptyset$ (in fact, we have something slightly stronger, that is, the relative interior of $V_\mu$ does not contain any points from $\cB$). Also, whenever $\eta_0>0$, the inequality $(\mu;\eta_0)$ separates the origin from $\clconv(\cS(A,\bR^n_+,\cB))$, and $0\in V_\mu$. Let us for a moment focus on the case of $\eta_0>0$, and without loss of generality assume that $\eta_0=1$. For example, when $\mu\in\cK^*$, and $(\mu;\eta_0)$ is a $\cK$-minimal inequality, by Proposition \ref{prop:minVIrhs} without loss of generality we can assume $\eta_0=1$. Then under the assumption that $0\in D_\mu$, we immediately observe that $V_\mu=D_\mu^o$, and conclude that $D_\mu^o$, the polar of the set $D_\mu$, is a $\cB$-free set. Thus, we arrive at the following result:

\begin{proposition} Suppose $\cK=\bR^n_+$ and let $(\mu;\eta_0)$ with $0\in D_\mu$ and $\eta_0>0$ be an $\bR^n_+$-sublinear inequality for $\cS(A,\bR^n_+,\cB)$. Then, the support function  $\sigma_{D_\mu}(\cdot)$ of $D_\mu$ is exactly the gauge function of its polar $D_\mu^o$, i.e., $\sigma_{D_\mu}=\gamma_{D_\mu^o}$. Thus, $\sigma_{D_\mu}(\cdot)$ is nonnegative implying $\vartheta(\mu)=\inf_{b\in\cB}\sigma_{D_\mu}(b)\geq 0$ and also, $D_\mu^o$ is a $\cB$-free set.
\end{proposition}
}

Valid inequalities $(\mu;\eta_0)$ with $\eta_0>0$ for disjunctive sets of form $\cS(A,\bR^n_+,\cB)$ have attracted specific attention in the MILP literature. 
For example, when we fix the dimension $n$ in the framework of  \cite{Conforti_Cornuejols_etal_13}, the set of interest  is exactly a disjunctive conic set with $\cK=\bR^n_+$.  Specifically, in \cite{Conforti_Cornuejols_etal_13}, the authors consider disjunctive sets of form $\cS(A,\bR^n_+,\cB)$ with $\cK=\bR^n_+$ for an arbitrary dimension $n$ (and thus $A$ is also arbitrary), but under the additional assumption that $\cB$ is a given nonempty, closed set satisfying $0\not\in\cB$. In this framework, the main focus is on cuts $\mu^Tx\geq \eta_0$ that separate the origin from $\clconv(\cS(A,\bR^n_+,\cB))$, and the properties of {cut generating functions}, that is $\psi:\bR^m\rightarrow \bR$, which takes as input $a^i$, the data pertaining to the variable $x_i$, and maps it to the corresponding cut coefficient $\mu_i$. Starting from a dominance relation among such functions, \cite{Conforti_Cornuejols_etal_13} establishes a minimality notion for cut generating functions and studies various structural properties of minimal finite-valued cut generating functions and their relations with  $\cB$-free sets. 

Let us examine the connection between our results and those from \cite{Conforti_Cornuejols_etal_13} by assuming that the dimension $n$ is fixed in advance in \cite{Conforti_Cornuejols_etal_13}. Under the assumption $0\notin\cB$ of \cite{Conforti_Cornuejols_etal_13}, it is easily seen that $0\notin\cS(A,\bR^n_+,\cB)$ (see \cite[Lemma 2.1]{Conforti_Cornuejols_etal_13}), and therefore, without loss of generality we can assume that the cuts separating the origin from $\clconv(\cS(A,\bR^n_+,\cB))$ have the form $(\mu;1)$, i.e., their right hand side value $\eta_0$ is 1. When the dimension $n$ is fixed, the main set of interest in \cite{Conforti_Cornuejols_etal_13} is exactly  our set $\cS(A,\bR^n_+,\cB)$ and the corresponding cuts separating the origin are a subset of inequalities  from $\cC(A,\bR^n_+,\cB)$. Furthermore, because these cuts have positive right hand sides, they are non-cone-implied, and thus, are all $\bR^n_+$-sublinear  \cite[Proposition 4.1 and Corollary 4.2]{KKS14}. 
Hence, the support functions of the corresponding sets $D_\mu$ do have a direct relation with the corresponding cut generating functions of interest from \cite{Conforti_Cornuejols_etal_13}. 

We note that whenever the function used to generate a valid inequality is finite-valued everywhere, it can be used for any data matrix $A$.  This underlies the cut generating function point of view. On the other hand, 
the support functions $\sigma_{D_\mu}(\cdot)$ associated with $\bR^n_+$-sublinear inequalities are not always guaranteed to be finite-valued. This  
indicates a distinction between our results and the ones from \cite{Conforti_Cornuejols_etal_13}.   
We believe that it is not necessary to require a function to be finite-valued everywhere in order to use it to generate cuts for a given problem instance with data matrix $A$. In particular, the functions that are not finite-valued everywhere, such as the support functions we are considering here, can still be  meaningful and interesting in terms of generating valid inequalities.  
Furthermore, given a problem instance $A,\,\cB$ and $\cK=\bR^n_+$, under further assumptions on $A$ and $\cB$, it may be possible to obtain an appropriate, nonempty, bounded set $\emptyset\neq\tilde{D}\subseteq D_\mu$ ensuring $\inf_{b\in\cB}\sigma_{\tilde{D}}(b)\geq\eta_0$ and $\sigma_{\tilde{D}}(a^i)=\mu_i$ for all $i=1,\ldots,n$. That is, the support function of $\tilde{D}$ is finite-valued everywhere and generates the same inequality $(\mu;\eta_0)$. 
{
Thus, under further technical assumptions we can in addition  ensure the finite-valuedness of the support functions $\sigma_{D_\mu}(\cdot)$, and then, they will lead to valid inequalities for an arbitrary selection of the columns $a^i$. That is, they will indeed be cut generating functions for the given set $\cB$.
}
Let us for example consider Example 6.1 of \cite{Conforti_Cornuejols_etal_13}.

\begin{example}\label{ex:SuppFiniteness}
Suppose $A$ is the $2\times 2$ identity matrix, $\cB=\{[0;1]\}\cup\{\Z;-1\}$ and $\cK=\bR^2_+$, which leads to $\cS(A,\bR^2_+,\cB)=\conv(\cS(A,\bR^2_+,\cB))=\{[0;1]\}$. This particular disjunctive conic set violates our \textbf{Assumption 1}, and therefore, none of the valid inequalities is $\bR^2_+$-minimal. Nevertheless, existence of $\bR^2_+$-sublinear inequalities is not based on \textbf{Assumption 1} (see \cite[Proposition 4.1]{KKS14}). Indeed, we next show that the particular inequality $(\mu;\eta_0)=([-1;1];1)$ considered in \cite{Conforti_Cornuejols_etal_13} is $\bR^2_+$-sublinear. It is easy to see that the sufficiency conditions for $\cK$-sublinearity established in Proposition \ref{prop:cutgenSubadditiveRel} are satisfied for this inequality.  Actually, the corresponding $D_\mu=\{(\lambda_1;\lambda_2)\in\bR^2:\, \lambda_1\leq -1,\, \lambda_2\leq 1\}$, and $\sigma_{D_\mu}(A\ \!e^1)=\sigma_{D_\mu}([1;0])=-1=\mu_1=\mu^T e^1$ and  $\sigma_{D_\mu}(A\ \!e^2)=\sigma_{D_\mu}([0;1])=1=\mu_2=\mu^T e^2$ and clearly $e^1+e^2\in\intt(\bR^2_+)$. Furthermore, $\inf_{b\in\cB}\sigma_{D_\mu}(b)=1=\eta_0$, proving that $(\mu;\eta_0)=([-1;1];1)$ is a tight $\bR^2_+$-sublinear inequality for this particular $\conv(\cS(A,\bR^2_+,\cB))$. On the other hand, the support function corresponding to this inequality is not finite valued everywhere. As a matter of fact, when we try to bound $D_\mu$ to obtain $\tilde{D}\subseteq D_\mu$ and use $\tilde D$ to generate a valid inequality, then we cannot ensure $\sigma_{\tilde{D}}(A\ \!e^i)=\mu^T e^i =\mu_i$ for $i=1,2$, and $\vartheta(\mu)=\inf_{\cB}\sigma_{\tilde{D}}(b)=1$ simultaneously. 
\epr
\end{example}

It was conjectured in \cite{Conforti_Cornuejols_etal_13} and later on proved in  \cite{Cornuejols_Wolsey_Yildiz_13}, that in addition to their  earlier assumption $0\notin \cB$, if we further suppose the following ``containment" assumption, 
$\cone(\{a^1,\ldots,a^n\})\supseteq \cB$,  
we can ensure the existence of finite-valued cut generating functions corresponding to every extreme inequality separating the origin from $\cS(A,\bR^n_+,\cB)$. 
Furthermore, it is shown in \cite[Proposition 4.3]{KKS14} that in the same setup and under the same containment assumption of \cite{Conforti_Cornuejols_etal_13,Cornuejols_Wolsey_Yildiz_13}, one can ensure that the support functions associated with all $\bR^n_+$-sublinear inequalities are finite-valued. Actually, there is a natural duality relation between the support functions we study here and the value functions used in the sufficiency proof of cut generating functions in \cite{Cornuejols_Wolsey_Yildiz_13}.  
We finish our discussion by examining a slight variant of Example \ref{ex:SuppFiniteness} obtained from setting $\bar{\cB}=\{[0;1]\}\cup\{(\Z^-;-1)\}$. 
Note that in this variant we still have $\cS(A,\bR^2_+,\cB)=\cS(A,\bR^2_+,\bar{\cB})$, and  $\cS(A,\bR^2_+,\bar{\cB})$ still violates the containment assumption of  \cite{Conforti_Cornuejols_etal_13,Cornuejols_Wolsey_Yildiz_13}. Nevertheless, we can show that $(\mu;\eta_0)=([-1;1];1)$ is generated by a finite-valued  cut generating function. Indeed, one can easily check that  the support function of the set $\tilde{D}:=\{(\lambda_1;\lambda_2)\in\bR^2:\, \lambda_1= -1,\, -1\leq\lambda_2\leq 1\}$ obtained from bounding $D_\mu$ will do the job. This indicates the possibility for weakening the containment assumption of \cite{Conforti_Cornuejols_etal_13,Cornuejols_Wolsey_Yildiz_13}.
  
{
Next, we comment on the fact that $0\in D_\mu$ is not always guaranteed. While $0\in D_\mu$ for all $\mu\in\cK^*$, the other cases of $\mu\in\Im(A^*)+\cK^*$ are also of interest. In such cases, by taking the polar of $D_\mu^o$, we obtain $D_\mu^{oo}:=(D_\mu^o)^o$, a closed convex set containing the origin.
In addition to this, we always have $D_\mu\subseteq D_\mu^{oo}$ and so $\sigma_{D_\mu}(r) \leq \sigma_{D_\mu^{oo}}(r) = \gamma_{D_\mu^o}(r)$, where the last equation follows from \cite[Proposition C.3.2.4]{Hiriart-Urruty_Lemarechal_01}. In general $\sigma_{D_\mu}(\cdot)$  and $\gamma_{D_\mu^o}(\cdot)$ may differ quite significantly, i.e., a support function can take negative values while a gauge function cannot. 
To address this issue of generating negative coefficients in cuts, in \cite{basu2011convex} the following subset of the relative boundary of $D_\mu^{oo}$ was considered:
\[
\widehat{D}_\mu^{oo} := \{\lambda\in D_\mu^{oo}:~\exists r\in D_\mu^o ~\mbox{ s.t. }~ \lambda^T r = 1  \}.
\]
Under the assumption $0\in\intt(D_\mu^o)$ (which does not necessarily hold in our setup), it was shown in \cite{basu2011convex}  that among the sublinear functions $\psi(\cdot)$ satisfying $D_\mu^o=\{r\in\bR^m:~\psi(r) \leq 1\}$, we have the following relation $\sigma_{\widehat{D}_\mu^{oo}}(r) \leq  \psi(r) \leq \gamma_{D_\mu^o} (r)$. Note that $\sigma_{\widehat{D}_\mu}(r) \leq \sigma_{\widehat{D}_\mu^{oo}}(r)$ holds for all $r$. Studying the cases when we have $\sigma_{\widehat{D}_\mu}(r)=\sigma_{\widehat{D}_\mu^{oo}}(r) $ and $\sigma_{D_\mu}(r)=\sigma_{\widehat{D}_\mu}(r)$ with or without the assumption $0\in\intt(D_\mu^o)$ is of independent interest for understanding the minimality of these support functions $\sigma_{D_\mu}(\cdot)$.
}

\begin{remark}\label{remark:cutgenfunc}
In the case of $\cK=\bR^n_+$, as discussed above, there are strong  connections between $\cK$-sublinear inequalities, cut generating functions \cite{Conforti_Cornuejols_Zambelli_11}, and the strong functional dual for MILPs \cite{Nemhauser_Wolsey_88}. 

Moving forward, one may be interested in extending the definition of a cut generating function from MILPs to MICPs. However, the situation seems to be much more complex for general regular cones $\cK$ other than the nonnegative orthant. In the MILP context, one of the main properties of a cut generating function is that the function acts \emph{locally} on each variable. Namely, the cut generating function takes as input solely the data associated with an individual variable $x_i$, i.e., the corresponding column $a^i$, and based on this input, it generates the individual cut coefficient $\mu_i$ associated with $x_i$. Imposing such a local view on cut generating functions is acceptable in the case of the nonnegative orthant because such cut generating functions are sufficient in the case of $\cK=\bR^n_+$. This, we believe, is strongly correlated with the fact that the underlying cone $\cK=\bR^n_+$ is decomposable in terms of individual variables. However, for general regular cones $\cK$ imposing the same local view requirement on cut generating functions turns out to be problematic, especially when the cone $\cK$ encodes non-trivial dependences among variables. 

In particular, Example \ref{ex:SupportAtLeastOne} reveals an important fact in this discussion: Unlike the case with $\cK=\bR^n_+$,  unless we make further  structural assumptions, for general $\solS$ with a regular cone $\cK$, even when  the cone $\cK$ is as simple as $\cL^3$, there are extreme (and also tight, $\cK$-minimal) valid linear inequalities such that there is no function acting locally on individual variables that can  generate precisely the vector defining the extreme inequality.    
Specifically, in Example \ref{ex:SupportAtLeastOne}, the linear map is given by $A=[1, 0, 0]$, and the class of valid inequalities $(\mu^{(t)};\eta^{(t)}_0)=([0;t;\sqrt{t^2+1}];1)$ parametrized by $t\in\bR$ are all extreme, and thus, necessary in the description of $\clconv(\solS)$. 
If one considers cut generating functions of the form that take as input the individual columns of $A$ and output the corresponding cut coefficient, then no such function $\psi(\cdot)$ will precisely generate the vector defining the inequality $(\mu^{(t)};\eta^{(t)}_0)=([0;t;\sqrt{t^2+1}];1)$ for any $t\in\bR$. This is because such a function $\psi(\cdot)$ will  inevitably need to satisfy $t=\mu^{(t)}_2=\psi(a^2)=\psi(0)=\psi(a^3)=\mu^{(t)}_3=\sqrt{t^2+1}$, which is impossible.

Therefore, this example demonstrates that for regular cones $\cK$ other than the nonnegative orthant, if we were to straightforwardly extend the notion of cut generating functions based on a local view from the MILP literature and rely only on such functions, we may completely miss large classes of nontrivial extreme inequalities necessary for the description of $\clconv(\solS)$. On the other hand, it may be possible to introduce and study \emph{cut generating maps} $\Gamma(\cdot)$,  which take a \emph{global} view and consider the entire data $A$ and  generate the cut coefficient vector $\mu$ at once, i.e., $\mu=\Gamma(A)$. We leave the questions around such cut generating maps, such as their existence, structural properties, sufficiency, etc., for future work.

{
On a positive note, for specific MICPs of form \eqref{eq:MICP2} discussed in Example \ref{ex:MICPform2}, Moran et al.\@ \cite{Moran_Dey_Vielma_12} show that a strong functional dual exists under a technical condition. Existence of strong MICP dual for these specific MICPs is equivalent to the sufficiency of (indeed, very specific classes of) finite-valued functions that generate the cut coefficients of all cuts for these sets. In fact, these functions from  \cite{Moran_Dey_Vielma_12} indeed act locally on each individual variable, and thus, naturally extend the standard cut generating function framework used in the MILP literature to specific MICPs of form \eqref{eq:MICP2}. Thus, in spite of the fact that Moran et al.\@ \cite{Moran_Dey_Vielma_12} do not refer to these functions as cut generating functions, they are indeed so.   
However, we highlight that the natural disjunctive conic representation $\solS$ for the specific class of MICPs from \cite{Moran_Dey_Vielma_12} discussed in Example \ref{ex:MICPform2}  impose further  structure. In particular, the underlying cone $\cK$ in the resulting equivalent disjunctive conic form representation $\cS(A,\bR^{2n}_+,\cB)$ of MICP given in \eqref{eq:MICP2}  is simply $\bR^{2n}_+$. 
On the other hand, the cone involved in Example \ref{ex:SupportAtLeastOne}  is $\cL^3$. 
On a related note, we do not know of the existence of a similar strong functional MICP dual result for MICPs of form \eqref{eq:MICP1} discussed in Example \ref{ex:MICPform1}. Example \ref{ex:SupportAtLeastOne}  suggests that such a result is not likely.
}
\epr
\end{remark}

\subsection{Connections to Conic Mixed Integer Rounding Cuts}
We start with the following simple remark.
\begin{remark}\label{remark:L2}
In the simple case of the polyhedral cone $\cK=\mathcal{L}^2=\{x\in\bR^2:~x_2\geq |x_1|\}$, there are only two extreme rays $\alpha^{(1)}=[1;1]$ and $\alpha^{(2)}=[-1;1]$. These extreme rays are orthogonal to each other, and thus condition \textbf{(A.1)} reduces to
\bse
(\textbf{A.1(i)}) &&0 \leq \sum_{i=1}^n \mu(a^i)u_i  ~~\mbox{ for all } u \mbox{ such that } Au=0 \mbox{ and } u+ \alpha^{(i)} \in\mathcal{L}^2 \quad \mbox{for $i=1,2$, }
\ese
where $a^i$ denotes the $i^{th}$ column of $A$. 
Following the same reasoning as in Proposition \ref{prop:RnSupportMain}, one can easily deduce that for any $\cK$-minimal valid inequality $(\mu;\eta_0)$ and any extreme ray $z$ of $\cK=\mathcal{L}^2$, we have $\sigma_{D_\mu}(Az)=\langle \mu, z\rangle$.
\epr
\end{remark}

Using Proposition \ref{prop:RnSupportMain} and Remark \ref{remark:L2}, we are ready to analyze the conic mixed integer rounding cuts  introduced in  \cite{Atamturk_Narayanan_10} for the following simple mixed integer set \begin{equation}\label{CMIRset}
\cS_0:=\{(x,y,w,t)\in\Z\times\bR^3_+:~ |x+y-w-b|\leq t \}.
\end{equation}
In \cite{Atamturk_Narayanan_10}, it is shown that when $b=\lfloor b \rfloor+f$ with $f\in(0,1)$, the valid inequality given by
\begin{equation}\label{CMIRcut}
(1-2f)(x-\lfloor b \rfloor) + f \leq t+y+w,
\end{equation}
together with the original conic inequality in $\cS_0$ gives $\overline{\conv}(\cS_0)$.

Here we will prove that \eqref{CMIRcut} is in fact a $\cK$-minimal inequality. The first step in this analysis is to transform $\cS_0$ into our normal form as
\begin{equation}\label{CMIRsetN}
\cS:=\left\{(y,w,t,\gamma)\in\bR^3_+\times\mathcal{L}^2:~ \left[ \begin{array}{c} y-w \\ t \end{array}\right] -\gamma = \left[ \begin{array}{c} b-x \\ 0 \end{array}\right]  \right\},
\end{equation}
which leads to $\cK=\bR^3_+\times\mathcal{L}^2$, which is a closed convex pointed cone with nonempty interior, and 
\[
A=\left[ \begin{array}{rrrrr}  1 & -1 & 0 & -1 & 0 \\ 0 & 0 & 1 & 0 & -1  \end{array}\right]  
\mbox{ and }
 \cB=\left\{ \underbrace{\left[ \begin{array}{c} f \\0  \end{array}\right] }_{:=b^+_1},  \underbrace{\left[ \begin{array}{c} 1+f \\0  \end{array}\right]}_{:=b^+_2} , \ldots, \underbrace{\left[ \begin{array}{c} f-1 \\0  \end{array}\right] }_{:=b^-_1},  \underbrace{\left[ \begin{array}{c} f-2 \\0  \end{array}\right]}_{:=b^-_2} , \ldots \right\}. 
\]

Before we proceed first note that \textbf{Assumption 1} is satisfied, i.e., for any $\epsilon_1,\epsilon_2>0$, $(y; w; t; \gamma_1; \gamma_2) = (f+\epsilon_1; \epsilon_1; \epsilon_2; 0; \epsilon_2)\in\intt(K)$ and also in $\solS$, therefore $\cK$-minimal inequalities exist. However, $\solS$ is not full dimensional, $t-\gamma_2=0$ is a valid equation. 
The set $D_e$ corresponding to this valid equation is simply $D_e=\{(\lambda_1,\lambda_2):~\lambda_1=1,\lambda_2=0\}=\{(1,0)\}$. The point $\bar{z}$ defined in the rest of this example works for this valid equation as well. Thus, the valid equation $t-\gamma_2=0$ satisfies the necessary condition for $\cK$-minimality. 

We can use the results of section \ref{sec:support} to verify that the inequality \eqref{CMIRcut} satisfies $\cK$-minimality conditions. Using the first equation in \eqref{CMIRsetN}, we get $y-w-\gamma_1 = b-x$, which implies that $x-\lfloor b \rfloor = -y+w+\gamma_1 + f $. By substituting $x-\lfloor b \rfloor $ with $-y+w+\gamma_1 + f $, in \eqref{CMIRcut}, we can rewrite it in terms of the variables in our representation as follows:
\bse
&&(1-2f)(-y+w+\gamma_1 + f) + f \leq t+y+w \\
&& (2-2f) y + 2f w + t + (2f-1) \gamma_1 + 0 \gamma_2 \geq f (2-2f).
\ese
This means, $\eta_0= f (2-2f)$, $\mu_1= 2-2f$, $\mu_2=2f$, $\mu_3=1$, $\mu_4=2f-1$ and $\mu_5=0$ in our usual notation. The necessary conditions for $\cK$-sublinearity state that for $D_\mu$ given by \eqref{eqD}, we should have $D_\mu\neq \emptyset$, and $\sigma_{D_\mu}(Az)=\langle \mu, z\rangle$ for all $z\in\Ext(\cK)$ (since all of the extreme rays of $\cK$ are orthogonal to each other). 

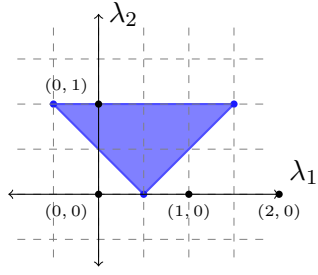
\begin{figure}[ht!]
\begin{center}
\begin{tikzpicture}[scale=1.2]
\fill[blue!50!white, thick] (-0.5,1) -- (1.5,1) -- (0.5,0) -- cycle;
\draw[blue!70!white, thick] (-0.5,1) -- (1.5,1) -- (0.5,0) -- cycle;
\fill[blue!90!white, thick] ( 0.5, 0 ) circle ( 0.04 );
\fill[blue!90!white, thick] ( -0.5, 1 ) circle ( 0.04 );
\fill[blue!90!white, thick] ( 1.5, 1 ) circle ( 0.04 );

\draw[ step = 0.5cm, gray, dashed, very thin] ( -0.9 ,-0.7 ) grid ( 1.9, 1.9 );
\draw[<->] ( 0, 2 ) node[ anchor = west ] {$ \lambda_2 $} -- ( 0, -0.8 ); 
\draw[<->] ( -1, 0 ) -- ( 2, 0 ) node[ anchor =  south west ] {$ \lambda_1 $};

\fill ( 0, 0 ) circle ( 0.04 );
\draw ( 0, 0 ) node[anchor = north east] {\tiny{$ ( 0, 0 ) $}};
\fill ( 0, 1 ) circle ( 0.04 );
\draw ( 0, 1 ) node[anchor = south east] {\tiny{$ ( 0, 1 ) $}};
\fill ( 1, 0 ) circle ( 0.04 );
\draw ( 1, 0 ) node[anchor = north] {\tiny{$ ( 1, 0) $}};
\fill ( 2, 0 ) circle ( 0.04 );
\draw ( 2, 0 ) node[anchor = north] {\tiny{$ ( 2, 0) $}};

\end{tikzpicture}
\end{center}
\caption{Feasible region corresponding to $D_\mu$ for $f=0.25$ in conic mixed integer rounding cut of \cite{Atamturk_Narayanan_10}.}
\label{fig:ex:CMIR}
\end{figure}

In our specific case, we have
\bse
D_\mu &=& \{\lambda \in\bR^2: ~ \exists \gamma \in \cK^*\mbox{ such that }  A^*\lambda +\gamma = \mu \} \\
& = & \left \{\lambda \in\bR^2: ~ \lambda_1 \leq \mu_1, ~ -\lambda_1\leq \mu_2, ~ \lambda_2 \leq \mu_3, ~  \left[ \begin{array}{c} -\lambda_1 \\ -\lambda_2  \end{array}\right]  \preceq_{\mathcal{L}^2}  \left[ \begin{array}{c} \mu_4 \\ \mu_5  \end{array}\right]  \right \} \\
& = &  \left \{\lambda \in\bR^2: ~ \lambda_1 \leq 2-2f, ~ -\lambda_1\leq 2f, ~ \lambda_2 \leq 1, ~ | 2f-1 +\lambda_1 | \leq \lambda_2  \right \}.
\ese
The set $D_\mu$ is plotted in Figure \ref{fig:ex:CMIR}.

Because $f\in(0,1)$, we have $D_\mu\neq\emptyset$, proving that $(\mu;\eta_0)$ is $\cK$-sublinear. Also the extreme rays of $\cK$ are precisely
$\Ext(\cK) = \{e^1, e^2, e^3, -e^4+e^5,  e^4+e^5\}$ where $e^i$ stands for the $i^{th}$ unit vector in $\bR^5$. Moreover, 
\[
\sigma_{D_\mu}(A e^1) = \sigma_{D_\mu}(a^1)= 2-2f = \mu_1 = \langle \mu, e^1\rangle,
\]
where $a^i$ denotes the $i^{th}$ column of the matrix $A$.  
Similarly we can show that $\sigma_{D_\mu}(Ae^i) = \mu_i = \langle \mu, e^i\rangle$ for $i=1,\ldots,3$. Moreover, we have $\sigma_{D_\mu}(A(-e^4+e^5)) =\sigma_{D_\mu}([1;-1]
) = 1- 2f = -\mu_4 +\mu_5
= \langle \mu, (-e^4+e^5)\rangle$ and $\sigma_{D_\mu}(A(e^4+e^5)) =\sigma_{D_\mu}([-1;-1]) = 2f-1 = \mu_4 +\mu_5 = \langle \mu, (e^4+e^5)\rangle$. 

Note that $\sigma_{D_\mu}(b^+_1)=f\cdot\sigma_{D_\mu}(e^1)=f(2-2f)$, and  for $i=1,2,\ldots$, we have $\sigma_{D_\mu}(b^+_{i+1})= (f+i) (2-2f)=(2-2f)i+2f-2f^2$.  Considering $f\in(0,1)$, we conclude $\sigma_{D_\mu}(b^+_1)<\sigma_{D_\mu}(b^+_2) < \ldots$ holds. Similarly $\sigma_{D_\mu}(b^-_1)=(1-f)\sigma_{D_\mu}(-e^1)=(1-f)(-2f)=2f(f-1)$, and  for $i=1,2,\ldots$, we have $\sigma_{D_\mu}(b^-_i)= (f-i) (-2f)=2fi-2f^2$, which implies   $\sigma_{D_\mu}(b^-_1)<\sigma_{D_\mu}(b^-_2)< \ldots$, and hence,
\[
\inf_{b\in\cB}\sigma_{D_\mu}(b)  = \min \left\{\sigma_{D_\mu}(b^+_1), ~\sigma_{D_\mu}(b^-_1) \right\}  = f (2-2f) = \eta_0.
\]

Finally, consider the following set of points
\[
\left\{z^1:=[f;0;0;0;0], z^2:=[0;1-f;0;0;0], z^3:=[0;0;f;-f;f], z^4:=[0;0;1-f;1-f;1-f]\right\}.
\]
Given $f\in(0,1)$, one can easily see that for $i=1,\ldots,4$, we have $z^i\in\solS$ and $\langle \mu, z^i \rangle=\eta_0=2f-2f^2$. Moreover, $\bar{z}:={1\over 4}\sum_{i=1}^4 z^i$ is in the interior of $\cK=\bR^3_+\times \mathcal{L}^2$. Therefore, using Proposition \ref{prop:subadditiveMinimalRel}, we have shown that the valid inequality given by $(\mu;\eta_0)=([2-2f;2f;1;2f-1;0];2f-2f^2)$, which is equivalent to \eqref{CMIRcut}, is a $\cK$-minimal inequality.

\section{Characterization of Valid Equations}\label{sec:equations}

Our results with regard to the existence of $\cK$-minimal inequalities was based on  \textbf{Assumption 1}, i.e., we assume that for all $\delta\in\cK^*\setminus\{0\}$, there exists $z_\delta\in\solS$ such that $\langle \delta, z_\delta \rangle > 0$. Under a stronger assumption, namely \textbf{Assumption 2} stated below, we can show that all valid equations $(\mu;\eta_0)$ satisfy $\mu\in\Im(A^*)$.  
\begin{quote}
\textbf{Assumption 2}: 
There exists $\hat{z}\in\solS$ such that $\hat{z}\in\intt(\cK)$ and $A\hat{z}=\hat{b}$ for some $\hat{b}\in\cB$. 
\end{quote}
Under \textbf{Assumption 2}, we can provide the following precise characterization of the valid equations.
 
\begin{theorem}\label{thm:ve}
Suppose that \textbf{Assumption 2} holds. 
Then $(\mu;\eta_0)$ is a valid equation if and only if there exists some $\bar{\lambda}\in\bR^m$ such that 
\[
A^*\bar{\lambda} = \mu \quad \mbox{ and }\quad   b^T\bar{\lambda}=\eta_0= \vartheta(\mu) \quad \mbox{for all } b\in\cB.
\]
\end{theorem}
\myproof{
\begin{itemize}
\item [$(\Leftarrow)$] It is easy to see that the condition in Theorem \ref{thm:ve} is sufficient. Suppose there exists $\bar{\lambda}\in\bR^m$ such that 
\[
A^*\bar{\lambda} = \mu ~\mbox{ and } ~ b^T\bar{\lambda} = \eta_0,
\]
for all $b\in\cB$. Then for any $z\in\solS$ we have
\[
\langle \mu, z\rangle = \langle A^*\bar{\lambda}, z \rangle = \bar{\lambda}^T Az = \bar{\lambda}^T b = \eta_0 = \vartheta(\mu),
\]
where the third equation follows because $z\in\solS$, and hence $Ax=b\in\cB$. This proves that $(\mu;\eta_0)$ is a valid equation. 

\item[$(\Rightarrow)$] To prove the necessity of the condition, suppose that $(\mu;\eta_0)$ is a valid equation. Then clearly $\eta_0=\vartheta(\mu)$. Let $\hat{b}$ and $\hat{z}$ be as described in \textbf{Assumption 2} preceding the theorem, and consider
\[
\inf_z\{\langle \mu, z\rangle:~ Az=\hat{b}, ~z\in\cK \}.
\]
This problem is strictly feasible because there exists $\hat{z}\in\intt(\cK)$ satisfying $A\hat{z}=\hat{b}$. Moreover, the solution set of this problem is contained in $\solS$. Thus, the fact that $(\mu;\vartheta(\mu))$ is a valid equation implies that its optimum value is equal to $\vartheta(\mu)$. By  strong conic duality we arrive at  
\[
\vartheta(\mu) = \sup_{\lambda \in \bR^m}\{\hat{b}^T \lambda: ~ A^*\lambda \preceq_{\cK^*} \mu \},
\]
which implies the existence of an optimal solution $\bar{\lambda}$ satisfying 
\[
A^*\bar{\lambda} \preceq_{\cK^*} \mu ~\mbox{ and }~ \hat{b}^T \bar{\lambda}=\vartheta(\mu).
\]
Note that any feasible solution to the primal problem is optimal including the  strictly feasible solution $\hat{z}$. Therefore, using the complementary slackness condition, we have
\[
\langle \hat{z}, \mu-A^*\bar{\lambda} \rangle = 0.
\]
Because $\hat{z}\in\intt(\cK)$, the above equation is possible if and only if $A^*\bar{\lambda} = \mu$.
Thus, we have established that there exists $\bar{\lambda}$ satisfying $A^*\bar{\lambda}=\mu$ and $\hat{b}^T\bar{\lambda} = \vartheta(\mu)$.
Now, for any $b\in\cB$, we have
\begin{equation}\label{mu+}
\vartheta(\mu)=\langle \mu, z_b\rangle 
~\geq~ \inf_z \{\langle \mu, z\rangle: ~ Az=b, ~z\in\cK\} 
~\geq~ \sup_{\lambda \in \bR^m}\{\hat{b}^T \lambda: ~ A^*\lambda \preceq_{\cK^*} \mu\} ~\geq \hat{b}^T \bar{\lambda}.
\end{equation}
Moreover, 
\begin{equation}\label{mu-}
-\vartheta(\mu)=\langle -\mu, z_b\rangle 
~\geq~ \inf_z \{\langle -\mu, z\rangle: ~ Az=b, ~z\in\cK\} 
~\geq~ \sup_{\lambda  \in \bR^m}\{\hat{b}^T \lambda: ~ A^*\lambda \preceq_{\cK^*} -\mu\} 
~\geq - \hat{b}^T \bar{\lambda},
\end{equation}
where the second inequality follows from weak duality and the last inequality follows because $-\bar{\lambda}$ is a feasible solution to the dual.
By combining \eqref{mu+} and \eqref{mu-}, we get $\vartheta(\mu) =\hat{b}^T \bar{\lambda}$,
which completes the proof.
\end{itemize}
}

Whenever \textbf{Assumption 2} holds, from Theorem \ref{thm:ve} we have $\mu\in\Im(A^*)$, and using the sufficient condition for $\cK$-minimality stated in Proposition \ref{prop:subadditiveMinimalRel}, we arrive at the following immediate corollary of Theorem \ref{thm:ve}.
\begin{corollary}
Suppose that \textbf{Assumption 2} holds. Then, any valid equation $(\mu; \vartheta(\mu))$ is $\cK$-minimal. 
\end{corollary}

In addition to the characterization of Theorem \ref{thm:ve}, we can relate each valid equation with its corresponding cut generating set $D_\mu$ given by \eqref{eqD} as follows:

\begin{corollary}\label{corr:ve}
Suppose that \textbf{Assumption 2} holds. Then,  for any valid equation $(\mu; \vartheta(\mu))$,  there exists $\lambda_\mu$ satisfying $D_\mu=\{\lambda: A^*\lambda \preceq_{\cK^*} \mu\} = \lambda_\mu+\{\lambda:~A^*\lambda \preceq_{\cK^*} 0\}$ and $\vartheta(\mu)=\inf_{b\in\cB} \sigma_{D_\mu}(b)=\sup_{b\in\cB} \sigma_{D_\mu}(b)$.
\end{corollary}
\myproof{
Suppose $(\mu; \vartheta(\mu))$ is a valid equation. Then by Theorem \ref{thm:ve} there exists $\bar{\lambda}=:\lambda_\mu$ such that $\mu=A^*\bar{\lambda}$ and $\vartheta(\mu)=b^T\bar{\lambda}$ for all $b\in\cB$. Thus, we have 
\[
D_\mu=\{\lambda: ~ A^*\lambda \preceq_{\cK^*} A^*\bar{\lambda} \} = \{\bar{\lambda}+\lambda: ~ A^*\lambda \preceq_{\cK^*} 0 \},
\]
and
\[
\inf_{b\in\cB} \sigma_{D_\mu}(b) = \inf_{b\in\cB} \sup_{\lambda \in \bR^m} \{ b^T (\bar{\lambda}+\lambda):~ A^*\lambda \preceq_{\cK^*} 0 \} 
= \inf_{b\in\cB} \left[ \underbrace{ b^T \bar{\lambda}}_{=\eta_0} + \underbrace{\sup_{\lambda} \{ b^T \lambda:~ A^*\lambda \preceq_{\cK^*} 0 \} }_{\in\{0,+\infty\}}\right] =\vartheta(\mu),
\]
where the last equation follows from the fact that $\vartheta(\mu)\in\bR$. Similarly, we can show that $\vartheta(\mu)=\sup_{b\in\cB} \sigma_{D_\mu}(b)$.
}

When $\cK=\bR^n_+$ (or any cone where each pair of its extreme rays is orthogonal), Corollary \ref{corr:ve} gives a complete characterization of valid equations.

\section{Conclusions and Further Research}
We introduce the class of $\cK$-minimal valid inequalities in the general disjunctive conic programming context and show that this class is a natural result of the dominance notion among valid inequalities, and thus,  contains a small yet essential set of nonredundant inequalities.  
In particular, under a mild technical assumption, we establish  that  the class of $\cK$-minimal inequalities together with the original constraint $x\in\cK$ are sufficient to describe $\clconv(\solS)$. This prompts an interest in $\cK$-minimal inequalities suggesting that an efficient cutting plane procedure for solving MICPs should at the least aim at separating inequalities from this class. Nevertheless, the definition of $\cK$-minimality reveals little about the structure of $\cK$-minimal inequalities.  In particular, testing $\cK$-minimality based on its definition is a non-trivial task. To address this issue, we show that the class of $\cK$-minimal inequalities is contained in a slightly larger class of so-called $\cK$-sublinear inequalities defined by algebraic conditions.  We establish a close connection between $\cK$-sublinear inequalities for disjunctive conic sets and the support functions of  convex sets with certain structure.  Using this connection, we show that when $\cK=\bR^n_+$, all $\cK$-sublinear inequalities are generated by sublinear (positively homogeneous, subadditive and convex) functions that are also piecewise linear. Thus, our  results naturally capture some of the earlier results from MILP setup, and generalize them to the conic case. Furthermore, this connection with support functions has led to practical ways of showing $\cK$-minimality  and/or $\cK$-sublinearity  properties of inequalities. To the best of our knowledge, these sufficient conditions for  $\cK$-minimality  and/or $\cK$-sublinearity of the valid inequalities are new even in the MILP setup.

Our work has shed some light on the structure of $\cK$-minimal and $\cK$-sublinear inequalities for disjunctive conic sets $\solS$ involving a regular cone $\cK$. However, many questions remain open when we start considering regular cones other than $\bR^n_+$. In particular, we find the following questions of interest:
\begin{itemize}
\item \textit{[Characterization of extreme valid inequalities]} Under a mild technical assumption, e.g., Assumption 1, we have shown that all extreme inequalities are $\cK$-minimal. However, not every $\cK$-minimal inequality is extreme (see e.g., Example \ref{ex:MinimalViRhs} and Proposition \ref{prop:MinimalViRhs}). Further characterizations of extreme inequalities beyond $\cK$-minimality are of great interest and importance. 
\item \textit{[Finiteness of the $\cK$-minimal conic inequalities]} When $\cK=\bR^n_+$ and $\cB$ is finite, Johnson \cite{Johnson_81} proved that the cone of $\cK$-minimal inequalities is finitely generated, i.e., $G_C$ is finite. Note that $G_L$ is always finite. 
For non-polyhedral regular cones, e.g., $\mathcal{L}^n,\cS^n_+$, in general, expecting $\clconv(\solS)$  to be given by finitely many linear inequalities is too much, and against the inherent nonlinear nature of these cones. Example \ref{ex:SupportAtLeastOne} shows that this is not possible even for $\mathcal{L}^3$, i.e., the resulting $\clconv(\solS)$ requires infinitely many  extreme linear inequalities. 
On the other hand, in that example, it is clear that the description of $\overline{\conv}(\solS)$ only involves two linear inequalities and two conic inequalities involving $\mathcal{L}^3$. While the $\cK$-minimality notion is seemingly  defined for linear inequalities, we can immediately extend it to a conic inequality by saying that a conic quadratic inequality is $\cK$-minimal if the associated (possibly infinite)  set of linear inequalities are  all $\cK$-minimal. 
We believe that instead of focusing on the finiteness of linear inequalities describing $\overline{\conv}(\solS)$, it is more natural and relevant to focus on the finiteness  of conic inequalities (of the same type of $\cK$) describing $\overline{\conv}(\solS)$. Therefore, we wonder what can be said in terms of the number of $\cK$-minimal conic inequalities required in the description of  $\overline{\conv}(\solS)$. Is it a finite number when $\cB$ is finite? Is it finite regardless of the size of $\cB$? Or, can we at least identify the cases where it is finite? In the very specific case of a two-term disjunctions on $\cL^n$, recent work of \cite{KKY14_ipco,KKY14} provide partial answers to some of these questions.

\item \textit{[Relations with valid inequalities for other nonconvex sets]} We showed that conic MIR inequalities introduced in \cite{Atamturk_Narayanan_10} can be interpreted in this framework. Moreover, in a recent series of papers \cite{KKY14_ipco,KKY14,YC14}, the characterization of tight $\cK$-minimal inequalities have played a critical role in the derivation of explicit expressions for convex valid inequalities for disjunctive conic sets associated with a two-term disjunction on $\cL^n$ and/or its cross-sections. These derivations relate back nicely to other recently developed valid inequalities for MICPs based on split or  disjunctive arguments in \cite{Andersen_Jensen_13,Belotti_Goez_Polik_Terlaky_12,Dadush_Dey_Vielma_11_2,Modaresi_Kilinc_Vielma_13}. 
Connecting our framework to other recent literature \cite{Bienstock_Michalka_13,BKK14,Modaresi_Kilinc_Vielma_13} covering more general setups involving nonconvex quadratic sets, and  extending our framework to cover these setups are also of interest.

\end{itemize}

\section*{Acknowledgements}
The author wishes to express her gratitude to the Associate Editor and anonymous referees for their constructive feedback, which led to substantial improvements on the presentation of the material in this paper. 
This work was first presented at UC Davis during Mixed Integer Programming Workshop in July 2012.

\bibliographystyle{abbrv} 
\bibliography{References,mypapers}

\end{document}